\documentclass[letterpaper]{article}

\usepackage{arxiv}

\usepackage{lineno,hyperref}

\usepackage{amsmath, amssymb, latexsym, multicol, amsthm, bm, placeins}
\usepackage[export]{adjustbox}
\usepackage{color}
\usepackage{multirow}

\usepackage{style}

\makeatletter\@addtoreset{equation}{section}\makeatother

\renewcommand{\shorttitle}{Enforcing Dirichlet boundary conditions in PINNs and VPINNs}

\title{Enforcing Dirichlet boundary conditions in physics-informed neural networks and variational physics-informed neural networks}
\author{S. Berrone\thanks{Dipartimento di Scienze Matematiche, Politecnico di Torino, Corso Duca degli Abruzzi 24, 10129 Torino, Italy. stefano.berrone@polito.it (S. Berrone), claudio.canuto@polito.it (C. Canuto), moreno.pintore@polito.it (M. Pintore).}
	\And
	C. Canuto\footnotemark[1]
	\And 
	M. Pintore\footnotemark[1]
	\And
	N. Sukumar\thanks{Department of Civil and Environmental Engineering, University of California, Davis, CA 95616, USA. nsukumar@ucdavis.edu}
}

\begin{document}

\maketitle

\begin{abstract}
In this paper, we present and compare four methods to enforce Dirichlet boundary conditions in Physics-Informed Neural Networks (PINNs) and Variational Physics-Informed Neural Networks (VPINNs). Such conditions are usually imposed by adding penalization terms in the loss function and properly choosing the corresponding scaling coefficients; however, in practice, this requires an expensive tuning phase. We show through several numerical tests that modifying the output of the neural network to exactly match the prescribed values leads to more efficient and accurate solvers. The best results are achieved by exactly enforcing the Dirichlet boundary conditions by means of an approximate distance function. We also show that variationally imposing the Dirichlet boundary conditions via Nitsche's method leads to suboptimal solvers.
\end{abstract}

\vspace{0.3cm}

\section{Introduction}
\label{intro_section}
Physics-Informed Neural Networks (PINNs), proposed in~\cite{raissi2019physics} after the initial pioneering contributions of Lagaris et al.~\cite{lagaris1997artificial, lagaris1998artificial, lagaris2000neural}, are rapidly emerging computational methods to solve partial differential equations (PDEs). In its basic formulation, a PINN is a neural network that is trained to minimize the PDE residual on a given set of collocation points in order to compute a corresponding approximate solution. In particular, the fact that the PDE solution is sought in a nonlinear space via a nonlinear optimizer distinguishes PINNs from classical computational methods. This provides PINNs flexibility, since the same code can be used to solve completely different problems by adapting the neural network loss function that is used in the training phase. Moreover, due to the intrinsic nonlinearity and the adaptive architecture of the neural network, PINNs can efficiently solve inverse~\cite{chen2020electromagnetic,guo2023high,mishra2021estimates}, parametric~\cite{gao2021phygeonet}, 
high-dimensional~\cite{han2018solving, LanthalerMishraKarniadakis2021} as well as nonlinear~\cite{jiang2021solving} problems. Another important feature characterizing PINNs is that it is possible to combine distinct types of information within the same loss function to readily modify the optimization process. This is useful, for instance, to effortlessly integrate (synthetic or experimental) external data into the training phase to obtain an approximate solution that is computed using both data and physics~\cite{chen2021physics}.

In order to improve the original PINN idea, several extensions have been developed. Some of these developments include the Deep Ritz method (DRM)~\cite{yu2018deep}, in which the energy functional of a variational problem is minimized; the conservative PINN (cPINN)~\cite{ameya2020conservative}, where the approximate solution is computed by a domain-decomposition approach enforcing flux conservation at the interfaces, as well as its improvement in the extended PINN (XPINN)~\cite{jagtap2020extended}; and the variational PINN 
(VPINN)~\cite{kharazmi2019variational,kharazmi2021hp}, in which the loss function is defined by exploiting the variational
structure of the underlying PDE.

Most of the existing PINN approaches enforce the essential (Dirichlet) boundary conditions by means of additional penalization terms that contribute to the loss function, these are each multiplied by constant weighting factors. See for instance~\cite{de2022error,de2022errorKolmogorov,demo2021extended, hu2022higher,pu2021solving,sirignano2018dgm,tartakovsky2018learning,yang2021bpinns,zhu2019physics}; note that this list is by no means exhaustive, therefore we also refer to \cite{beck2022overview,cuomo2022scientific,lawal2022physics} for more detailed overviews of the PINN literature. However, such an approach may lead to poor approximation, and therefore several techniques to improve it have been proposed. In~\cite{mcclenny2023self} and~\cite{wang2021understanding}, adaptive scaling parameters are proposed to balance the different terms in the loss functions. In particular, in~\cite{mcclenny2023self} the parameters are updated during the minimization to maximize the loss function via backpropagation, whereas in~\cite{wang2021understanding} a fixed learning rate annealing procedure is adopted. Other alternatives are related to adaptive sampling strategies (e.g.,~\cite{wight2021solving,tang2022adaptive,feng2021solving}) or to specific techniques such as the Neural Tangent Kernel~\cite{wang2022when}.

Note that although it is possible to automatically tune these scaling parameters during the training, such techniques require more involved implementations and in most cases lead to intrusive methods since the optimizer has to be modified. Instead, in this paper, we focus on three simple and non-intrusive approaches to impose Dirichlet boundary conditions and we compare their accuracy and efficiency. The proposed approaches are tested using standard PINN and 
interpolated VPINN which have been proven to be more stable than standard VPINNs~\cite{berrone2022variational}. 

The main contributions of this paper are as follows:
\begin{enumerate}
\item We present three non-standard approaches to enforce Dirichlet boundary conditions on PINNs and VPINNs, and discuss their mathematical formulation and their pros and cons. Two of them, based on the use of an approximate distance function, modify the output of the neural network to exactly impose such conditions, whereas the last one enforces them approximately by a weak formulation of the equation.
\item 
The performance of the distinct approaches to impose Dirichlet boundary conditions is assessed on various test cases. On average, we find that exactly imposing the boundary conditions leads to more efficient and accurate solvers. We also compare the interpolated VPINN to the standard PINN, 
and observe 
that the different approaches used to enforce the boundary conditions affect these 
two models in similar ways.
\end{enumerate}

The structure of the remainder of this paper is as follows. In 
Section~\ref{sec:pinn_description}, the PINN and VPINN formulations are described: first, we describe the neural network architecture in 
Section~\ref{sec:nn_description} and then focus on the loss functions that characterize the two models in Section~\ref{sec:loss_definition}. Subsequently, 
in Section~\ref{sec:used_approaches}, we present the four approaches to enforce the imposition of Dirichlet boundary conditions; three of them can be used with both PINNs and VPINNs, whereas the last one 
is used to 
enforce the required boundary conditions only on VPINNs because it relies on the variational formulation. Numerical results are presented in Section~\ref{sec:numerical_results}. In Section~\ref{sec:conv_rates}, we first analyze for a second-order elliptic problem the convergence rate of the VPINN with respect to
mesh refinement.
In doing so, we demonstrate that when the neural network is properly trained, identical optimal convergence rates are realized by all approaches only if the PDE solution is simple enough. Otherwise, only enforcing the Dirichlet boundary conditions with Nitsche's method or by exactly imposing them via approximate distance functions ensure the theoretical convergence rate. In addition, we compare the behavior of the loss function and the $H^1$ error while increasing the number of epochs, as well as the behavior of the error when the network architecture is varied. In Section~\ref{sect:nonlinear_parametric}, 
we show that it is also possible to efficiently solve second-order parametric nonlinear elliptic problems. 
Furthermore, in Sections~\ref{sect:elasticity}--~\ref{sect:transport}, we 
compare the performance of all approaches on PINNs and VPINNs by solving a linear elasticity problem and a stabilized Eikonal equation over an L-shaped domain, and a convection problem. Finally, in Section~\ref{sec:conclusion}, we close with our main findings and present a few perspectives for future work.

\section{PINNs and interpolated variational PINNs}\label{sec:pinn_description}
In this section, we describe the PINN and VPINN that are used in Section~\ref{sec:numerical_results}. In particular, in Section~\ref{sec:nn_description} the neural network 
architecture is presented, and the construction of the loss functions is discussed in Section~\ref{sec:loss_definition}.

\subsection{Neural network description}\label{sec:nn_description}
In this work we compare the efficiency of four approaches to enforce Dirichlet boundary conditions in 
PINN and VPINN. The main difference between these two numerical models is the training loss function; the architecture of the neural network is the same and is independent of the way
the boundary conditions are imposed. 

In our numerical experiments we only consider fully-connected feed forward neural networks with a fixed architecture. Such neural networks can be represented as nonlinear parametric functions $u^{\mathcal{N\!N}}: \mathbb{R}^{N_{\text{in}}} \rightarrow \mathbb{R}^{N_{\text{out}}}$ that can be evaluated via the following recursive formula:

\begin{equation}\label{eq:nn_rule}
\mathbf x_i^* = \sigma_i\left(A_i \mathbf x_{i-1}^* + b_i\right), \hspace{0.5cm} i = 1,2,\dots,L.
\end{equation}

In particular, with the notation of \eqref{eq:nn_rule}, $\mathbf x_0^*\in \mathbb{R}^{N_{\text{in}}}$ is the neural network input vector,  $\mathbf x_L^*\in \mathbb{R}^{N_{\text{out}}}$ is the neural network output vector, the neural network architecture consists of an input layer, $L-1$ hidden layers and one output layer, $A_i$ and $b_i$ are matrices and vectors containing the neural network weights, and $\sigma_i:\mathbb R\rightarrow\mathbb R$ is the activation function of the $i$-th layer and is element-wise applied to its input vector. We also remark that the $i$-th layer is said to contain $\text{dim}(\mathbf x_i^*)$ neurons and that $\sigma_i$ has to be nonlinear for any $i=1,2,\dots,L-1$. Common nonlinear activation functions are the rectified linear unit ($\textrm{ReLU}(x):=\text{max}(0,x)$), the hyperbolic tangent and the sigmoid function. In this work, we take $\sigma_L$ to be the identity function in order to avoid imposing any constraint on the neural network output. 

The weights contained in $A_i$ and $b_i$ can be logically reorganized in a single vector $\mathbf w^{\mathcal{N\!N}}$. The goal of the training phase is to find a vector $\mathbf w^{\mathcal{N\!N}}$ that minimizes the loss function; however, since such a loss function is nonlinear with respect to $\mathbf w^{\mathcal{N\!N}}$ and the corresponding manifold is extremely complicated, we can at best find good local minima.

\subsection{PINN and interpolated VPINN loss functions}\label{sec:loss_definition}
For the sake of simplicity, the loss function for PINN and interpolated VPINN is stated for second-order elliptic 
boundary-value problems.
However, the discussion can be directly generalized to different PDEs, and in Section \ref{sec:numerical_results}, numerical results associated with other problems are also presented.

Let us consider the model problem:
 \begin{equation}\label{eq:elliptic-model-pb}
\begin{cases}
Lu:=-\nabla \cdot (\mu \nabla u) + \boldsymbol{\beta}\cdot \nabla u + \sigma u =f & \text{in \ } \Omega , \\
u=g & \text{on \ } \Gamma_D , \\
\mu \frac{\partial u}{\partial n} = \psi  & \text{on \ } \Gamma_N , 
\end{cases}
\end{equation}
where $\Omega\subset\mathbb R^n$ is a bounded domain whose Lipschitz boundary $\partial\Omega$ is partitioned as $\partial\Omega=\Gamma_D\cup\Gamma_N$, with $\text{meas}_{n-1}(\Gamma_D)>0$. For the well-posedness of the boundary-value problem we require $\mu$, $\sigma\in \text{L}^\infty(\Omega)$ and $\boldsymbol \beta\in (\text{W}^{1,\infty}(\Omega))^n$ satisfying, in the entire domain $\Omega$, $\mu\ge\mu_0$ for some strictly positive constant $\mu_0$ and $\sigma-\frac{1}{2}\nabla\cdot\boldsymbol\beta\ge0$. Moreover, $f\in\text{L}^2(\Omega)$, $\psi\in\text{L}^2(\Gamma_N)$ and $g={\overline{u}}_{|\Gamma_D}$ for some ${\overline{u}}\in\text{H}^1(\Omega)$. We point out that even if these assumptions ensure the well-posedness of the problem, PINNs and VPINNs often struggle to compute low regularity solutions. We refer to \cite{taylor2023deep} for a recent example of a neural network based model that overcomes this issue.

In order to train a PINN, one introduces a set of collocation points $\{x_1,\dots,x_{N_I}\}$ and evaluates the corresponding equation residuals $\{r_1^\text{PINN},\dots,r_{N_I}^\text{PINN}\}$. Such residuals, for problem \eqref{eq:elliptic-model-pb}, are defined as:
\begin{equation}\label{eq:pinn_residuals}
r_i^\text{PINN}(u) = -\nabla \cdot (\mu \nabla u)(x_i) + \boldsymbol{\beta}\cdot \nabla u(x_i) + \sigma u(x_i) -f(x_i) 
\ \ \forall i=1,2,\dots,N_I.
\end{equation}
Since we are interested in a neural network that satisfies the PDE in a discrete sense, the loss function minimized during the PINN training is:
\begin{equation}\label{eq:pinn_loss}
R_{\text{PINN}}^2(w) = \sum_{i=1}^{N_I} \left\vert r_i^\text{PINN}(w) \right\vert^2.
\end{equation}
In~\eqref{eq:pinn_loss}, 
when $N_I$ is sufficiently large 
and $R_{\text{PINN}}^2(u^\NN) $ is close to zero, 
the function $u^\NN$ represented by the neural network output approximately satisfies the PDE and can thus be
considered a good approximation of the exact solution. Other terms are often added to impose the boundary conditions or improve the training, which are discussed in Section~\ref{sec:used_approaches}.

Let us now focus on the interpolated VPINN proposed in~\cite{berrone2022variational}. We introduce the function spaces $U:=H^1(\Omega)$ and $V:=\{v\in H^1(\Omega):v_{|\Gamma_D}=0\}$, the bilinear form $a:U\times V\rightarrow \R$ and the linear form $F:V\rightarrow \R$,
\begin{equation*}\label{eq:a_f}
a(w,v) = \int_\Omega \mu\nabla w\cdot\nabla v + \beta \nabla wv+\sigma wv, \hspace{1.0cm} F(v)=\int_\Omega fv + \int_{\Gamma_N}\psi v.
\end{equation*}
The variational counterpart of problem \eqref{eq:elliptic-model-pb} thus reads: Find $u\in U$ such that:
\begin{equation}\label{eq:var_probl_exact}
\begin{aligned}
a(u,v) &= F(v) \ \ \forall v\in V,\\
u &= g \ \ \text{on } \Gamma_D \,.
\end{aligned}
\end{equation}

In order to discretize problem \eqref{eq:var_probl_exact}, we use two discrete function spaces. Inspired by the Petrov-Galerkin framework, we denote the discrete trial space by $U_h\subset U$ and the discrete test space by $V_h\subset V$. The functions comprising such spaces are generated on two conforming, shape-regular and nested partitions $\mathcal T_H$ and $\mathcal T_h$ with compatible meshsizes $H$ and $h$, respectively. Assuming that $\mathcal T_h$ is the finer 
mesh, one can claim that $H\lesssim h < H$ and that every element of $\mathcal T_h$ is strictly contained in an element of $\mathcal T_H$.

\begin{figure}
\centering
\tabskip=0pt
\valign{#\cr
  \hbox{%
    \begin{subfigure}{.56\textwidth}
    \centering
    \includegraphics[width=\textwidth,keepaspectratio]{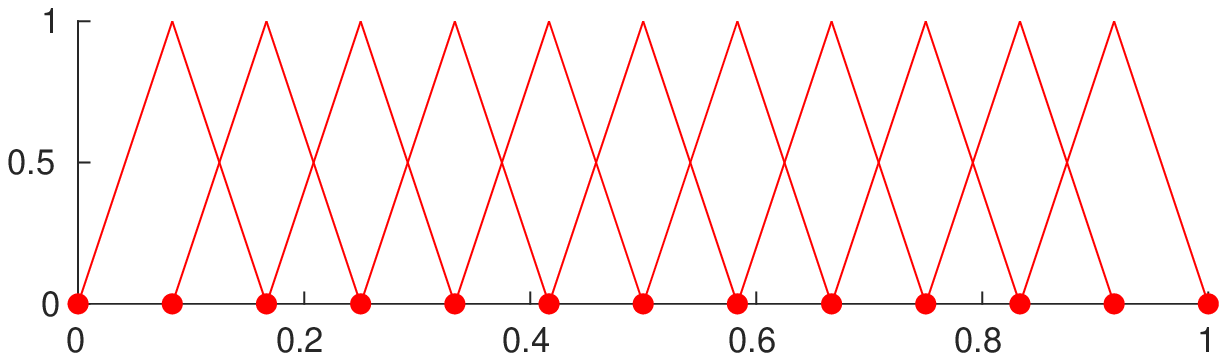} 
    \caption{}
    \label{fig:vh}
    \end{subfigure}%
  }\vfill
  \hbox{%
    \begin{subfigure}{.56\textwidth}
    \centering
    \includegraphics[width=\textwidth,keepaspectratio]{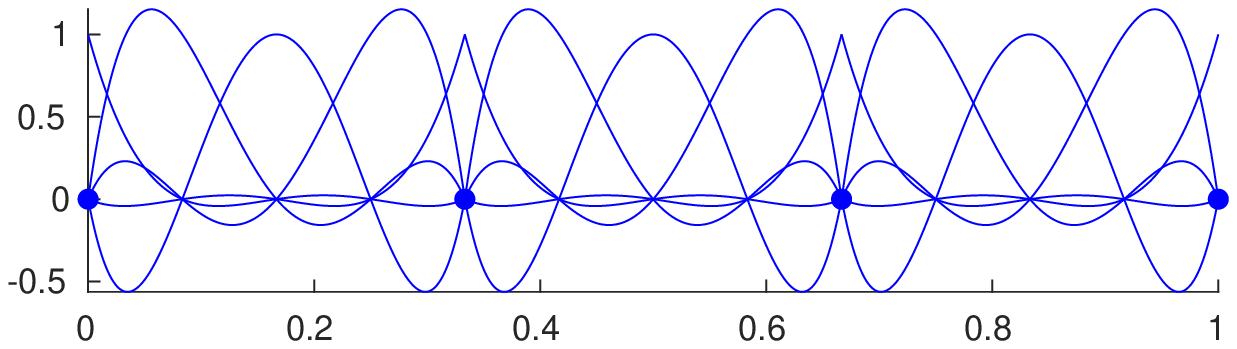}
    \caption{}
    \label{fig:uh}
    \end{subfigure}%
  }\cr
  \noalign{\hfill}
  \hbox{%
    \begin{subfigure}[t]{.4\textwidth}
    \centering
    \includegraphics[width=\textwidth]
    {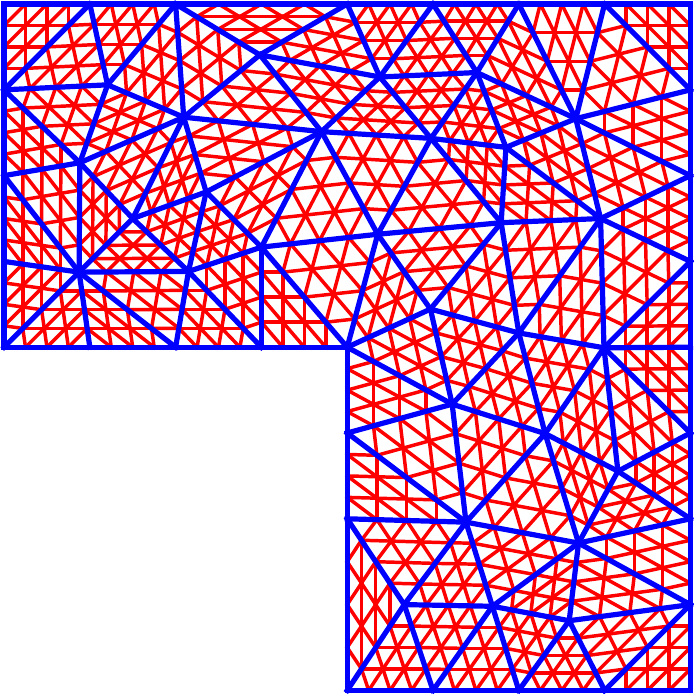} 
    \caption{}
    \label{fig:l-mesh}
    \end{subfigure}%
  }\cr
}
\caption{Pair of meshes and corresponding basis functions of a one-dimensional discretization (left) and nested meshes $\mathcal T_H$ and $\mathcal T_h$ in a two-dimensional domain (right). 
(a) Basis functions of $V_h$. The filled circles (red) are the nodes of the corresponding mesh $\mathcal T_h$;
(b) Basis functions of $U_H$. The filled circles (blue) are the vertex nodes that define the elements of the corresponding mesh $\mathcal T_H$;
and 
(c) Meshes used in the numerical experiments of Sections \ref{sect:elasticity} and \ref{sect:eikonal}. The blue mesh is $\mathcal T_H$, the red one is $\mathcal T_h$.
All the figures are obtained with $q=3$, $k_\text{test}=1$, $k_\text{int}=4$.}
\label{fig:discretizations}
\end{figure}

Denoting by $U_H:=\text{span}\{\varphi_i^u:i\in I_H\}\subset U$ the space of piecewise polynomial functions of order $k_{\text{int}}$ over $\mathcal T_H$ and $V_h:=\text{span}\{\varphi_i^v:i\in I_h\}\subset V$ the space of piecewise polynomial functions of order $k_{\text{test}}$ over $\mathcal T_h$ that vanish on $\Gamma_D$, we define the discrete variational problem as: Find $u\in U_H$ such that:
\begin{equation}\label{eq:global_discr_var_problem}
\begin{aligned}
a(u,v) &= F(v) \ \ \forall v \in 
V_h,\\
u &= g_H \ \ \text{on } \Gamma_D,
\end{aligned}
\end{equation}
where $g_H$ is a suitable piecewise polynomial approximation of $g$. 
A representation of the spaces $U_H$ and $V_h$ in a one-dimensional domain is provided in Figures \ref{fig:vh} and \ref{fig:uh}. Examples of pair of meshes $\mathcal T_H$ and $\mathcal T_h$ are shown in Fig.~\ref{fig:l-mesh}.

In order to obtain computable forms $a_h$ and $F_h$, we introduce elemental quadrature rules of order $q$ and define 
$a_h(\cdot,\cdot)$ and $F_h(\cdot)$ 
as the approximations of 
$a(\cdot,\cdot)$ and $F(\cdot)$ 
computed with such quadrature rules. In \cite{berrone2022variational}, under suitable assumptions, an a priori error estimate with respect to mesh refinement has been proved when $q=k_{\text{int}}+k_{\text{test}}-2$. It is then possible to define the computable variational residuals associated with the basis functions of $V_h$ as:
\begin{equation}\label{eq:residuals}
r_{h,i}(w) = F_h(\varphi_i^v) - a_h(w, \varphi_i^v), 
\quad i\in I_h.
\end{equation}
Consequently, in order to compute an approximate solution of problem \eqref{eq:global_discr_var_problem}, one seeks a function $w\in U_H$ that minimizes the quantity:
\begin{equation}\label{eq:var_loss}
R_h^2(w) = \sum_{i\in I_h}r_{h,i}^2(w),
\end{equation}
and satisfies the imposed boundary conditions. We refer to Section \ref{sec:used_approaches} for a detailed description of different approaches used to impose Dirichlet boundary conditions. It should be noted that, since in 
Sections~\ref{sect:nonlinear_parametric}--\ref{sect:transport} we consider problems other than~\eqref{eq:elliptic-model-pb}, the residuals 
in~\eqref{eq:residuals} have to be suitably modified, while the loss function structure defined in \eqref{eq:var_loss} is maintained.

We are interested in using a neural network to find the minimizer of $R_h^2$. We thus denote by $\mathcal I_H:C^0(\overline\Omega)\rightarrow U_H$ an interpolation operator used to map the function $u^{\mathcal N\!\mathcal N}$ associated with the neural network to its interpolating element in $U_H$, and train the neural network to minimize the quantity $R_h^2(\mathcal I_Hu^{\mathcal N\!\mathcal N})$. We highlight that in order to construct the function $\mathcal I_Hu^{\mathcal N\!\mathcal N}$, the neural network has to be evaluated only on $\text{dim}(U_H)$ interpolation points $\{x_1^{\mathcal I},\dots,x_{\text{dim}(U_H)}^{\mathcal I}\}\subset\overline\Omega$. Then, assuming that $\{\varphi_i^u:i\in I_H\}$ is a Lagrange basis such that $\varphi_i^u(x_j^{\mathcal I})=\delta_{ij}$ for every $i,j\in I_H$, it holds:
\begin{equation}\label{eq:interpolation_formula}
\mathcal I_Hu^{\mathcal N\!\mathcal N} = \sum_{i\in I_H} u^{\mathcal N\!\mathcal N}(x_i^{\mathcal I})\varphi_i^u.
\end{equation}

We remark that the approaches proposed in Section \ref{sec:used_approaches} can also be used on non-interpolated VPINNs. However, we restrict our analysis to interpolated VPINNs because of their better 
stability properties (see Fig.~\ref{fig:net_dims} and the corresponding discussion).

\section{Mathematical formulation} \label{sec:used_approaches}
We compare four methods to impose Dirichlet boundary conditions on PINNs and VPINNs. We do not consider Neumann or Robin boundary conditions since 
they can be weakly enforced by the trained VPINN due to the chosen variational formulation (computations using PINNs 
is discussed in~\cite{sukumar2022exact}). We also highlight that method $\mathbf{M_D}$ below can be used only with VPINNs because it relies on the variational formulation of the PDE. 
We analyze the following methods:
\begin{itemize}
\item[$\mathbf{M_A}$:] Incorporation of an additional cost in the loss function that penalizes unsatisfied boundary conditions; this is the standard approach in PINNs and VPINNs because of its simplicity and effectiveness. In fact, it is possible to choose $N_B$  control points $\{x_1^g,\dots,x_{N_B}^g\}\subset\Gamma_D$ and modify the loss functions defined in \eqref{eq:pinn_loss} or \eqref{eq:var_loss} as follows:
\begin{equation}\label{eq:ma_loss_pinn}
R_{\text{PINN}}^2(w) = \sum_{i=1}^{N_I} \left\vert r_i^\text{PINN}(w) \right\vert^2+ \lambda \sum_{i=1}^{N_B} \left( w(x_i^g)- g(x_i^g)\right)^2,
\end{equation}
or
\begin{equation}\label{eq:ma_loss}
R_h^2(w) = \sum_{i\in I_h}r_{h,i}^2(w) + \lambda \sum_{i=1}^{N_B} \left( w(x_i^g)- g(x_i^g)\right)^2,
\end{equation}
where $\lambda>0$ is a model hyperparameter. Note that on
considering the interpolated VPINN and exploiting the solution structure in~\eqref{eq:interpolation_formula}, it is possible to ensure the uniqueness of the numerical solution by choosing the control points $\{x_1^g,\dots,x_{N_B}^g\}$ as the $N_B$ interpolation points belonging to $\Gamma_D$.

We also highlight that such a method can be easily adapted to impose other types of boundary conditions just by adding suitable terms to \eqref{eq:ma_loss_pinn} and \eqref{eq:ma_loss}. On the other hand, despite its simplicity, the main drawback of this approach is that it leads to a more complex multi-objective optimization problem.

\item[$\mathbf{M_B}$:] Exactly imposing the Dirichlet boundary conditions as described in \cite{sukumar2022exact} and \cite{berrone2022variational}. In this method we add a non-trainable layer $B$ at the end of the neural network to modify its output $w$ according to the rule:
\begin{equation}\label{eq:b_def}
Bw = \overline g + \phi w,
\end{equation}
where $\overline g\in C^0(\overline \Omega)$ is an extension of the function $g$ inside the domain $\Omega$ (i.e., $\overline g_{|\Gamma_D}=g$) and $\phi\in C^0(\overline\Omega)$ is an approximate distance function (ADF) to the boundary $\Gamma_D$, 
i.e., $\phi(\bm{x})=0$ if and only if $\bm x\in\Gamma_D$, and it is positive elsewhere. During the training phase one minimizes the quantity $R_{\text{PINN}}^2(Bw)$ or $R_h^2(Bw)$.

For the sake of simplicity, we only consider ADFs for two-dimensional unions of segments, even though the approach generalizes to more complex geometries. Following the derivation of $\overline g$ and $\phi$ in \cite{sukumar2022exact}, we start by defining $d$ as the signed distance function from $\bm{x}:=(x,y)$ to the line defined by the segment $AB$ of length $L$ with vertices $A=(x_A, y_A)$ and $B=(x_B,y_B)$:
\[
d(\bm{x}) = \dfrac{(x-x_A)(y_B-y_A) - (y-y_A)(x_B-x_A)}{L}.
\]
Then, we denote $(x_c,y_c) := \bigl( (x_A+x_B)/2, (y_A+y_B)/2 \bigr)$ to be the center of $AB$ and define $t$ 
as the following trimming function:
\[
t(\bm{x}) = \frac 1L \left[\left(\frac L2\right)^2 - \left\Vert (x,y) - (x_c,y_c)\right\Vert^2\right].
\]
Note that $t\ge0$ defines a circle of center $(x_c,y_c)$. Finally, the ADF to $AB$ is defined as
\[
\phi(\bm{x}) = \sqrt{d^2 + \left(\frac{\sqrt{t^2+d^4}-t}{2}\right)^2}.
\]
A graphical representation of $d(\bm{x})$, $t(\bm{x})$ and $\phi(\bm{x})$ for an inclined line segment is shown 
in Figures~\ref{fig:dtphia},~\ref{fig:dtphib} and\ref{fig:dtphic}, respectively.

\begin{figure}[t!]
\centering 
\begin{subfigure}[t]{0.32\linewidth}
  \includegraphics[width=0.98\columnwidth,keepaspectratio,clip]{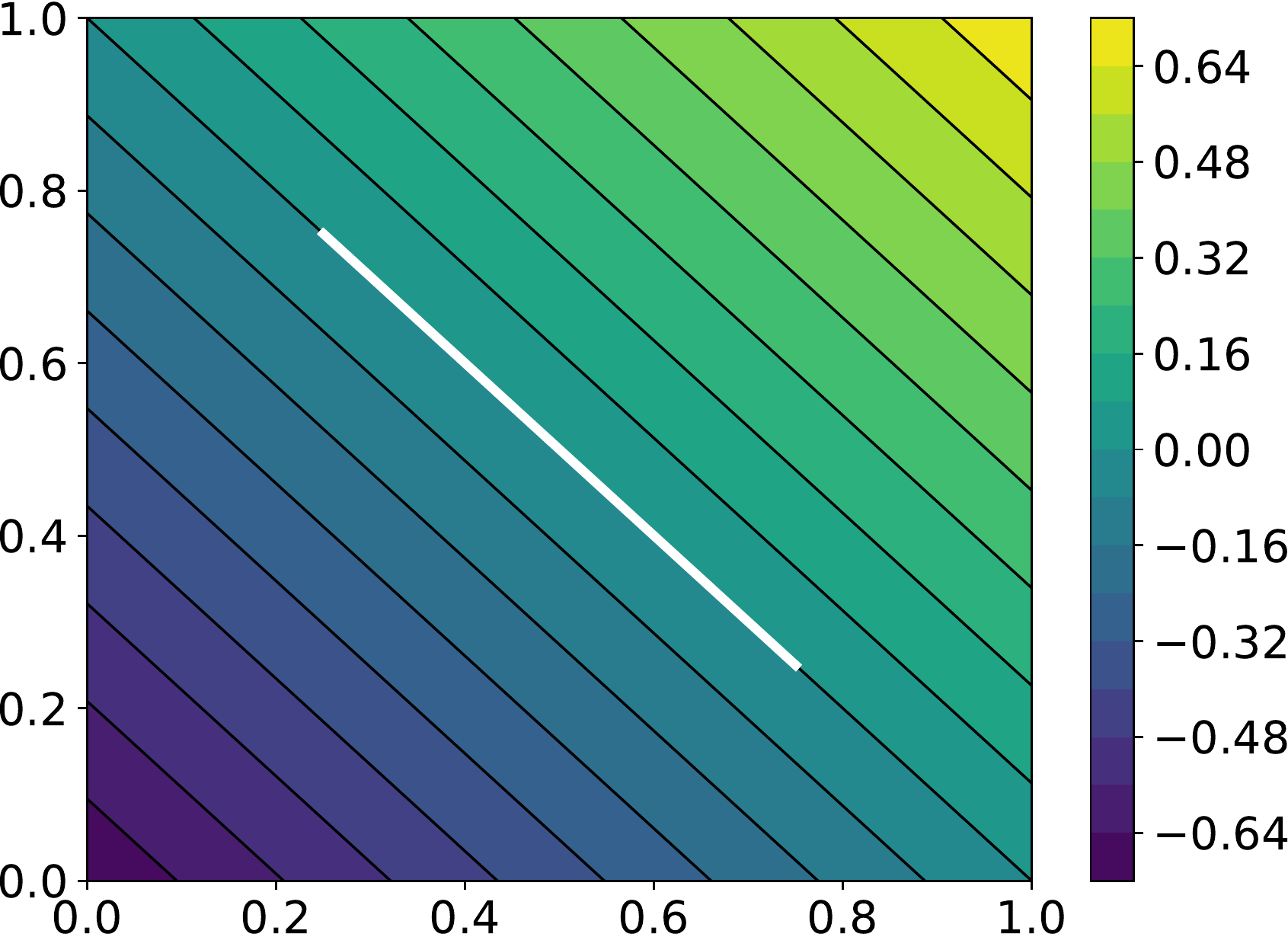} 
   \subcaption{$d(\mathbf x)$.}
  \label{fig:dtphia}
\end{subfigure}
\begin{subfigure}[t]{0.32\linewidth}
  \includegraphics[width=0.98\columnwidth,keepaspectratio,clip]{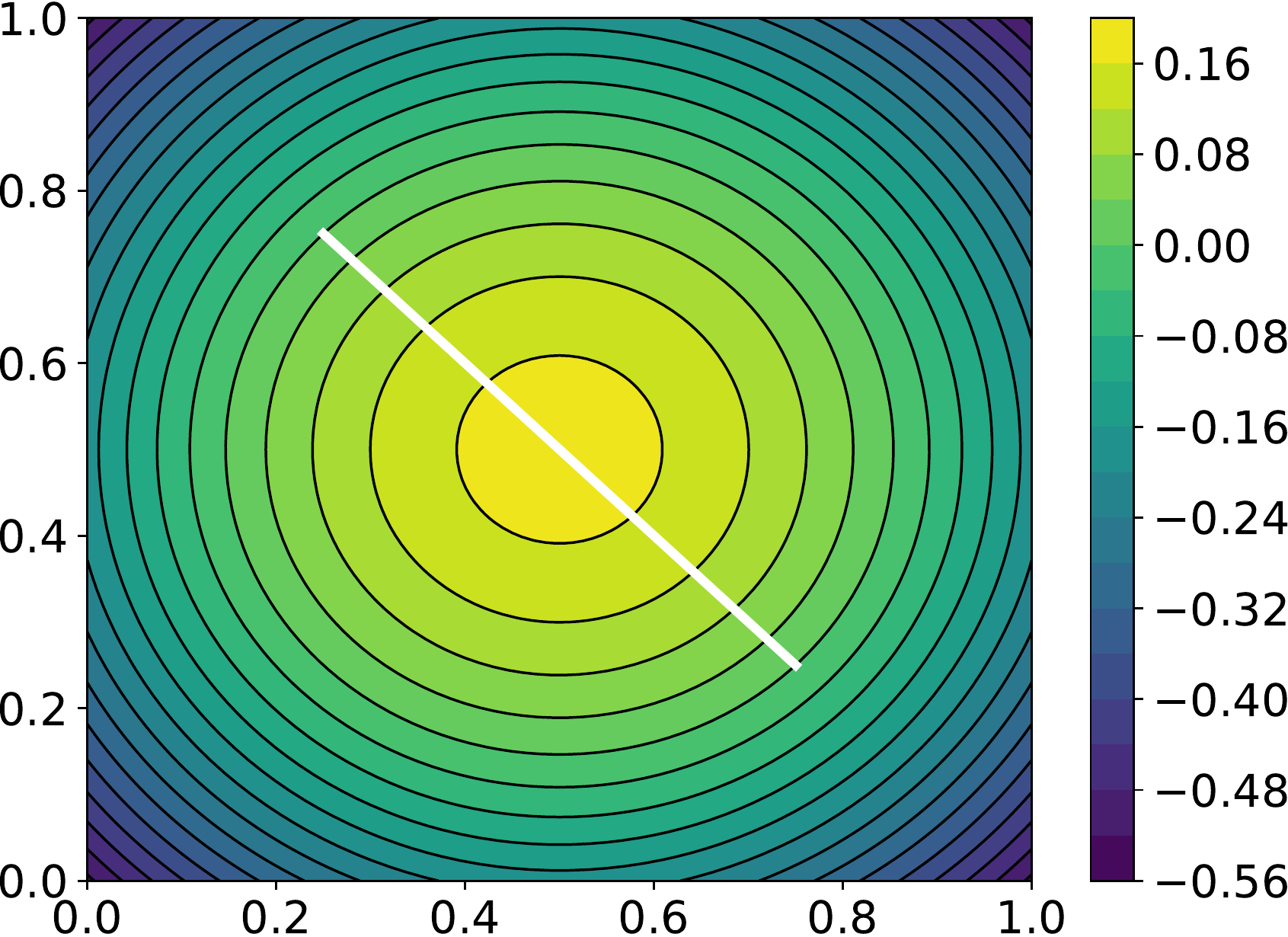} 
   \subcaption{$t(\mathbf x)$.}
  \label{fig:dtphib}
\end{subfigure}
\begin{subfigure}[t]{0.32\linewidth}
  \includegraphics[width=0.98\columnwidth,keepaspectratio,clip]{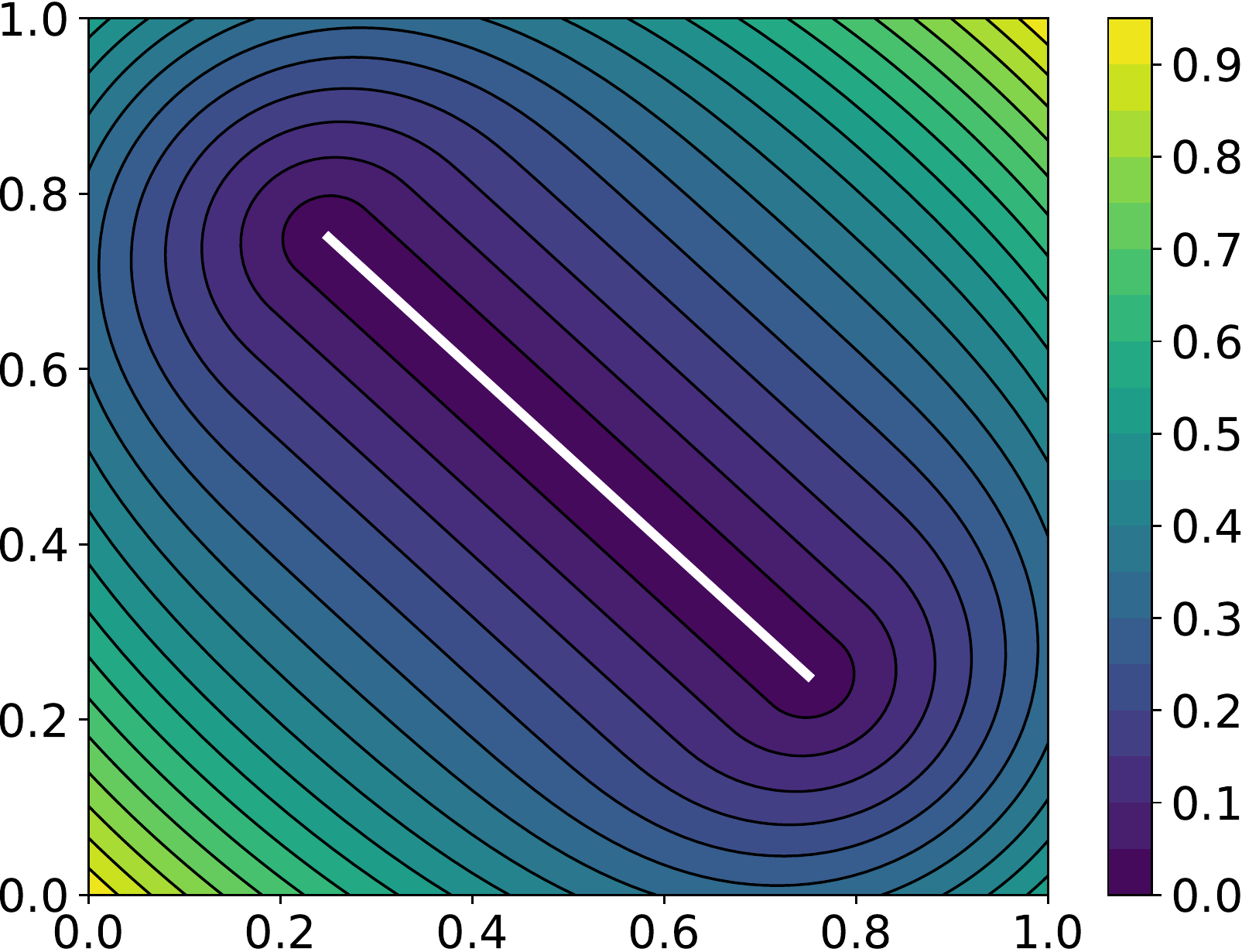}
   \subcaption{$\phi(\mathbf x)$.}
  \label{fig:dtphic}
\end{subfigure}
\caption{Representation of the signed distance function $d(\mathbf x)$ to a straight line (left), the trimming function $t(\mathbf x)$ (middle) and the approximate distance function $\phi(\mathbf x)$ to a segment (right).}
  \label{fig:dtphi}
\end{figure}

Assuming that $\Gamma_D$ can be expressed as the union of $n_s$ segments $\{s_1,\dots,s_{n_s}\}$, then the ADF to $\Gamma_D$, normalized up to order $m\ge1$, is defined as:
\begin{equation}\label{eq:adf_gammaD}
\phi=\dfrac{1}{\sqrt[m]{\frac{1}{\phi_1^m} + \frac{1}{\phi_2^m} + \dots + \frac{1}{\phi_{n_s}^m}}},
\end{equation}
where $\phi_i$ is the ADF to the segment $s_i$ (see \cite{biswas2004approximate}). We remark that an ADF normalized up to order $m\ge1$ is an ADF such that, for every regular point of $\Gamma_D$, the following holds:
\[
\phi = 0,\hspace{0.5cm} \dfrac{\partial \phi}{\partial n} = 1,\hspace{0.5cm} \dfrac{\partial^k\phi}{\partial n^k}=0 \hspace{0.2cm}(k=2,3,\dots,m).
\]
Such a normalization is useful to impose constraints associated with the solution derivatives and to obtain ADFs with about the same order of magnitude in every region of the domain $\Omega$. However, one of the main limitations of this approach with collocation-PINN is that $\Delta\phi$ tends to infinity near the vertices of $\Gamma_D$ (see Appendix \ref{sec:appendix} for an example). This phenomenon produces oscillations in the numerical solutions, hence collocation points that are close to such vertices should not be selected.
On the other hand, when only first derivatives are present in the weak formulation of second-order problems (as in the present study), then one can choose quadrature points 
that are very close to the vertices of $\Gamma_D$.

When a function $\overline g$ is not known, it is possible to construct it using transfinite interpolation. Let $g_i$ be a function such that ${g_i}_{|s_i} = g_{|s_i}$, then $\overline g$ can be defined as:
\[
\overline g = \sum_{i=1}^{n_s}w_i g_i,
\]
where $w_i$ is defined as:
\[
w_i = \dfrac{\prod_{j=1; j\ne i}^{n_s} \phi_j}{\sum_{k=1}^{n_s} \prod_{j=1;j\ne k}^{n_s} \phi_j} .
\]
Note that since $s_i$ is a segment, a function $g_i$ can be readily defined at any arbitrary point $(x,y)$ just by evaluating $g$ at the orthogonal projection of $(x,y)$ onto $s_i$.

\item[$\mathbf{M_C}$:] Exactly imposing the Dirichlet boundary conditions as in $\mathbf{M_B}$ but without normalizing the ADF. Therefore, we consider a different function $\phi$ in \eqref{eq:b_def}, namely
\[
\phi = \prod_{i=1}^{n_s}\phi_i.
\]
This ensures that $\phi$ and all its derivatives exist and are bounded in $\overline\Omega$, although $\phi$ may be very small in regions close to many segments $s_i$.

\item[$\mathbf{M_D}$:] Using Nitsche's method \cite{nitscheweak}.  The goal of this method is to variationally impose the Dirichlet boundary conditions. In doing so, the network architecture is not modified with additional layers (as in $\mathbf{M_B}$ and $\mathbf{M_C}$) and a single objective function suffices for network training.

To do so, one enlarges the space $V_h$ to contain all piecewise polynomials of order $k_{\text{test}}$ defined on $\mathcal T_h$  and modifies the residuals defined in~\eqref{eq:residuals} in the following way:
\begin{equation}\label{eq:nitsche_residuals}
\begin{split}
r_{h,i}(w) = F_h(\varphi_i^v) - a_h(w, \varphi_i^v) &+ \int_{\Gamma_D}(w-g)\frac{\partial \phi_i^v}{\partial n} \\
&+ \gamma \int_{\Gamma_D}h^{-1}(g-w)\phi_i^v, 
\quad i\in I_h,
\end{split}
\end{equation}
where $\gamma$ is a positive constant satisfying $\gamma\ge\gamma_0$ for a suitable $\gamma_0>0$ and $I_h$ is now an enlarged index set corresponding to the enlarged basis $\{\phi_i^v: i\in I_H\}$. Thanks to the scaling term $h^{-1}$ that magnifies the quantity $\int_{\Gamma_D}(g-w)\phi_i^v$ when fine meshes are used, the choice of $\gamma$ is not as important as the one of $\lambda$ in method $\mathbf{M_A}$. This property is confirmed by numerical results shown in Figures \ref{fig:eq_2_hole_reg_hard_ma} and \ref{fig:eq_2_hole_reg_hard_md}. Since there is no ambiguity, we maintain the same symbols $V_h$, $I_h$ and $\{\phi_i^v\}$ introduced in Section \ref{sec:loss_definition}; they always represent the enlarged sets when method $\mathbf M_D$ is considered.
Note that when $w$ satisfies the Dirichlet boundary conditions, the terms added in \eqref{eq:nitsche_residuals} vanish.
\end{itemize}

We point out that method $\mathbf{M_A}$ is often referred to as {\em soft boundary condition imposition}, whereas $\mathbf{M_B}$ and $\mathbf{M_C}$ are known as {\em hard boundary condition impositions}. Hence, we can treat $\mathbf{M_D}$ as
{\em weak boundary condition imposition}.

\section{Numerical results} \label{sec:numerical_results}
In this section, the methods $\mathbf{M_A}$, $\mathbf{M_B}$, $\mathbf{M_C}$ and $\mathbf{M_D}$ discussed in Section \ref{sec:used_approaches} are analyzed and compared. In each numerical experiment the neural network is a fully-connected feed-forward neural network as described in Section \ref{sec:nn_description}. The corresponding architecture is composed of 4 hidden layers with 50 neurons in each layer and with the hyperbolic tangent as the activation function, while the output layer is a linear layer with one or two neurons.

In order to properly minimize the loss function we use the first-order ADAM optimizer \cite{kingma2014adam} with an exponentially decaying learning rate, and after a prescribed number of epochs, the second-order BFGS method \cite{wright1999numerical} is used until a maximum number of iterations is reached, or it is not possible to further improve the objective function (i.e. when two consecutive iterates are identical, up to machine precision). When the interpolated VPINN is used, the training set consists of all the interpolation nodes  $\{x_1^{\mathcal I},\dots,x_{\text{dim}(U_H)}^{\mathcal I}\}$ and no regularization is applied since the interpolation operator already filters the neural network high frequencies out. Instead, when the PINN is used, the training set contains a set of $\text{dim}(U_H)$ control points inside the domain $\Omega$, and when $\mathbf{M_A}$ is employed, a set of approximately $\sqrt{\text{dim}(U_H)}$ control points on the boundary $\partial \Omega$. Moreover, in order to stabilize the PINN, the $L^2$ regularization term 
\begin{equation}\label{eq:l2_reg}
\mathcal L_{\text{reg}}(\mathbf{w}^{\mathcal{N\!N}}) = \lambda_{\text{reg}}\Vert \mathbf{w}^{\mathcal{N\!N}}\Vert_2^2 
\end{equation}
is added to the loss function, where $\mathbf{w}^{\mathcal{N\!N}}$ is the vector containing all the neural network weights defined in Section \ref{sec:nn_description} and $\lambda_{\text{reg}}=10^{-6}$. The value of this parameter has been chosen through several numerical experiments to minimize the $H^1$ norm of the error.

The computer code to perform the numerical experiments is written in Python, while the neural networks and the optimizers are implemented using the open-source Python package Tensorflow \cite{tensorflow2015-whitepaper}. The loss function gradient with respect to the neural network weights and the PINN output gradient with respect to the spatial coordinates are always computed with automatic differentiation that is 
available in Tensorflow \cite{baydin2018automatic}. 
On the other hand,
the VPINN output gradient with respect to its input is computed by means of suitable projection matrices as described in \cite{berrone2022variational}.

\subsection{Rate of convergence for second-order elliptic problems}\label{sec:conv_rates}
\begin{figure}[t!]
\centering 
\begin{subfigure}[t]{0.495\linewidth}
  \includegraphics[width=0.98\columnwidth,keepaspectratio,clip]{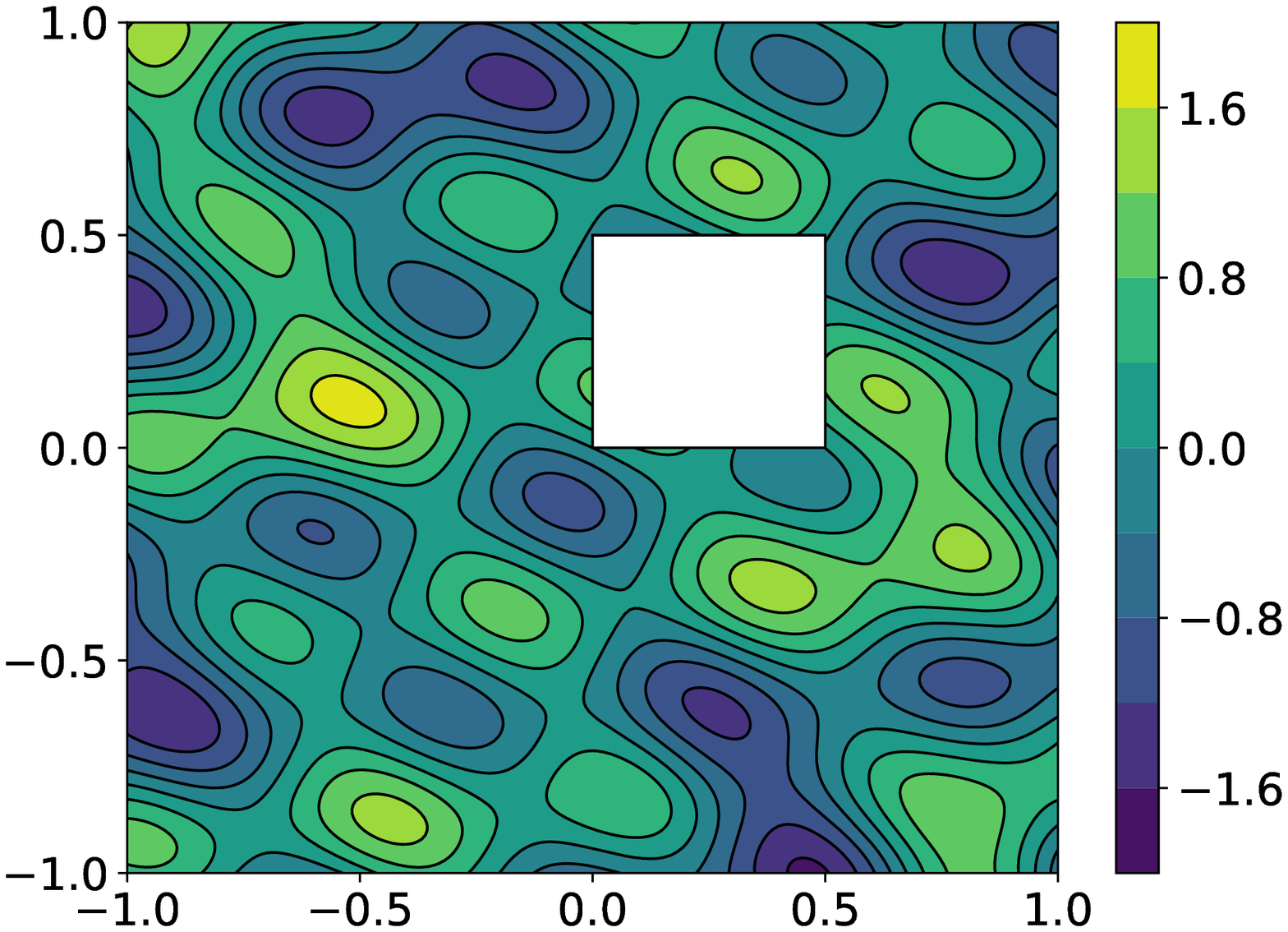} 
  \caption{}
  \label{fig:exact_sol_and_error_a}
\end{subfigure}
\begin{subfigure}[t]{0.495\linewidth}
  \includegraphics[width=0.98\columnwidth,keepaspectratio,clip]{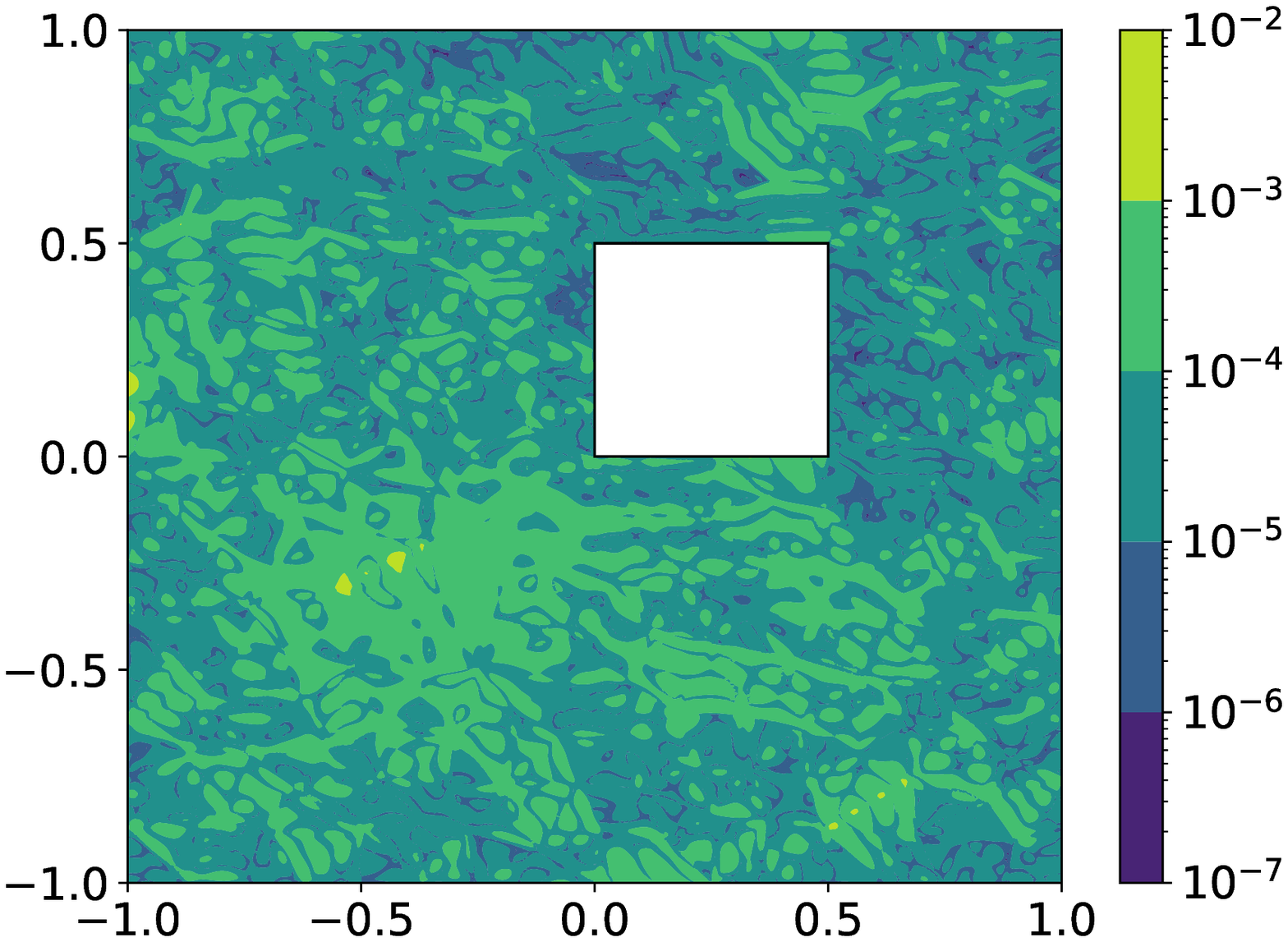} 
  \caption{}
  \label{fig:exact_sol_and_error_b}
\end{subfigure}
  \caption{Exact solution $u$ (left) and 
  a plot of the absolute error with VPINN and method $\bf M_B$ in which the Dirichlet boundary conditions are imposed on every edge of $\partial \Omega$ (right).}
  \label{fig:exact_sol_and_error}
\end{figure}

We focus on the VPINN model and show that the a priori error estimate proved in \cite{berrone2022variational} for second-order elliptic problems holds even on varying the way in which the boundary conditions are imposed. 
On letting $\bm{x} := (x,y)$, we 
consider problem \eqref{eq:elliptic-model-pb} in 
the domain $\Omega = (-1,1)^2\setminus(0,0.5)^2$ with the physical parameters
\[
\mu(\bm{x}) = 2+\sin(x+2y), \quad
\beta(\bm{x}) = \left\{\sqrt{x-y^2+5},\sqrt{y-x^2+5}\right\}, \quad
\sigma(\bm{x})=e^{\frac x2-\frac y3}+2 .
\]
We consider two test cases. In the first one the Dirichlet boundary conditions and forcing term are chosen so that the exact solution is
\begin{equation}\label{eq:sol2}
u(\bm{x}) = \cos(5(x+y/2)) + (x+y/2)^2,
\end{equation}
whereas in the second one they are chosen such that the exact solution is more oscillatory. Its expression is:
\begin{equation}\label{eq:sol5}
u(\bm{x}) = \sin\left[3x(x-y)\right]\cos(4y+x)+\sin\left[5(x+2y)\right]\cos\left[3(y-2x)\right].
\end{equation}
Such a solution is shown in 
Fig.~\ref{fig:exact_sol_and_error_a}, whereas an example of numerical error corresponding to the VPINN in which Dirichlet boundary conditions are imposed using method $\mathbf{M_B}$ is shown in Fig.~\ref{fig:exact_sol_and_error_b}; it exhibits a rather uniform distribution of the error, which is not localized near boundaries. We remark that in these numerical tests and in the subsequent ones, the function $\overline{g}$ used in $\mathbf{M_B}$ and $\mathbf{M_C}$ is computed via transfinite interpolation.

\begin{figure}[t!]
\centering 
\begin{subfigure}[t]{0.32\linewidth}
  \includegraphics[width=0.98\columnwidth,keepaspectratio,clip]{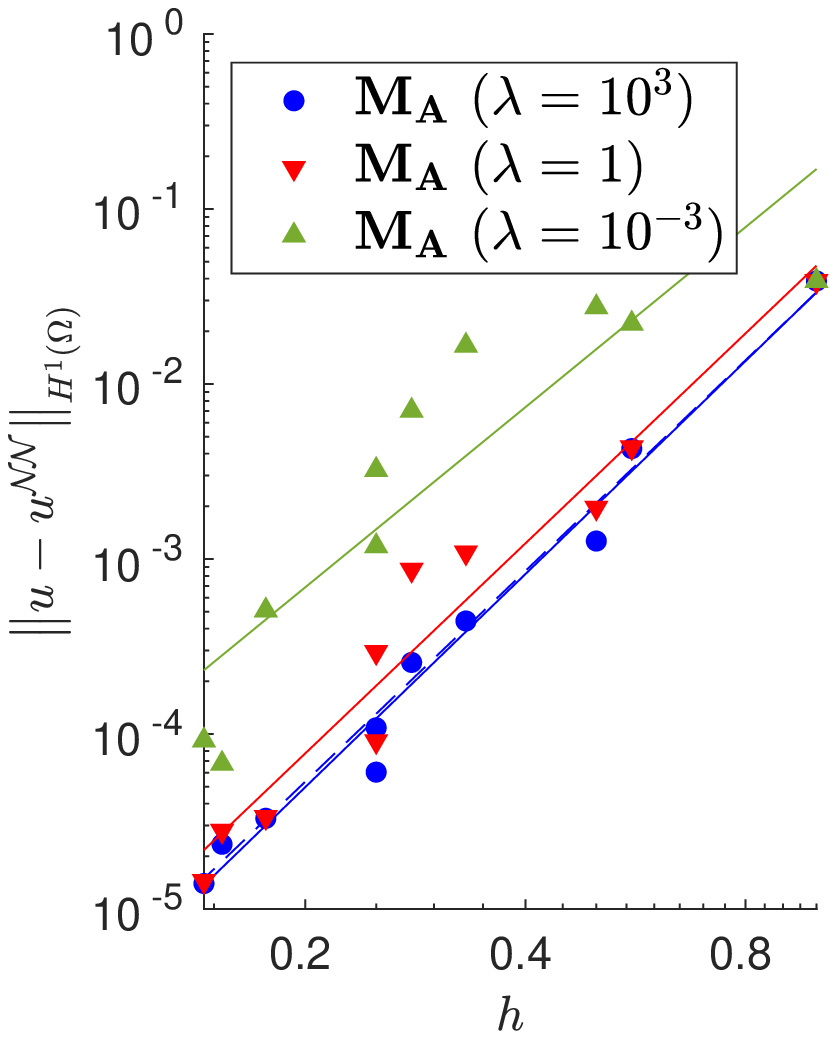} 
   \subcaption{$q=3$, $k_{\text{test}}=1$, $k_\text{int}=4$.}
   \label{fig:eq_2_hole_regular_ma_a}
\end{subfigure}
\begin{subfigure}[t]{0.32\linewidth}
  \includegraphics[width=0.98\columnwidth,keepaspectratio,clip]{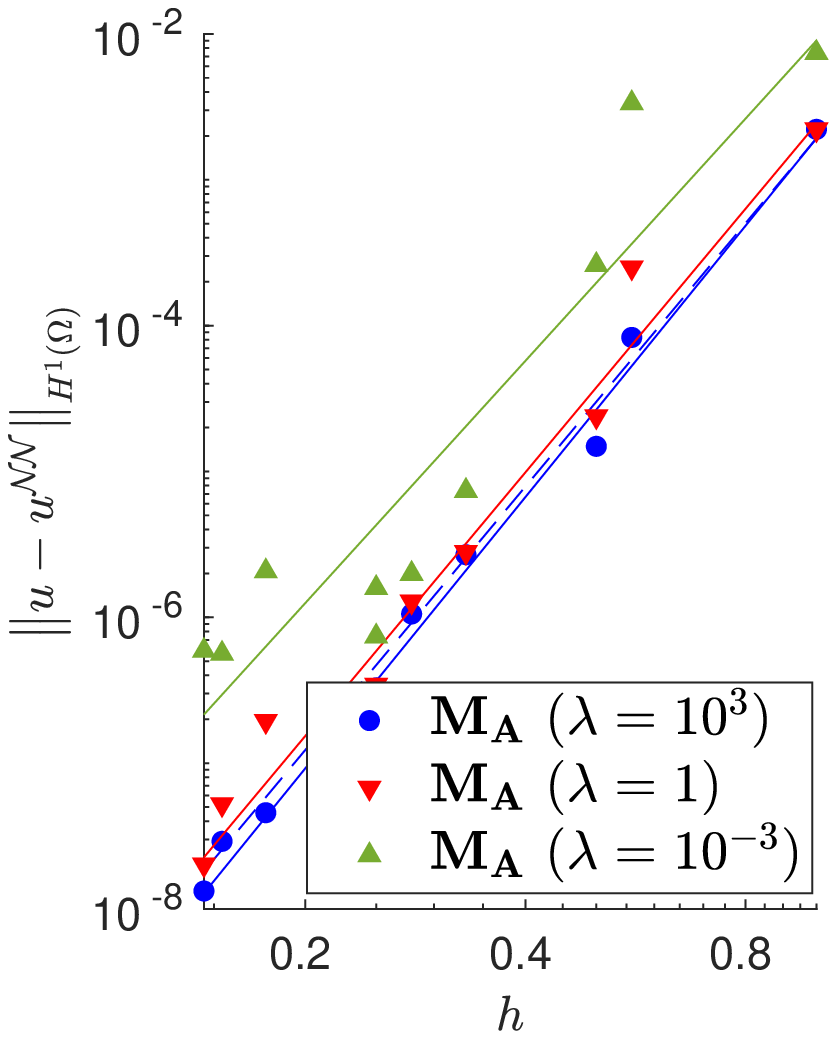} 
   \subcaption{$q=5$, $k_{\text{test}}=1$, $k_\text{int}=6$.}
   \label{fig:eq_2_hole_regular_ma_b}
\end{subfigure}
\begin{subfigure}[t]{0.32\linewidth}
  \includegraphics[width=0.98\columnwidth,keepaspectratio,clip]{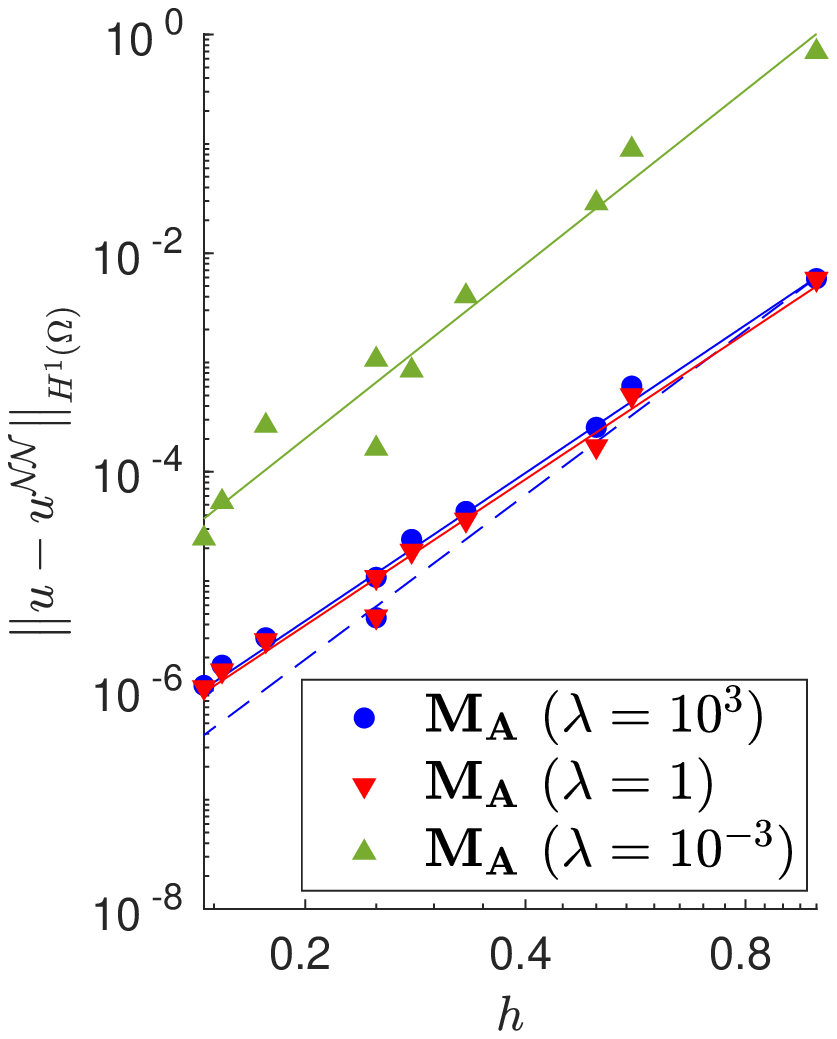}
   \subcaption{$q=5$, $k_{\text{test}}=2$, $k_\text{int}=5$.}
   \label{fig:eq_2_hole_regular_ma_c}
\end{subfigure}
  \caption{Error decay obtained with $\mathbf{M_A}$ and different values of $\lambda$. Forcing term and Dirichlet boundary conditions are set such that the exact solution is \eqref{eq:sol2}. The theoretical convergence rate is $k_\text{int}$. (a) Convergence rates: 4.04 $(\lambda=10^{3})$, 3.99 $(\lambda=1)$, 3.42 $(\lambda=10^{-3})$. (b) Convergence rates: 6.18 $(\lambda=10^{3})$, 6.00 $(\lambda=1)$, 5.52 $(\lambda=10^{-3})$. (c) Convergence rates: 4.50 $(\lambda=10^{3})$, 4.44 $(\lambda=1)$, 5.29 $(\lambda=10^{-3})$. }
  \label{fig:eq_2_hole_regular_ma}
\end{figure}

\begin{figure}[t!]
\centering 
\begin{subfigure}[t]{0.32\linewidth}
  \includegraphics[width=0.98\columnwidth,keepaspectratio,clip]{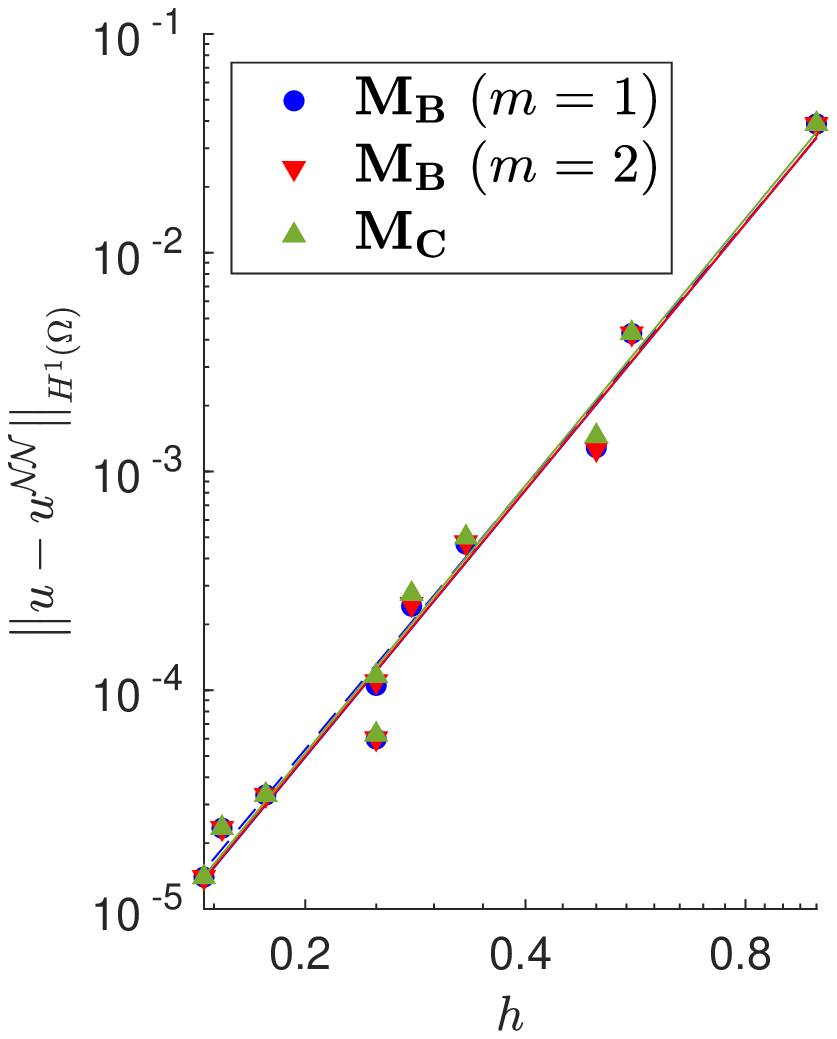} 
   \subcaption{$q=3$, $k_{\text{test}}=1$, $k_\text{int}=4$.}
  \label{fig:eq_2_hole_regular_mbmc_a}
\end{subfigure}
\begin{subfigure}[t]{0.32\linewidth}
  \includegraphics[width=0.98\columnwidth,keepaspectratio,clip]{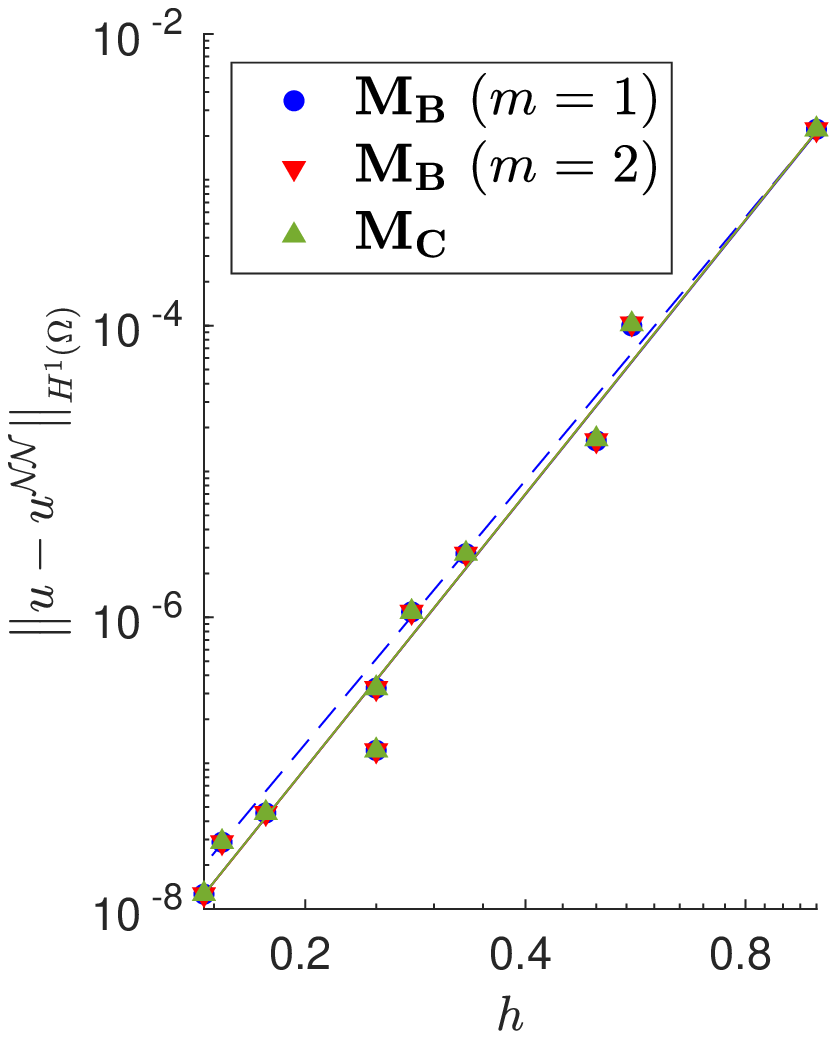} 
   \subcaption{$q=5$, $k_{\text{test}}=1$, $k_\text{int}=6$.}
  \label{fig:eq_2_hole_regular_mbmc_b}
\end{subfigure}
\begin{subfigure}[t]{0.32\linewidth}
  \includegraphics[width=0.98\columnwidth,keepaspectratio,clip]{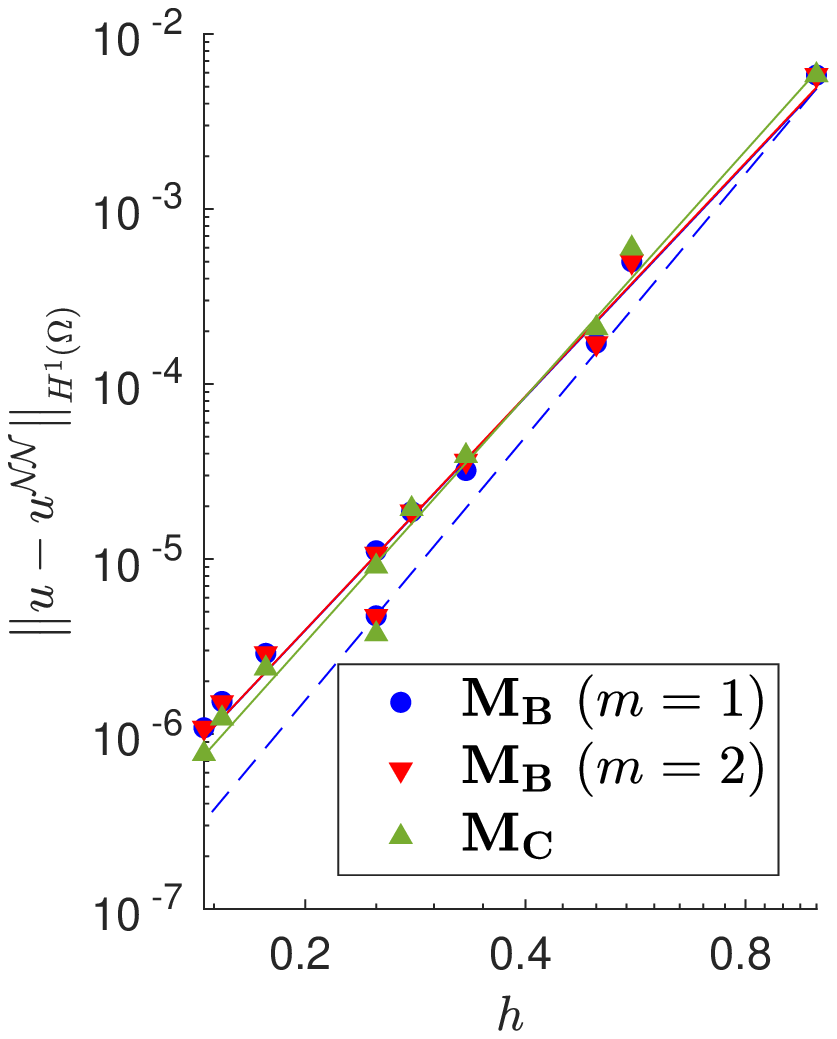}
   \subcaption{$q=5$, $k_{\text{test}}=2$, $k_\text{int}=5$.}
  \label{fig:eq_2_hole_regular_mbmc_c}
\end{subfigure}
  \caption{Error decay obtained with $\mathbf{M_B}$, with different values of $m$, and $\mathbf{M_C}$. Forcing term and Dirichlet boundary conditions are set such that the exact solution is \eqref{eq:sol2}. The theoretical convergence rate is $k_\text{int}$. (a) Convergence rates: 4.05 $(\mathbf{M_B}, m=1)$, 4.05 $(\mathbf{M_B}, m=2)$, 4.06 $(\mathbf{M_C})$. (b) Convergence rates: 6.24 $(\mathbf{M_B}, m=1)$, 6.25 $(\mathbf{M_B}, m=2)$, 6.25 $(\mathbf{M_C})$. (c) Convergence rates: 4.43 $(\mathbf{M_B}, m=1)$, 4.43 $(\mathbf{M_B}, m=2)$, 4.67 $(\mathbf{M_C})$. }
  \label{fig:eq_2_hole_regular_mbmc}
\end{figure}

\begin{figure}[t!]
\centering 
\begin{subfigure}[t]{0.32\linewidth}
  \includegraphics[width=0.98\columnwidth,keepaspectratio,clip]{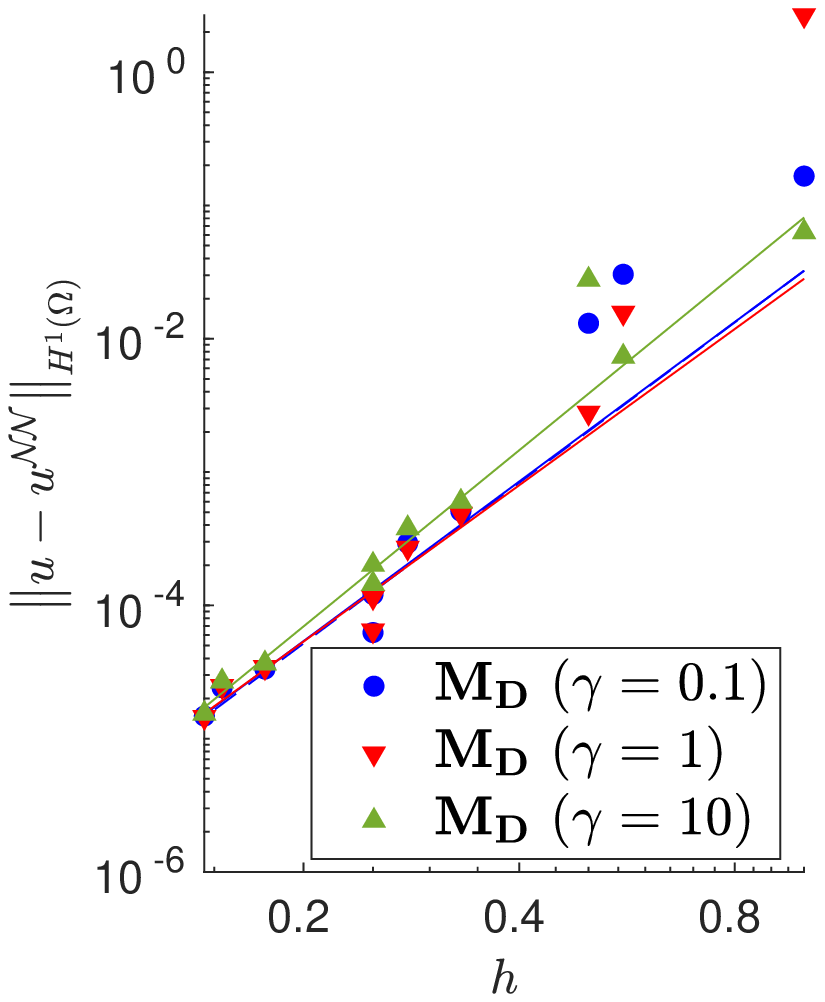} 
   \subcaption{$q=3$, $k_{\text{test}}=1$, $k_\text{int}=4$.}
  \label{fig:eq_2_hole_regular_md_a}
\end{subfigure}
\begin{subfigure}[t]{0.32\linewidth}
  \includegraphics[width=0.98\columnwidth,keepaspectratio,clip]{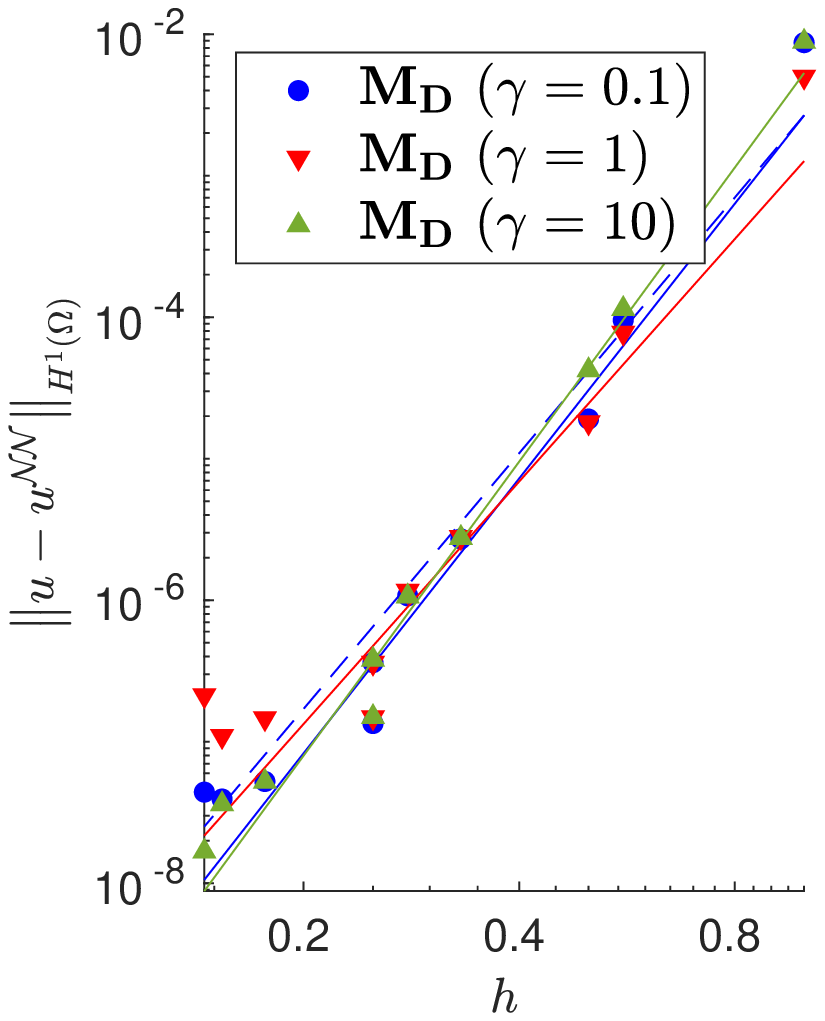} 
   \subcaption{$q=5$, $k_{\text{test}}=1$, $k_\text{int}=6$.}
  \label{fig:eq_2_hole_regular_md_b}
\end{subfigure}
\begin{subfigure}[t]{0.32\linewidth}
  \includegraphics[width=0.98\columnwidth,keepaspectratio,clip]{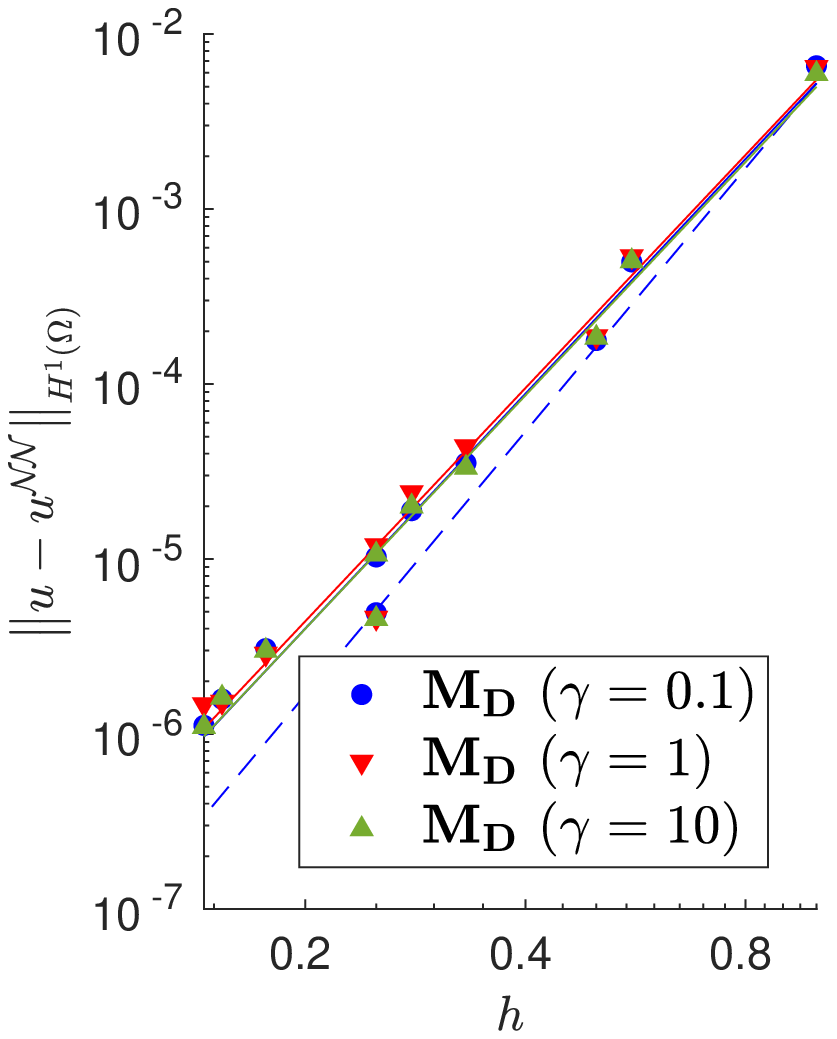}
   \subcaption{$q=5$, $k_{\text{test}}=2$, $k_\text{int}=5$.}
  \label{fig:eq_2_hole_regular_md_c}
\end{subfigure}
  \caption{Error decay obtained with $\mathbf{M_D}$ and different values of $\gamma$. Forcing term and Dirichlet boundary conditions are set such that the exact solution is \eqref{eq:sol2}. The theoretical convergence rate is $k_\text{int}$. (a) Convergence rates: 3.98$(\gamma=0.1)$, 3.89 $(\gamma=1)$, 4.39 $(\gamma=10)$. (b) Convergence rates: 6.45 $(\gamma=0.1)$, 5.69 $(\gamma=1)$, 6.90 $(\gamma=10)$. (c) Convergence rates: 4.46 $(\gamma=0.1)$, 4.43 $(\gamma=1)$, 4.43 $(\gamma=10)$. }
  \label{fig:eq_2_hole_regular_md}
\end{figure}

We vary both the order of the quadrature rule and the degree of the test functions, and train the same model with different meshes and impose the Dirichlet boundary conditions with the proposed approaches. In Figures 
\ref{fig:eq_2_hole_regular_ma_a}--\ref{fig:eq_2_hole_regular_ma_c}, 
\ref{fig:eq_2_hole_regular_mbmc_a}--\ref{fig:eq_2_hole_regular_mbmc_c} and 
\ref{fig:eq_2_hole_regular_md_a}--\ref{fig:eq_2_hole_regular_md_c}, 
in which the exact solution is the one in \eqref{eq:sol2}, we observe close agreement with the results shown in \cite{berrone2022variational}. In fact, when the loss is properly minimized, all the approaches perform comparably and the corresponding empirical convergence rates are always close to the theoretical rate of $k_{\text{int}}=q+2-k_{\text{test}}$. We point out that in \cite{berrone2022variational} we prove that, when the solution is regular enough and a method similar to $\mathbf{M_C}$ is used to enforce the boundary conditions, the convergence rate is $k_{\text{int}}=q+2-k_{\text{test}}$. Here, instead we show that the same behaviour is observed even if the boundary conditions are enforced in different ways. Note in particular that the choice $m=1$ or $m=2$ in $\mathbf{M_B}$, and the choice $\gamma=0.1$, $\gamma=1$ or $\gamma=10$ in $\mathbf{M_D}$ yields nearly identical results (see Figures \ref{fig:eq_2_hole_regular_mbmc} and \ref{fig:eq_2_hole_regular_md}).

\begin{figure}[t!]
\centering 
\begin{subfigure}[t]{0.49\linewidth}
  \includegraphics[width=0.98\columnwidth,keepaspectratio,clip]{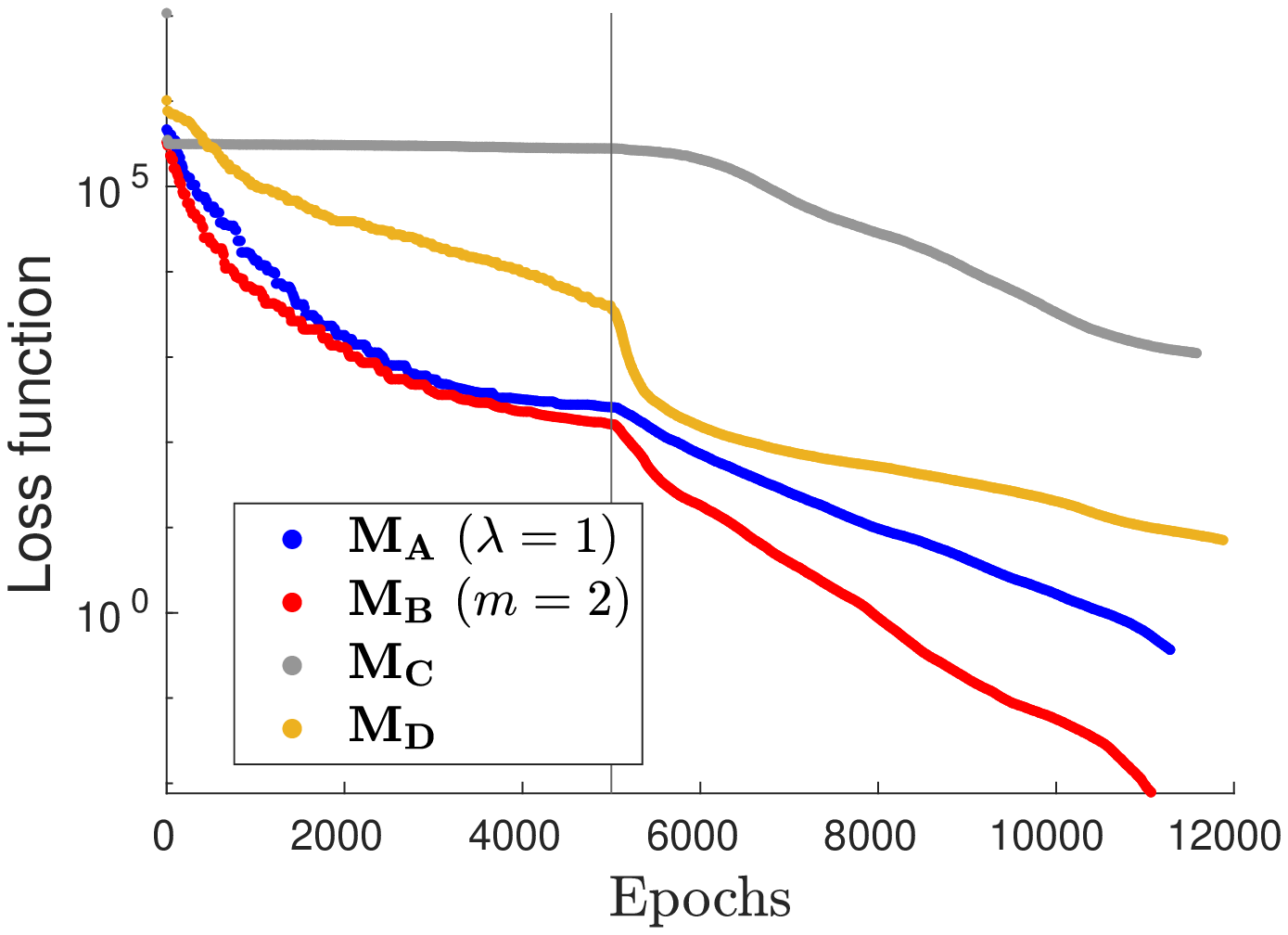} 
  \caption{Loss function behaviour during training.}
  \label{fig:eq2_sol5_loss_}
\end{subfigure}
\begin{subfigure}[t]{0.49\linewidth}
  \includegraphics[width=0.98\columnwidth,keepaspectratio,clip]{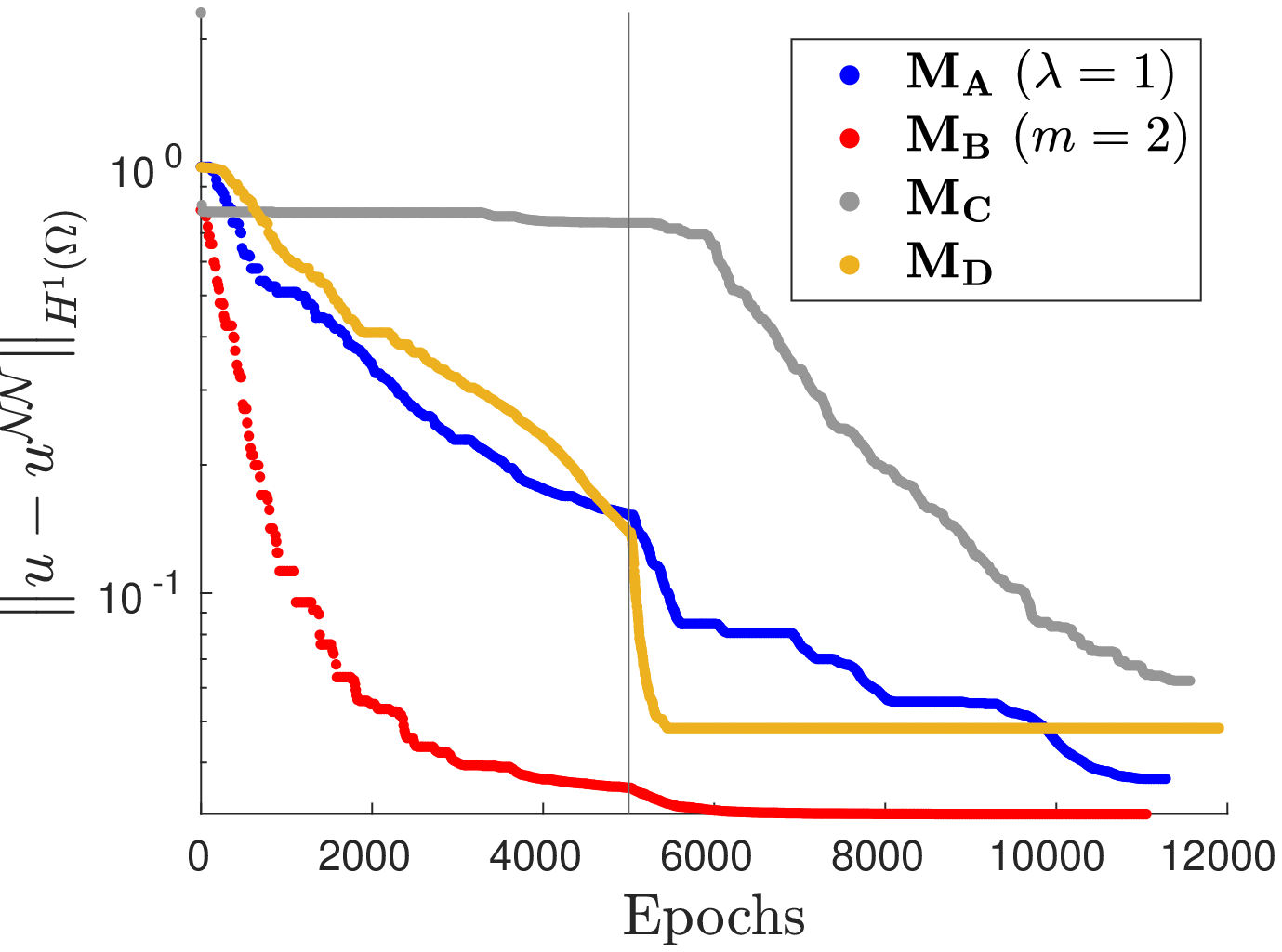} 
  \caption{$H^1$ error during training.}
  \label{fig:eq2_sol5_error_}
\end{subfigure}
  \caption{Training loss (left) and $H^1$ error
  prediction (right) for the VPINN.
  The first 5000 epochs are performed with a standard ADAM optimizer, the remaining ones with the BFGS optimizer. The exact solution is given in \eqref{eq:sol5}.}
  \label{fig:eq2_sol5_loss_and_errors}
\end{figure}

We highlight that the different methods, while delivering similar empirical convergence rates with respect to mesh refinement, exhibit very different performance during training. To observe this phenomenon, let us train multiple identical neural networks on the same mesh but impose 
the Dirichlet boundary conditions in different ways. Here we only consider quadrature rules of order $q=3$ and piecewise linear test functions. The values of the loss function and of the $H^1$ error prediction during training are presented 
in Figures~\ref{fig:eq2_sol5_loss_} and~\ref{fig:eq2_sol5_error_}, respectively. A vertical line separates the epochs where the ADAM optimizer is used from the ones where the BFGS optimizer is used. 

It can be noted that the most efficient method is $\mathbf{M_B}$, as it converges faster and to a more accurate solution, while method $\mathbf{M_D}$ is characterized by very fast convergence only when the BFGS optimizer is adopted. Such an optimizer is also crucial to train the VPINN with $\mathbf{M_C}$; in fact the corresponding error does not decrease when the ADAM optimizer is used. Instead, the convergence obtained using method $\mathbf{M_A}$ seems independent of the choice of the optimizer. It is important to remark that all the loss functions are decreasing even when the error is constant. This implies that there exist other sources of error that dominate and that a very small loss function does not ensure a very accurate solution; this phenomenon is also observed in Fig.~3 of \cite{berrone2022solving} and is discussed in greater detail therein.

\begin{figure}[t!]
\centering 
\begin{subfigure}[t]{0.32\linewidth}
  \includegraphics[width=0.98\columnwidth,keepaspectratio,clip]{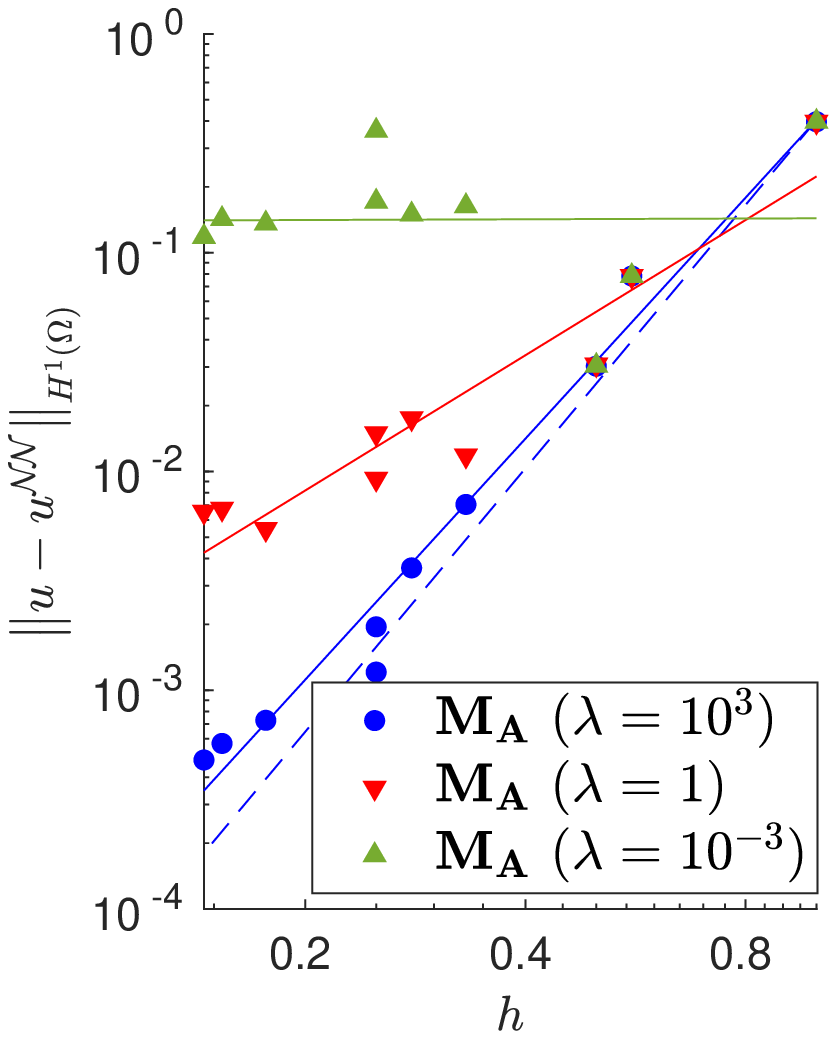} 
   \subcaption{$q=3$, $k_{\text{test}}=1$, $k_\text{int}=4$.}
  \label{fig:eq_2_hole_reg_hard_ma_a}
\end{subfigure}
\begin{subfigure}[t]{0.32\linewidth}
  \includegraphics[width=0.98\columnwidth,keepaspectratio,clip]{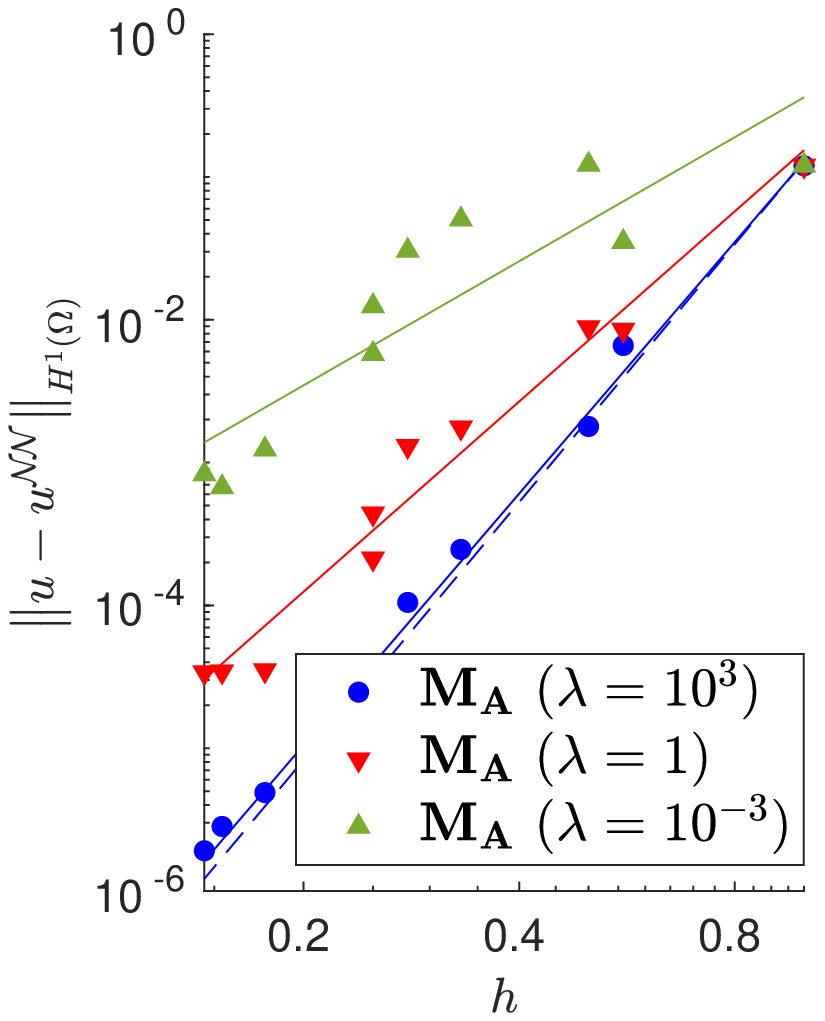} 
   \subcaption{$q=5$, $k_{\text{test}}=1$, $k_\text{int}=6$.}
  \label{fig:eq_2_hole_reg_hard_ma_b}
\end{subfigure}
\begin{subfigure}[t]{0.32\linewidth}
  \includegraphics[width=0.98\columnwidth,keepaspectratio,clip]{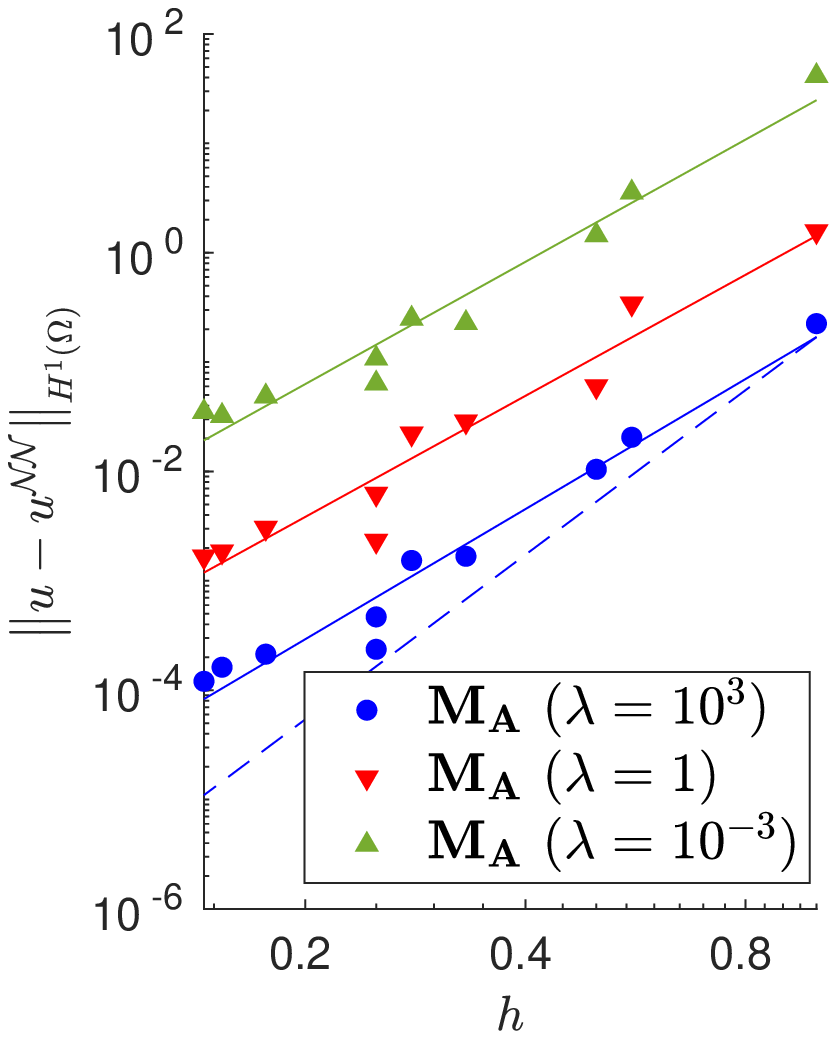}
   \subcaption{$q=5$, $k_{\text{test}}=2$, $k_\text{int}=5$.}
  \label{fig:eq_2_hole_reg_hard_ma_c}
\end{subfigure}
  \caption{Error decay obtained with $\mathbf{M_A}$ and different values of $\lambda$. Forcing term and Dirichlet boundary conditions are set such that the exact solution is \eqref{eq:sol5}. The theoretical convergence rate is $k_\text{int}$. (a) Convergence rates: 3.66 $(\lambda=10^{3})$, 2.05 $(\lambda=1)$, 0.01 $(\lambda=10^{-3})$. (b) Convergence rates: 5.85 $(\lambda=10^{3})$, 4.42 $(\lambda=1)$, 2.89 $(\lambda=10^{-3})$. (c) Convergence rates: 3.95 $(\lambda=10^{3})$, 3.68 $(\lambda=1)$, 3.71 $(\lambda=10^{-3})$. }
  \label{fig:eq_2_hole_reg_hard_ma}
\end{figure}

\begin{figure}[t!]
\centering 
\begin{subfigure}[t]{0.32\linewidth}
  \includegraphics[width=0.98\columnwidth,keepaspectratio,clip]{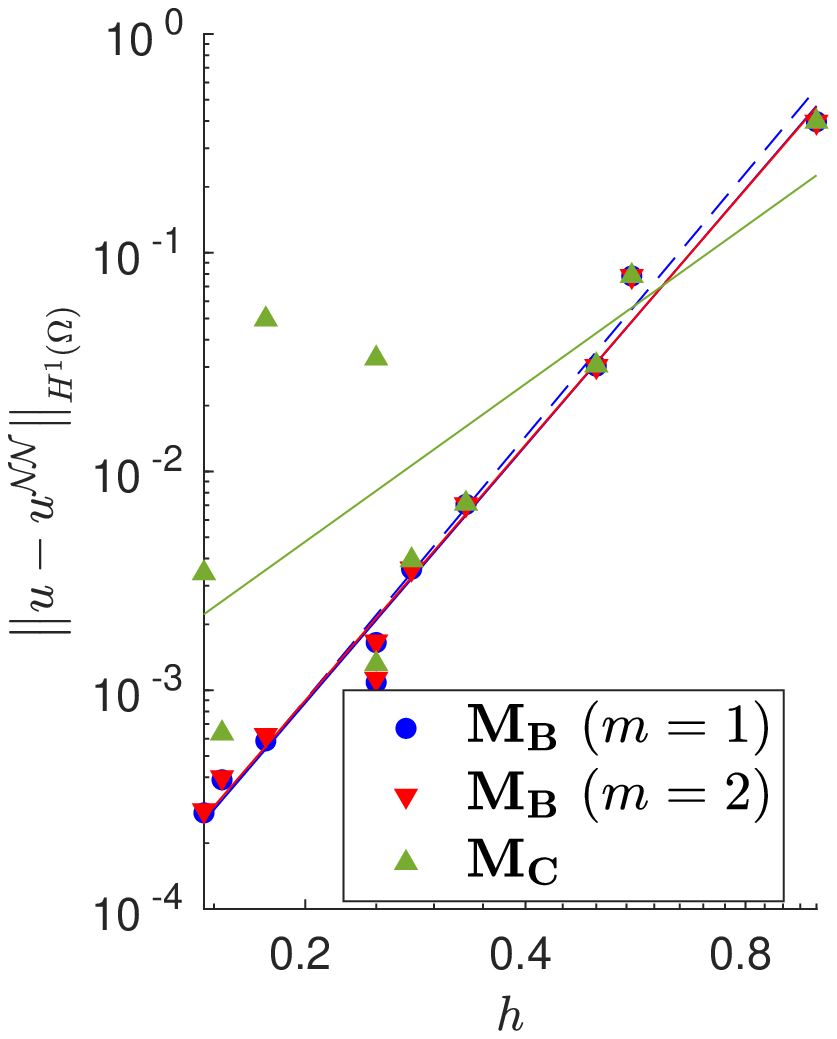} 
   \subcaption{$q=3$, $k_{\text{test}}=1$, $k_\text{int}=4$.}
  \label{fig:eq_2_hole_reg_hard_mbmc_a}
\end{subfigure}
\begin{subfigure}[t]{0.32\linewidth}
  \includegraphics[width=0.98\columnwidth,keepaspectratio,clip]{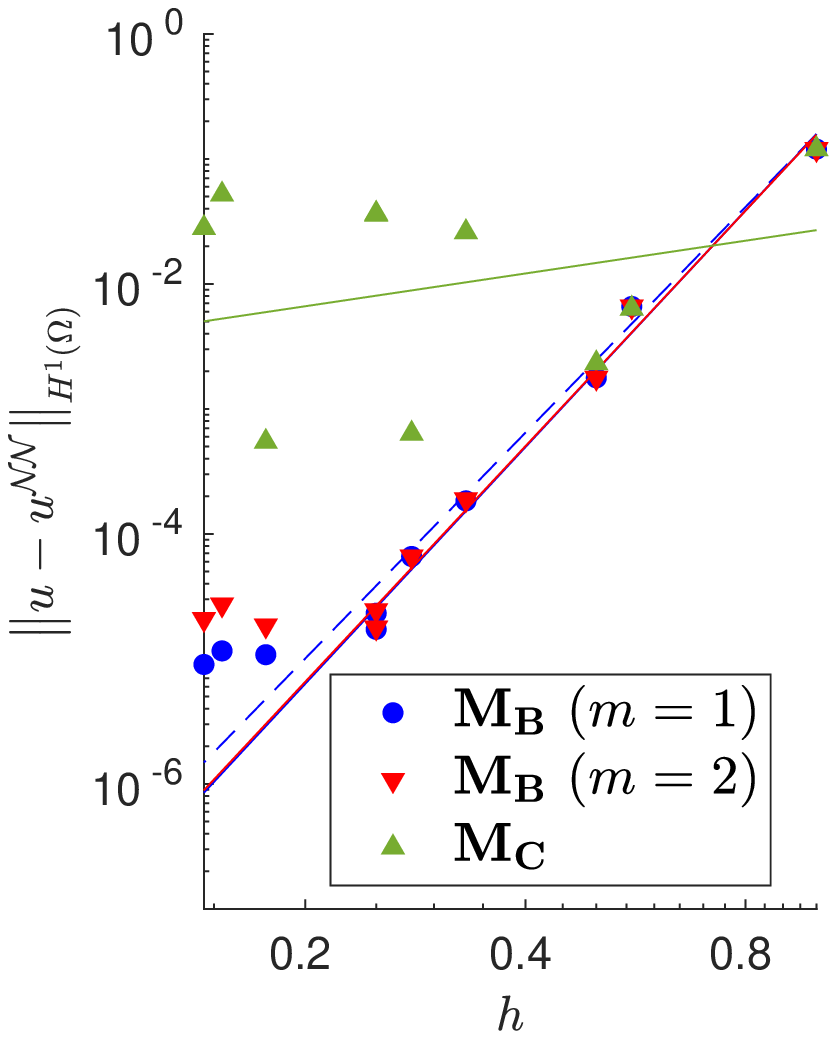} 
   \subcaption{$q=5$, $k_{\text{test}}=1$, $k_\text{int}=6$.}
  \label{fig:eq_2_hole_reg_hard_mbmc_b}
\end{subfigure}
\begin{subfigure}[t]{0.32\linewidth}
  \includegraphics[width=0.98\columnwidth,keepaspectratio,clip]{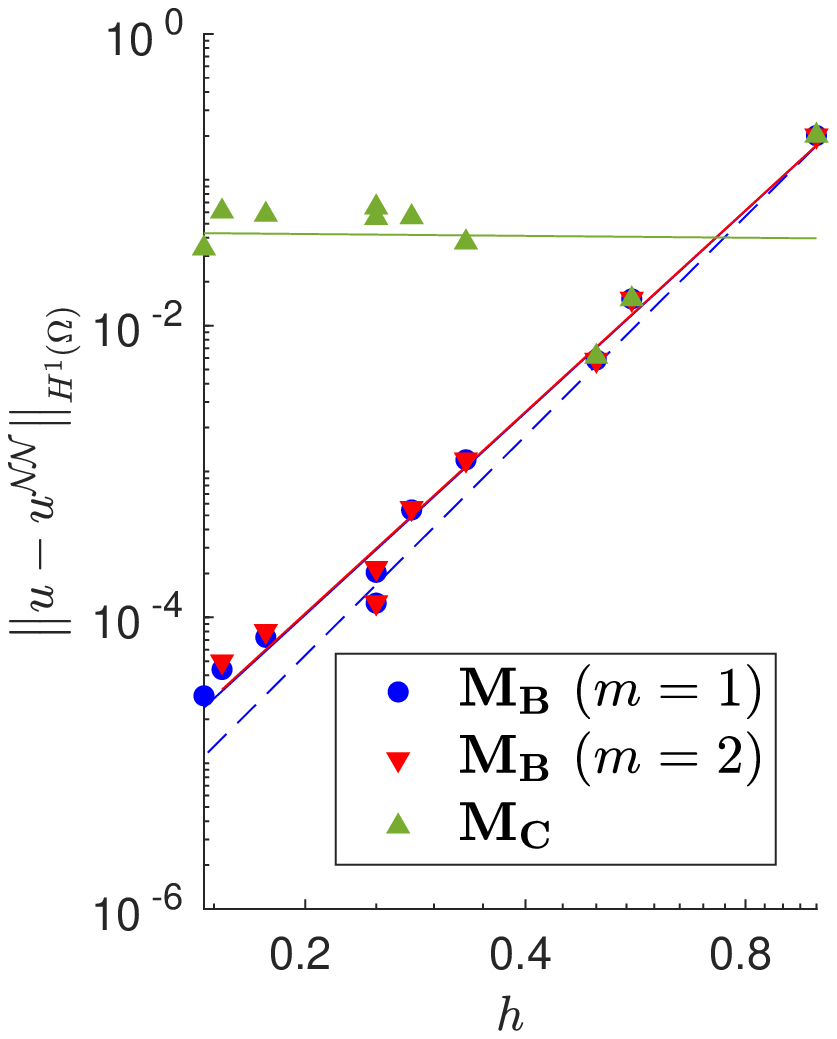}
   \subcaption{$q=5$, $k_{\text{test}}=2$, $k_\text{int}=5$.}
  \label{fig:eq_2_hole_reg_hard_mbmc_c}
\end{subfigure}
  \caption{Error decay obtained with $\mathbf{M_B}$, with different values of $m$, and $\mathbf{M_C}$. Forcing term and Dirichlet boundary conditions are set such that the exact solution is \eqref{eq:sol5}. The theoretical convergence rate is $k_\text{int}$. (a) Convergence rates: 3.90 $(\mathbf{M_B}, m=1)$, 3.88 $(\mathbf{M_B}, m=2)$, 2.36 $(\mathbf{M_C})$. (b) Convergence rates: 5.85 $(\mathbf{M_B}, m=1)$, 4.42 $(\mathbf{M_B}, m=2)$, 2.89 $(\mathbf{M_C})$. (c) Convergence rates: 4.60 $(\mathbf{M_B}, m=1)$, 4.59 $(\mathbf{M_B}, m=2)$, -0.43 $(\mathbf{M_C})$. }
  \label{fig:eq_2_hole_reg_hard_mbmc}
\end{figure}

\begin{figure}[t!]
\centering 
\begin{subfigure}[t]{0.32\linewidth}
  \includegraphics[width=0.98\columnwidth,keepaspectratio,clip]{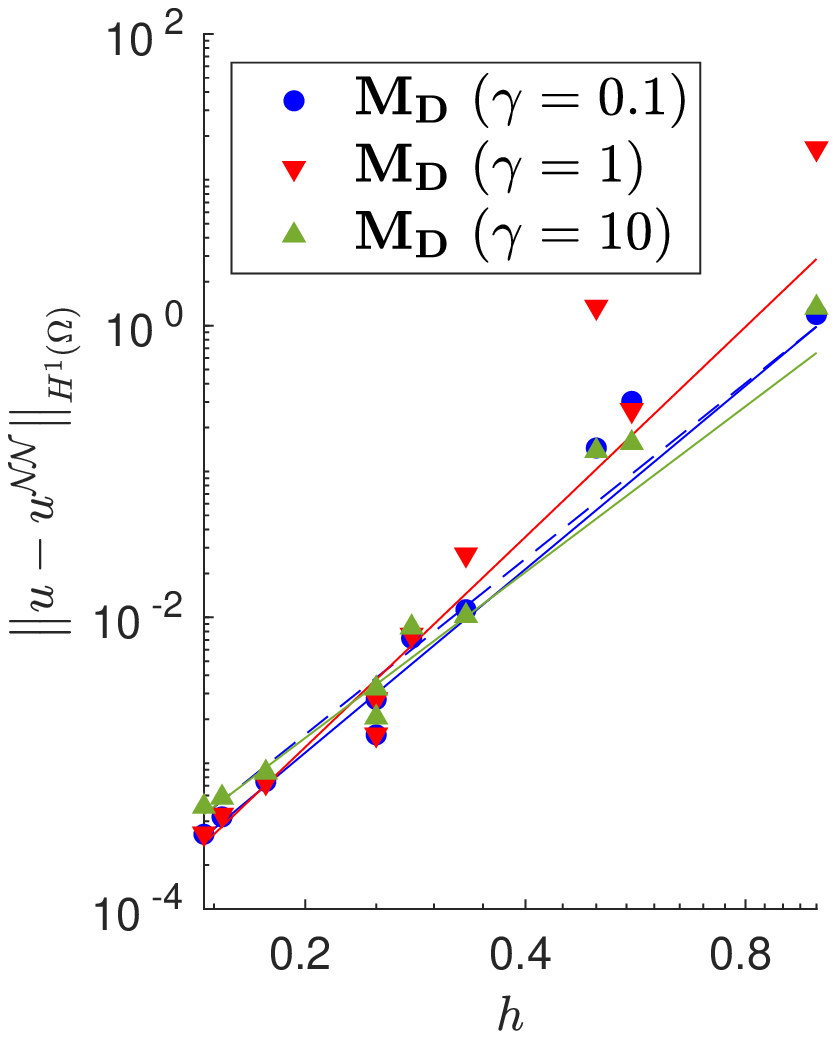} 
   \subcaption{$q=3$, $k_{\text{test}}=1$, $k_\text{int}=4$.}
  \label{fig:eq_2_hole_reg_hard_md_a}
\end{subfigure}
\begin{subfigure}[t]{0.32\linewidth}
  \includegraphics[width=0.98\columnwidth,keepaspectratio,clip]{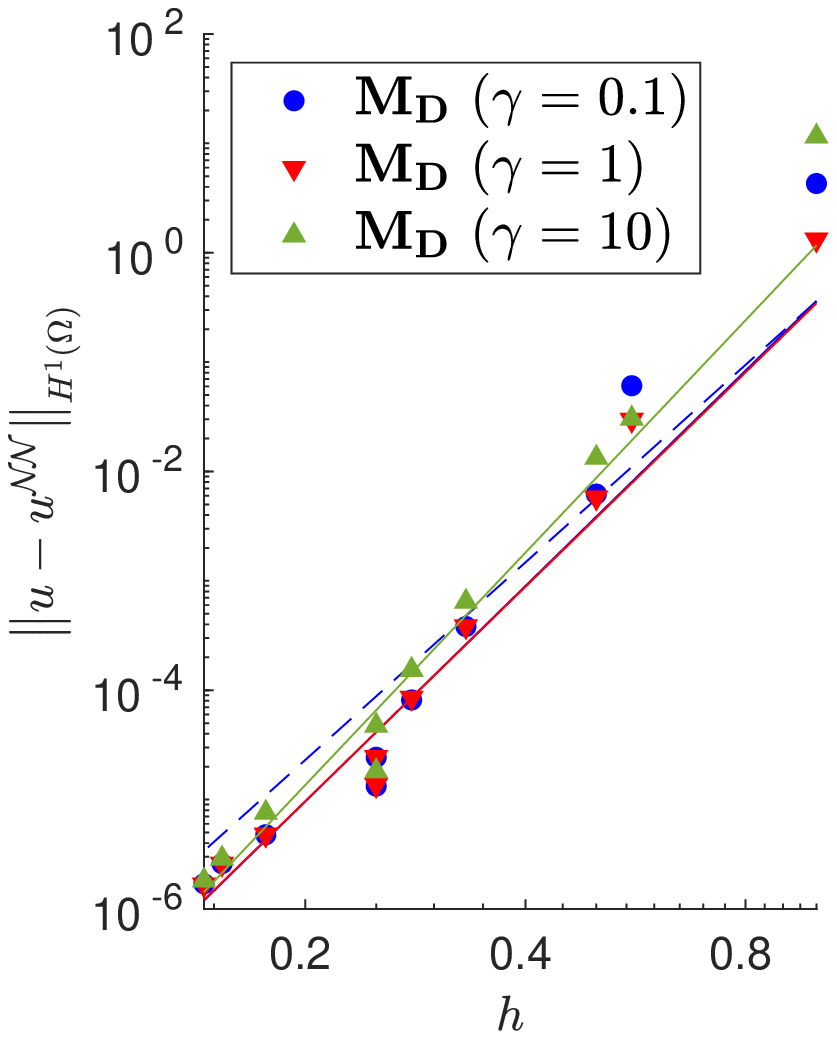} 
   \subcaption{$q=5$, $k_{\text{test}}=1$, $k_\text{int}=6$.}
  \label{fig:eq_2_hole_reg_hard_md_b}
\end{subfigure}
\begin{subfigure}[t]{0.32\linewidth}
  \includegraphics[width=0.98\columnwidth,keepaspectratio,clip]{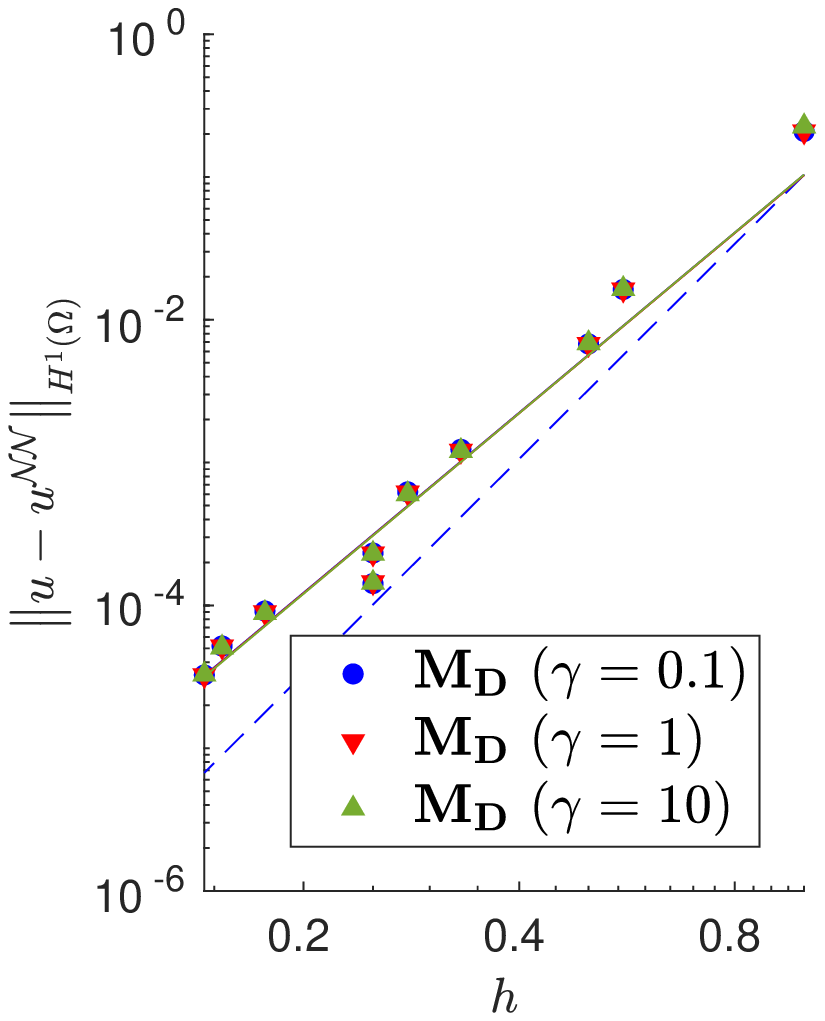}
   \subcaption{$q=5$, $k_{\text{test}}=2$, $k_\text{int}=5$.}
  \label{fig:eq_2_hole_reg_hard_md_c}
\end{subfigure}
  \caption{Error decay obtained with $\mathbf{M_D}$ and different values of $\gamma$. Forcing term and Dirichlet boundary conditions are set such that the exact solution is \eqref{eq:sol5}. The theoretical convergence rate is $k_\text{int}$. (a) Convergence rates: 4.18 $(\gamma=0.1)$, 4.79 $(\gamma=1)$, 3.78 $(\gamma=10)$. (b) Convergence rates: 6.54 $(\gamma=0.1)$, 6.51 $(\gamma=1)$, 7.06 $(\gamma=10)$. (c) Convergence rates: 4.19 $(\gamma=0.1)$, 4.19 $(\gamma=1)$, 4.20 $(\gamma=10)$. }
  \label{fig:eq_2_hole_reg_hard_md}
\end{figure}

Note that, if we change the forcing term and Dirichlet boundary conditions to consider the more oscillatory exact solution in \eqref{eq:sol5}, some approaches does not ensure the theoretical convergence rate (see Figures 
\ref{fig:eq_2_hole_reg_hard_ma_a}--\ref{fig:eq_2_hole_reg_hard_ma_c}, 
\ref{fig:eq_2_hole_reg_hard_mbmc_a}--\ref{fig:eq_2_hole_reg_hard_mbmc_c} and 
\ref{fig:eq_2_hole_reg_hard_md_a}--\ref{fig:eq_2_hole_reg_hard_md_c}). In fact, in Fig. \ref{fig:eq_2_hole_reg_hard_ma} it is evident that, in this case, large values of $\lambda$ are required to properly enforce the Dirichlet boundary conditions. In Fig. \ref{fig:eq_2_hole_reg_hard_mbmc}, instead, we can observe that the VPINN trained with method $\mathbf{M_C}$ is often inaccurate and the corresponding error decay is very noisy. The performance of methods $\mathbf{M_B}$ and $\mathbf{M_D}$ seems independent of the complexity of the forcing term and boundary conditions.

\begin{figure}[!h]
\centering 
\begin{subfigure}[t]{0.49\linewidth}
  \includegraphics[width=\columnwidth,keepaspectratio,clip]{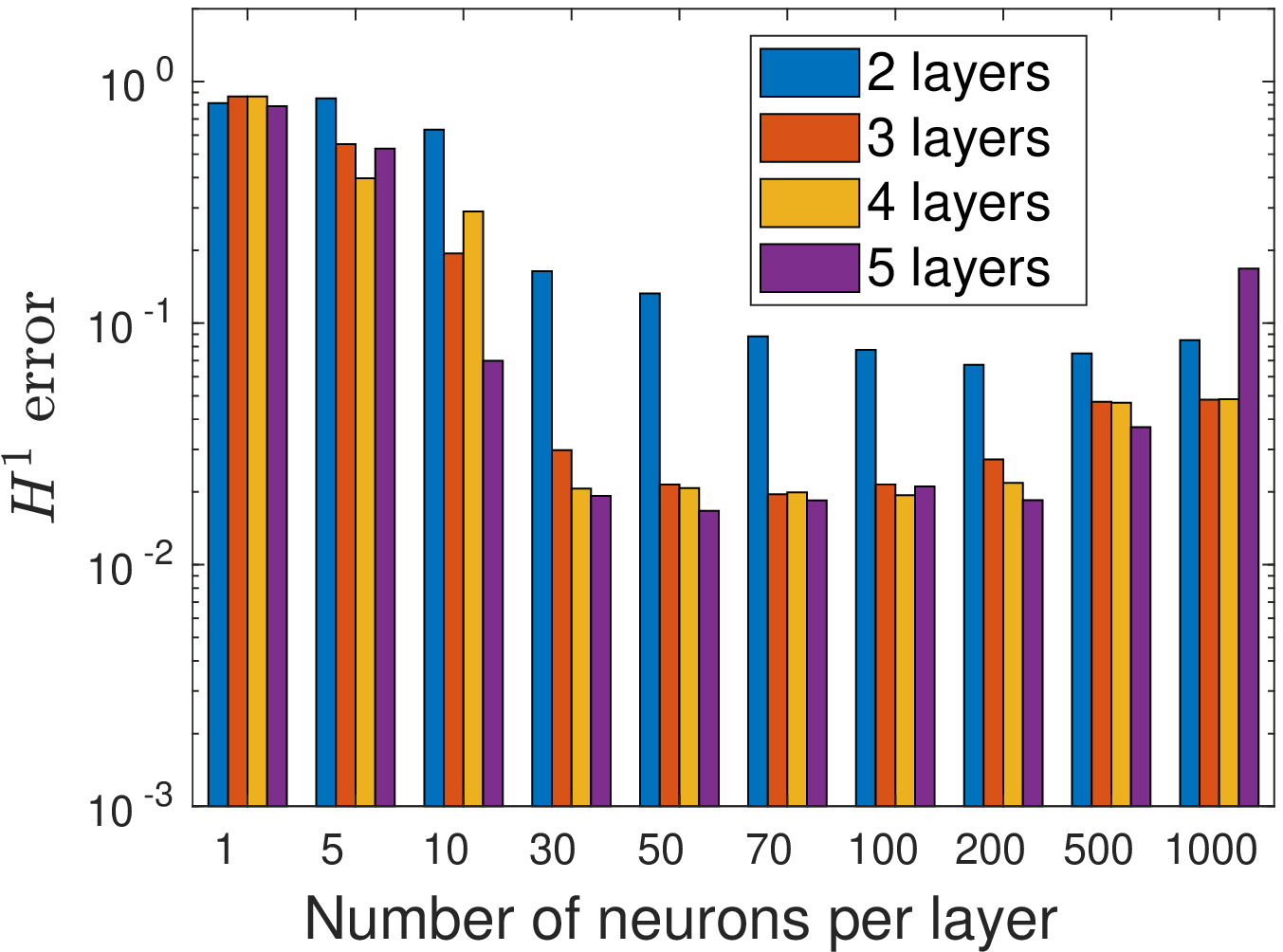} 
  \caption{VPINN error.}
  \label{fig:net_dims_vpinn}
\end{subfigure}
\begin{subfigure}[t]{0.49\linewidth}
  \includegraphics[width=\columnwidth,keepaspectratio,clip]{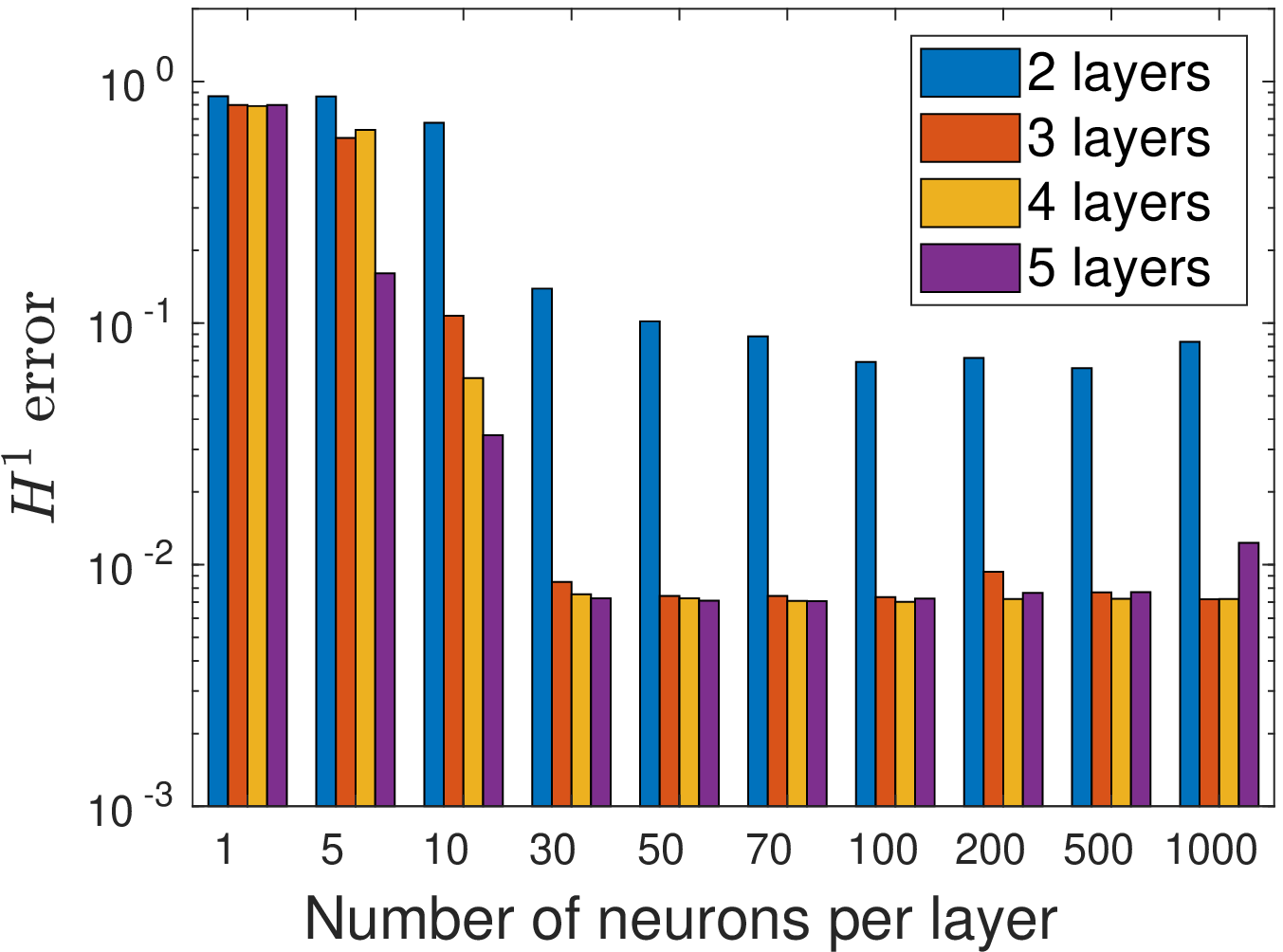} 
  \caption{Interpolated VPINN error.}
  \label{fig:net_dims_ivpinn}
\end{subfigure}
  \caption{$H^1$ errors using method $\bf M_B$ for
  standard (left) and interpolated VPINNs (right)
  as a function of the hyperparameters.}
  \label{fig:net_dims}
\end{figure}

In order to show that interpolation acts as a stabilization, we fix a mesh and vary the number of layers and neurons of the neural network. The boundary conditions are imposed using method $\mathbf{M_B}$ with $m=2$ and the exact solution is the one in \eqref{eq:sol2}; the results are shown in Fig.~\ref{fig:net_dims}. The number $L$ of layers varies in $\{2,3,4,5\}$, whereas the number of neurons in each hidden layer belongs to the set $\{1,5,10,30,50,70,100,200,500,1000\}$. In Fig.~\ref{fig:net_dims_vpinn} we show the performance of a non-interpolated VPINN trained with the $L^2$ regularization in \eqref{eq:l2_reg}, where $\lambda_{\text{reg}}=10^{-6}$. It can be noted that the error is high when the neural network is small because of its poor approximation 
capability, and that it decreases with intermediate values of the two hyperparameters. However, when the neural networks contain more than 100 neurons in each layer the error increases because of uncontrolled spurious 
zero-energy
modes and the fact that we are looking for good local minima in a very high-dimensional space. On the other hand, when the VPINN is interpolated and the neural network is sufficiently rich, the error is constant and independent of the network dimension (see Fig.~\ref{fig:net_dims_ivpinn}).
In addition, note that the average accuracy of an interpolated VPINN is better than its non-interpolated counterpart.

\subsection{Application to nonlinear parametric problem}\label{sect:nonlinear_parametric}
Let us now extend our analysis to nonlinear and parametric PDEs. Since in the previous section we observed that method $\mathbf{M_B}$ performs the best, in this example we 
do not consider $\mathbf{M_C}$ and $\mathbf{M_D}$. We focus on the following problem:
 \begin{equation}\label{eq:parametric-model-pb}
\left\{\begin{aligned}
N(u;p) := -\nabla \cdot (\mu \nabla u) + \boldsymbol{\beta}\cdot \nabla u + \sigma\sin(pu) u &= f \ \ \text{in } \Omega=(0,1)^2 , \\
u &= g \ \ \text{on } \Gamma_D .
\end{aligned}\right.
\end{equation}
It has been observed in \cite{berrone2022variational} that considering constant or variable coefficients does not influence VPINN convergence. Hence, we choose $\mu=1$, $\beta=[2,3]$, $\sigma=4$ and assume that the exact solution is
\begin{equation}\label{eq:sol8_param}
u(\bm{x};p) = \sin(p\pi x)\sin\left(\frac 1p\pi y\right),
\end{equation}
where $p\in\mathcal I_p = [0.5,2]$ is a scalar parameter.

In order to train the VPINN to solve problem \eqref{eq:parametric-model-pb}, we minimize
\begin{equation*}\label{eq:parametric-loss-function_ma}
R_h^2(w) = \sum_{p\in\mathcal I_p^\#}\left[\sum_{i \in I_h} r_{h,i;p}^2(w) + \lambda \sum_{i=1}^{N_B} \left( w(x_i^g)- g(x_i^g;p)\right)^2\right]
\end{equation*}
when $\mathbf{M_A}$ is used, or
\begin{equation*}\label{eq:parametric-loss-function_mb}
R_h^2(Bw) = \sum_{p\in\mathcal I_p^\#}\sum_{i \in I_h} r_{h,i;p}^2(Bw) 
\end{equation*}
when $\mathbf{M_B}$ is used instead. Here $\mathcal I_p^\#=\{p_1,\dots,p_{N_p^{\text{train}}}\}\subset\mathcal I_p$ is a finite set of parameter values and $ r_{h,i;p}$ is the residual obtained using the $i$-th test function and the parameter $p$. In the numerical computations, we use $N_p^{\text{train}}=13$ and the VPINN is trained with $q=3$ and $k_{\text{test}}=1$. 

\begin{figure}[!bht]
\centering 
\begin{subfigure}[t]{0.48\linewidth}
  \includegraphics[width=0.99\columnwidth,keepaspectratio,clip]{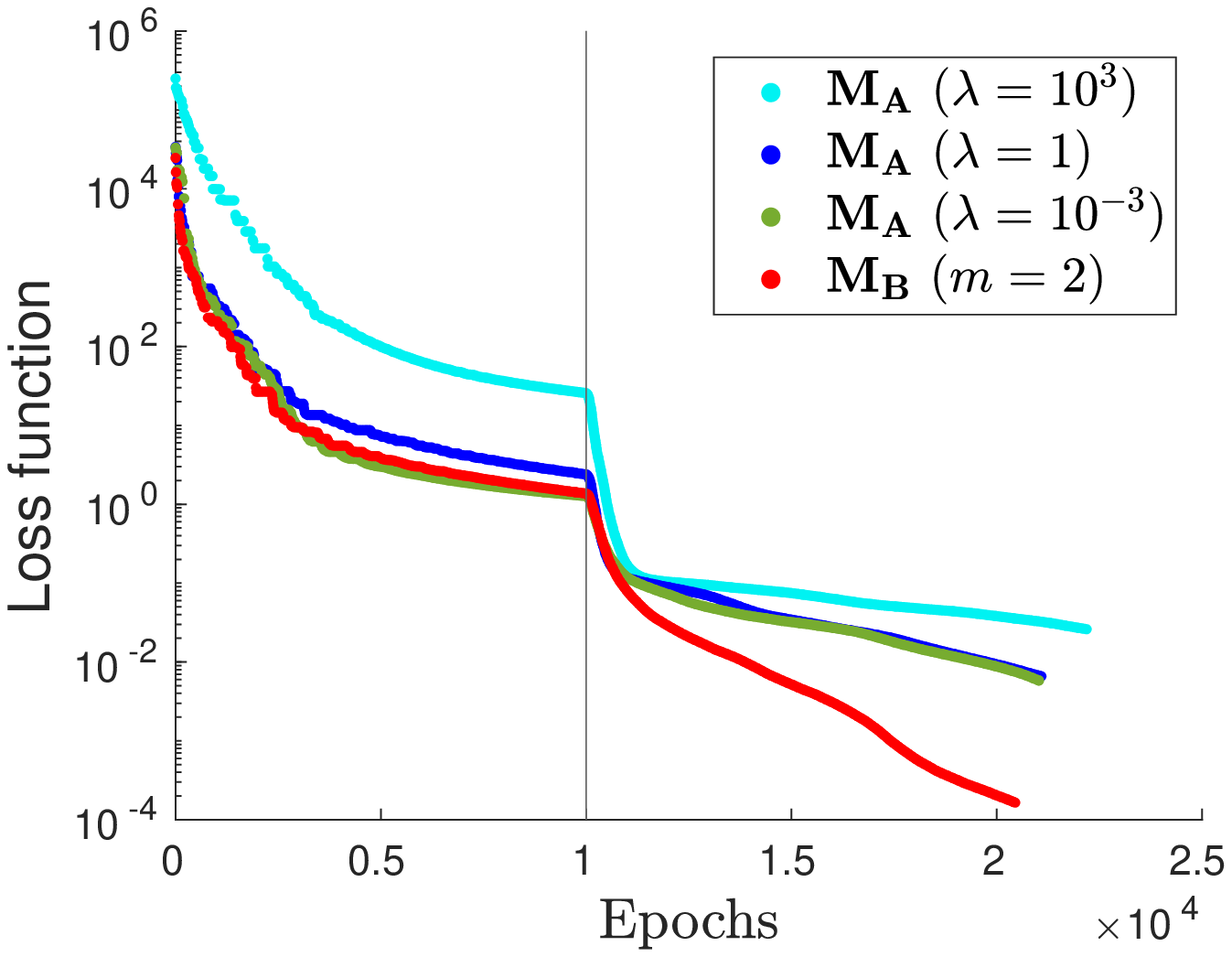} 
  \caption{}
  \label{fig:param_loss}
\end{subfigure}
\hspace{0.015\linewidth}
\begin{subfigure}[t]{0.48\linewidth}
  \includegraphics[width=0.99\columnwidth,keepaspectratio,clip]{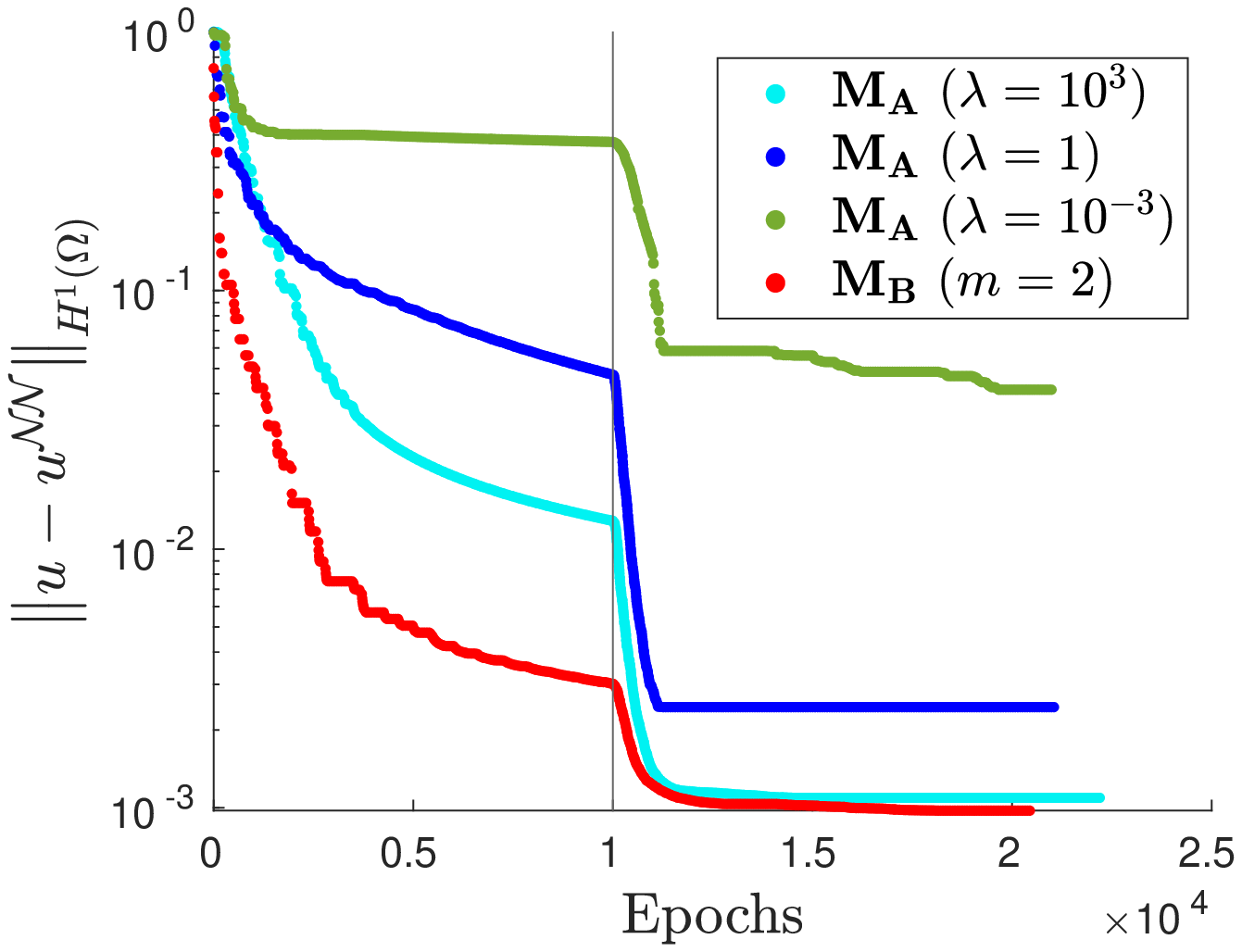} 
  \caption{}
  \label{fig:param_error}
\end{subfigure}
  \caption{(a) Training loss and (b) $H^1$ error
  prediction for the VPINN. 
  The first 10000 epochs are performed with a standard ADAM optimizer, the remaining ones with the BFGS optimizer. The exact solution is given in \eqref{eq:sol8_param}.}
  \label{fig:eq1param_sol8_loss_and_errors}
\end{figure}

In Fig.~\ref{fig:eq1param_sol8_loss_and_errors}, we report the behavior of the loss function and the average $H^1$ error:
 \[
 \frac{1}{N_p^{\text{test}}} \sum_{i=1}^{N_p^{\text{test}}} \Vert u(\cdot;p_i) - u^\NN(\cdot;p_i)\Vert_{H^1(\Omega)},
 \]
 where ${N_p^{\text{test}}}=100$. It is 
 noted that the loss functions behaves qualitatively 
 similarly (see Fig.~\ref{fig:param_loss}). On the other hand, higher values of $\lambda$ lead to lower errors when $\mathbf{M_A}$ is adopted, but the most stable and accurate approach remains $\mathbf{M_B}$.
 
\begin{figure}[!t]
\centering 
  \includegraphics[width=0.5\columnwidth,keepaspectratio,clip]{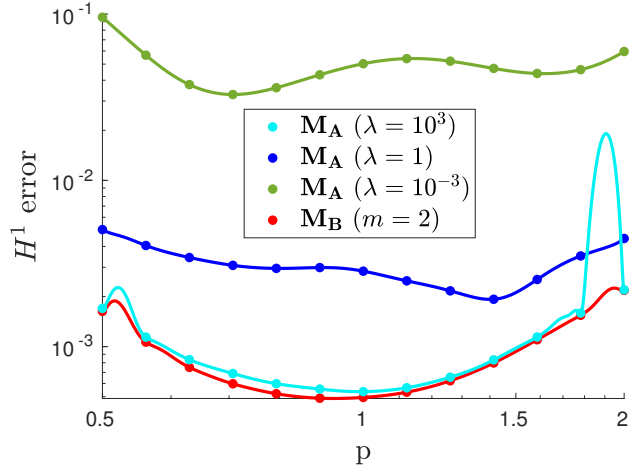} 
  \caption{$H^1$ error for different parameter values in problem~\eqref{eq:parametric-model-pb}.}
  \label{fig:global_param_error}
\end{figure}

 Moreover, when the VPINN is trained, it can be evaluated at arbitrary locations in the parameter domain $\mathcal I_p$, yielding the error plot shown in Fig.~\ref{fig:global_param_error}. Therein, for each dot and for each point belonging to the solid lines, given the parameter value $\hat p$ represented on the horizontal axis, its value on the vertical axis represents the $H^1$ error between the VPINN solution and the exact solution $u(\cdot;\hat p)$. Note that dots are associated with parameter values that are chosen in $\mathcal I_p^\#$ during the training, whereas solid lines are the predictions to assess the accuracy of the models for intermediate values of $p$. Such lines thus show the $H^1$ error for values of the parameter not used during the training.
 

\subsection{Deformation of an elastic body}\label{sect:elasticity}
We consider the deformation of a linear elastic solid in the region $\Omega_L=(-1,1)^2\setminus[-1,0]^2$, which is subjected to a body force field $\mathbf f$ and Dirichlet boundary conditions imposed on $\Gamma_D=\partial \Omega$. The elastostatic boundary-value problem is:
\begin{subnumcases}{\label{eq:le_system}}
-\nabla\cdot \bm\sigma =\mathbf f & \text{in \ } $\Omega_L$\,, \label{eq:elastic-model-pb}\\
\bm\varepsilon = \frac12\left(\nabla \mathbf u + \left(\nabla\mathbf u^T\right)\right)& \text{in \ } $\Omega_L$\,, \label{eq:symm_grad_u}\\
\bm\sigma = 2\mu \bm\varepsilon + \lambda \hspace{0.075cm}\text{trace}(\bm\varepsilon)\mathbf I& \text{in \ } $\Omega_L$\,, \label{eq:hooke}\\
\mathbf u=\mathbf g & \text{on \ }$ \Gamma_D$ \,.\label{eq:el_bc}
\end{subnumcases}

In \eqref{eq:le_system}, $\bm\sigma:=\bm\sigma(\mathbf u)$ is the Cauchy stress tensor, $\bm\varepsilon := \varepsilon(\mathbf u)$ is the small strain tensor and~\eqref{eq:hooke} is the isotropic linear elastic constitutive relation.
The Lam\'e parameters $\lambda$ and $\mu$ are related to the Young modulus $E$ and the Poisson ratio $\nu$ via
\[
\mu = \frac{E}{2(1+\nu)}, \quad 
\lambda = \frac{E\nu}{(1+\nu)(1-2\nu)}.
\]
For the numerical experiments, we choose $E=117$, $\nu=1/3$ and the following body force field and boundary data:

\[
\mathbf f = (\mu+\lambda)\left[ x e^y, \,\,y\sqrt{x+2}\right], \hspace{1cm}
\mathbf g = \left[ \sin(\pi(x+y)), \,\,e^{x-y}\right]xy.
\]
The variational formulation of problem \eqref{eq:elastic-model-pb} reads as: Find $\mathbf u\in\overline{\mathbf u}+ \left(H_0^1(\Omega)\right)^2$ such that:
\begin{equation*}\label{eq:variational_elasticity}
\int_\Omega \bm{\sigma}(\mathbf u):\bm{\varepsilon}(\mathbf v) = \int_\Omega \mathbf f\mathbf v  
\ \ \forall \mathbf v\in\left(H_0^1(\Omega)\right)^2,
\end{equation*}
where $\overline{\mathbf{u}}=\mathbf g$ is the natural lifting of the boundary data. Such a formulation is used to compute the quantity $R_h^2$ in \eqref{eq:var_loss}, where the modified residuals
\[
r_{h,i}(\mathbf w) = \int_\Omega \mathbf f\bm{ \varphi}_i^v - \int_\Omega \bm\sigma(\mathbf w):\bm\varepsilon(\bm{\varphi}_i^v ), 
\quad i\in I_h,
\]
replace the ones defined in \eqref{eq:residuals}. For this and the subsequent test cases, we will also provide a comparison with the results obtained by a PINN, in order to give a more complete view of the performance of the methods. The modified residuals required in the PINN loss function are defined as:
\[
r_i^\text{PINN}(\mathbf u) = \nabla\cdot \bm\sigma(x_i) +\mathbf f(x_i) 
\ \ \forall i=1,2,\dots,N_I.
\]

\begin{figure}[t!]
\centering 
\begin{subfigure}[t]{0.49\linewidth}
 \centering
 \includegraphics[width=\columnwidth,keepaspectratio,clip]{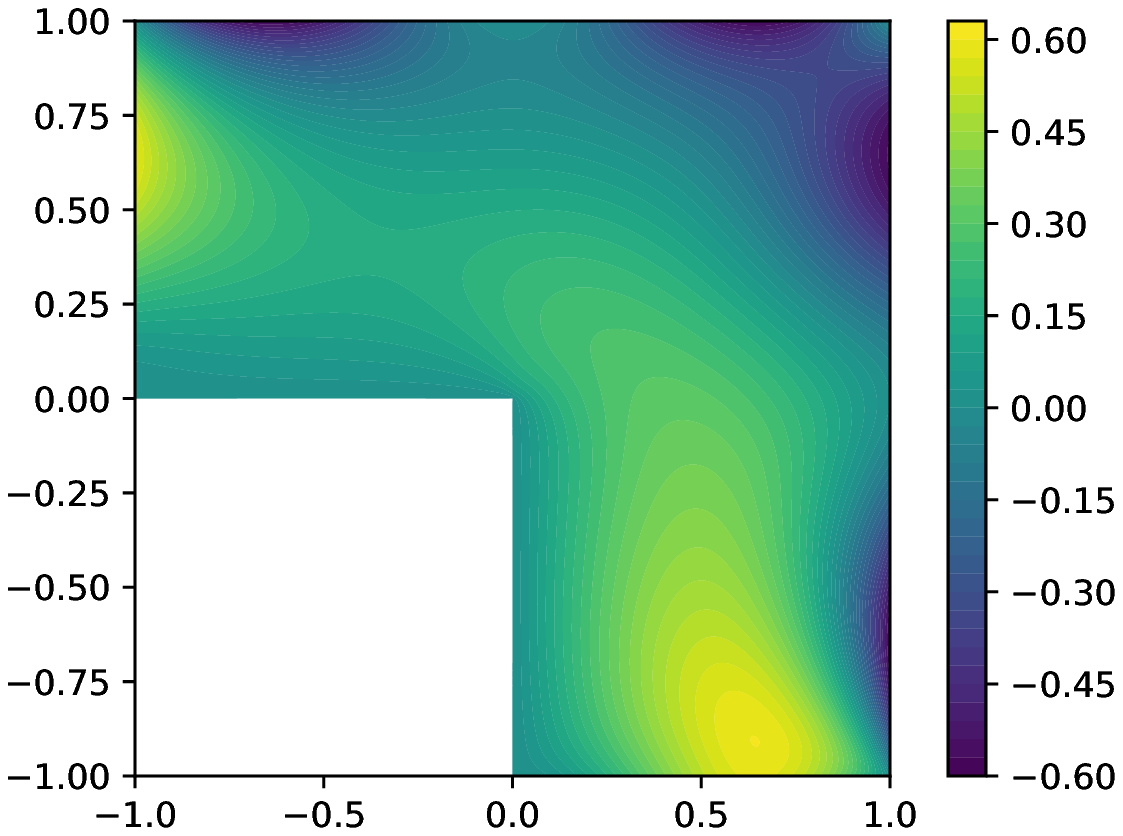}
 \caption{$x$-component of $\mathbf{u}^h$}
 \label{fig:el_sol_ux}
\end{subfigure}
\begin{subfigure}[t]{0.49\linewidth}
  \centering
  \includegraphics[width=0.95\columnwidth,keepaspectratio,clip]{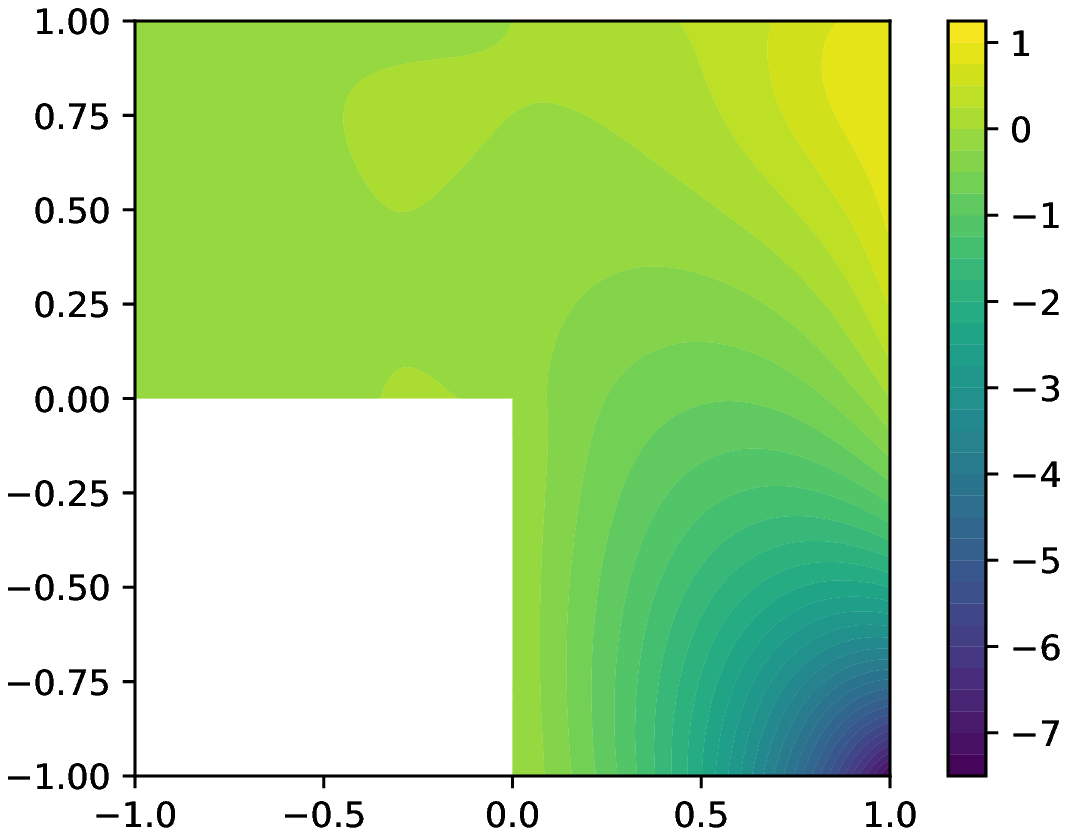} 
 \caption{$y$-component of $\mathbf{u}^h$}
 \label{fig:el_sol_uy}
\end{subfigure}
  \caption{Reference finite element displacement field solution 
  for problem~\eqref{eq:le_system}. The $x$-component
  (left) and $y$-component (right) of $\mathbf{u}^h$ are 
  shown.}
  \label{fig:el_sol}
\end{figure}

Since the exact solution is not known, we produce a very accurate numerical solution for comparison (shown in Figures~\ref{fig:el_sol_ux} and~\ref{fig:el_sol_uy}), using the open-source FEM solver FEniCS~\cite{AlnaesEtal2015}.


Problem \eqref{eq:le_system} is solved by training a VPINN on the mesh shown in Fig.~\ref{fig:l-mesh} with $q=3$ and $k_{\text{test}}=1$. Then, in order to compare the performance of PINN and VPINN, a standard PINN is trained to solve the same problem. In order to verify if the distribution of the collocation points affects the PINN accuracy, we firstly train it by choosing as collocation points the interpolation nodes used in the VPINN training, and then we train it with the same number of uniformly distributed collocation points.

\begin{figure}[t!]
\centering 
\begin{subfigure}[t]{0.8\linewidth}
  \includegraphics[width=0.99\columnwidth,keepaspectratio,clip]{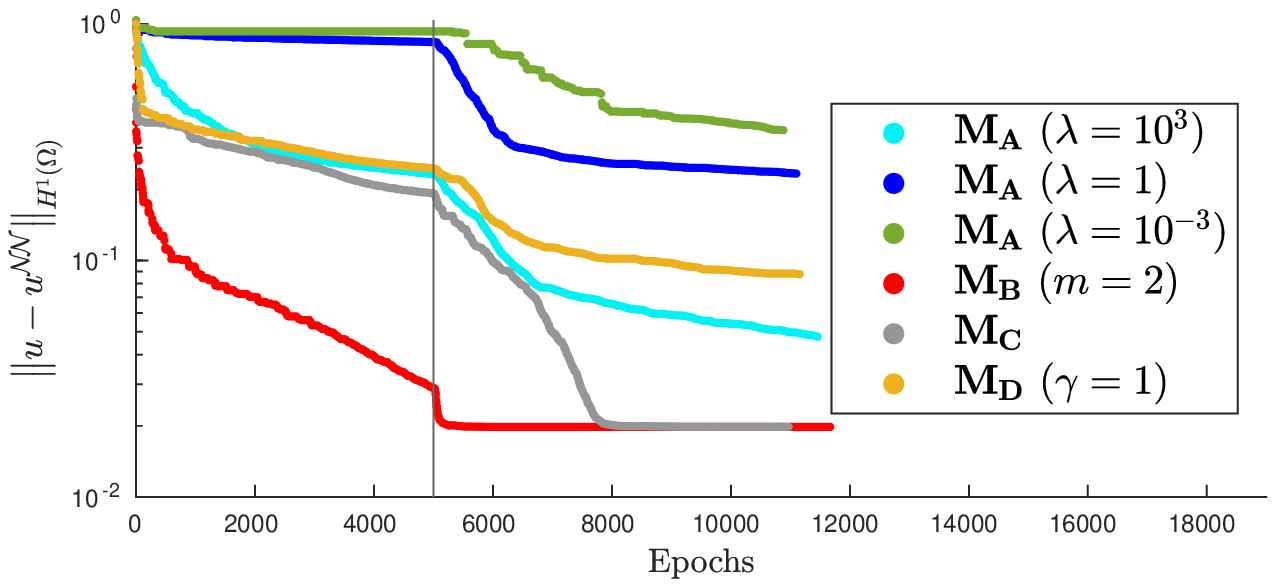} 
  \caption{}
  \label{fig:le_vpinn_errors}
\end{subfigure}

\begin{subfigure}[t]{0.48\linewidth}
  \includegraphics[width=0.99\columnwidth,keepaspectratio,clip]{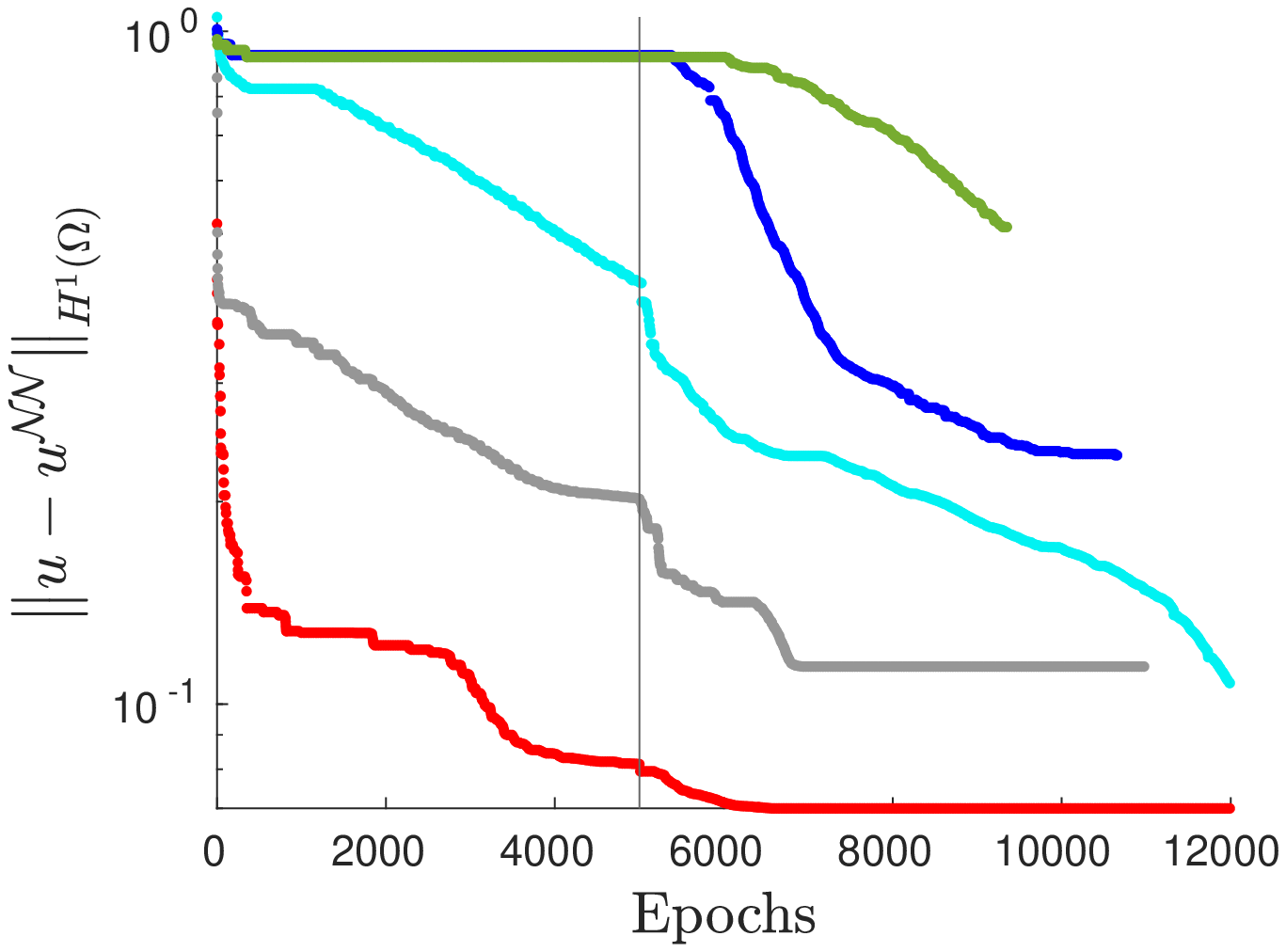} 
  \caption{}
  \label{fig:le_pinn_errors_intppts}
\end{subfigure}
\hspace{0.015\linewidth}
\begin{subfigure}[t]{0.48\linewidth}
  \includegraphics[width=0.99\columnwidth,keepaspectratio,clip]{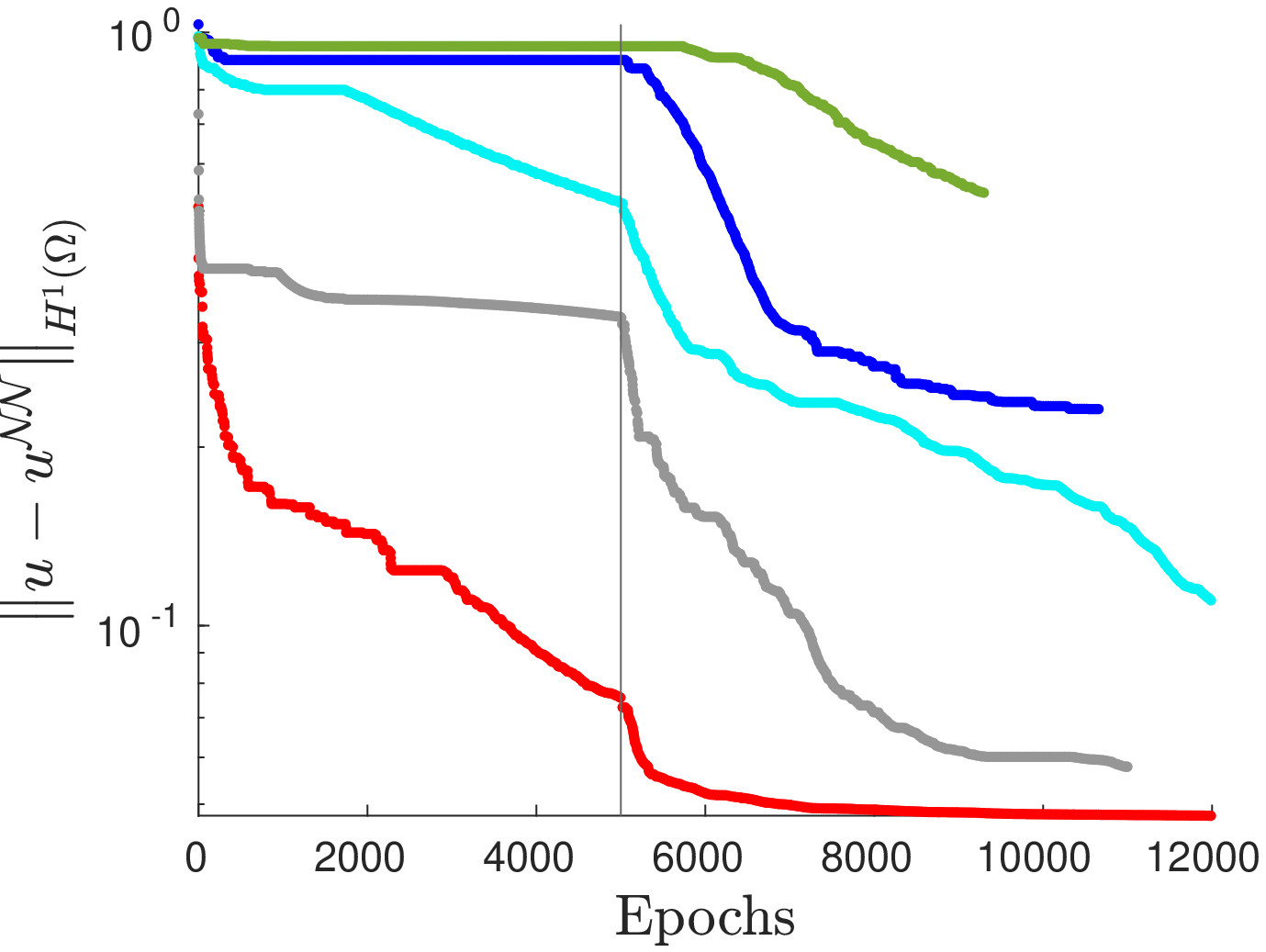} 
  \caption{}
  \label{fig:le_pinn_errors_rndpts}
\end{subfigure}
  \caption{$H^1$ error decay during the neural network training when solving 
  problem~\eqref{eq:le_system}. (a) VPINN error: $H^1$ error of the most accurate solution is 0.020; (b) PINN error: model is trained with collocation points distributed on a Delaunay mesh and the $H^1$ error of the most accurate solution is 0.070; and 
  (c) PINN error: model is trained with 
  collocation points from a uniform distribution and the $H^1$ error of the most accurate solution is 0.047. The legend in (a) also applies to (b) and (c). }
  \label{fig:le_errors}
\end{figure}

For these three methods we analyze the $H^1$ error during the neural network training for a fixed training set dimension; we report the results in Figures {\ref{fig:le_vpinn_errors}}--{\ref{fig:le_pinn_errors_rndpts}}. Observing that Figures {\ref{fig:le_pinn_errors_intppts}} and {\ref{fig:le_pinn_errors_rndpts}} are very similar, we deduce that, in this case, the choice of control points in the PINN training is not strictly related to the efficacy of the different approaches.

It can be observed that method $\mathbf{M_B}$ is always the most efficient approach and leads to convergence to more accurate solutions. Exactly imposing the Dirichlet boundary conditions via $\mathbf{M_C}$ can be considered a good alternative since the solutions at convergence obtained with the VPINN and the PINN trained with random control points are very similar to the ones computed using $\mathbf{M_B}$, although the convergence is slower. The most commonly used approach, $\mathbf{M_A}$, is instead dependent on the choice of the non-trainable parameter $\lambda$. In this case, large values of $\lambda$ ensure accurate solutions and acceptably efficient training phases, but the correct values are problem dependent and can be often found only after a potentially expensive tuning. Indeed, choosing the wrong values of $\lambda$ can ruin the efficiency and the accuracy of the method, as it can be observed in Fig.~\ref{fig:le_errors} when $\lambda=10^{-3}$ or $\lambda=1$. We also highlight that the performance of method $\mathbf{M_D}$ is very similar to method $\mathbf{M_A}$ when reasonable values of $\lambda$ are chosen.
 
\subsection{Stabilized Eikonal equation}\label{sect:eikonal}
In this section we consider the stabilized Eikonal equation, which is a nonlinear second-order PDE and reads as:
 \begin{equation*}\label{eq:eikonal-model-pb}
\begin{cases}
-\varepsilon\Delta u + \Vert \nabla u\Vert_2 = f & \text{in } \Omega_L , \\
u = g & \text{on } \Gamma_D ,
\end{cases}
\end{equation*}
where $\varepsilon$ is a small positive constant. Note that when $\varepsilon=0$, $f=1$ and $g=0$, the exact solution  is the distance function to the boundary and the problem can be efficiently solved by the fast sweeping method \cite{zhao2005fast} or by the fast marching method \cite{sethian1999level}. In our numerical computations we set $f=1$ and $g=0$ and we introduce a weak diffusivity with $\varepsilon=0.1$ to guarantee uniqueness of the solution.

The PINN and VPINN residuals associated with problem \eqref{eq:eikonal-model-pb} that extend the residuals in \eqref{eq:pinn_residuals} and \eqref{eq:residuals}, respectively, are defined as:
\[
r_i^\text{PINN}(w) =-\varepsilon\Delta w(x_i) + \Vert \nabla w(x_i)\Vert_2 -f(x_i) 
\ \ \forall i=1,\dots,N_I,
\]
and 
\[
r_{h,i}(w) = \int_\Omega f\varphi_i^v - \int_\Omega \varepsilon \nabla w \nabla \varphi_i^v - \int_\Omega\Vert\nabla w\Vert_2 \varphi_i^v, 
\quad i\in I_h.
\]

\begin{figure}[t!]
\centering 
  \includegraphics[width=0.65\columnwidth,keepaspectratio,clip]{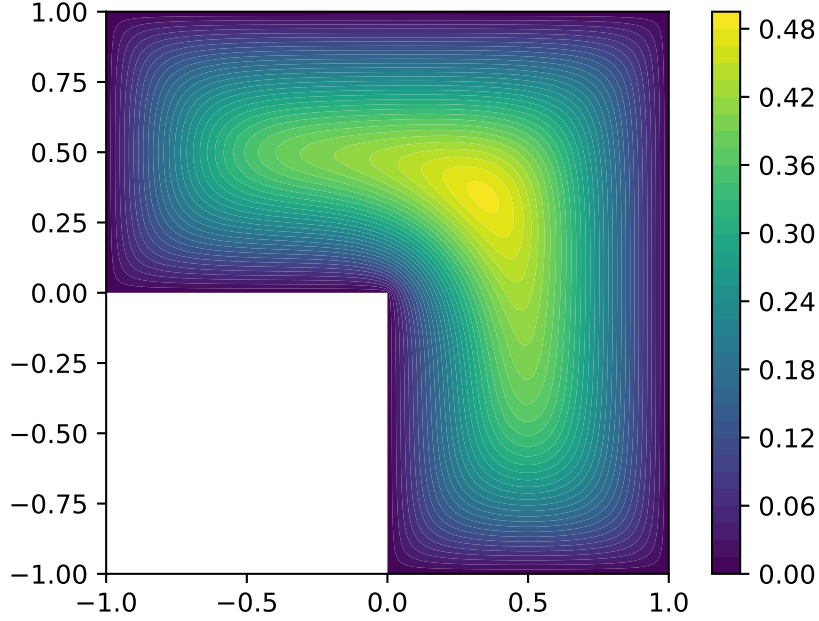} 
  \caption{Reference finite element solution for
  problem~\eqref{eq:eikonal-model-pb}.}
  \label{fig:eik_sol}
\end{figure}

We compute the VPINN and PINN numerical solutions as described in Section~\ref{sect:elasticity} and compute the corresponding $H^1$ errors using a finite element reference solution that is 
computed on a much finer 
mesh (see Fig.~\ref{fig:eik_sol}).

\begin{figure}[t!]
\centering 
\begin{subfigure}[t]{0.89\linewidth}
  \includegraphics[width=0.99\columnwidth,keepaspectratio,clip]{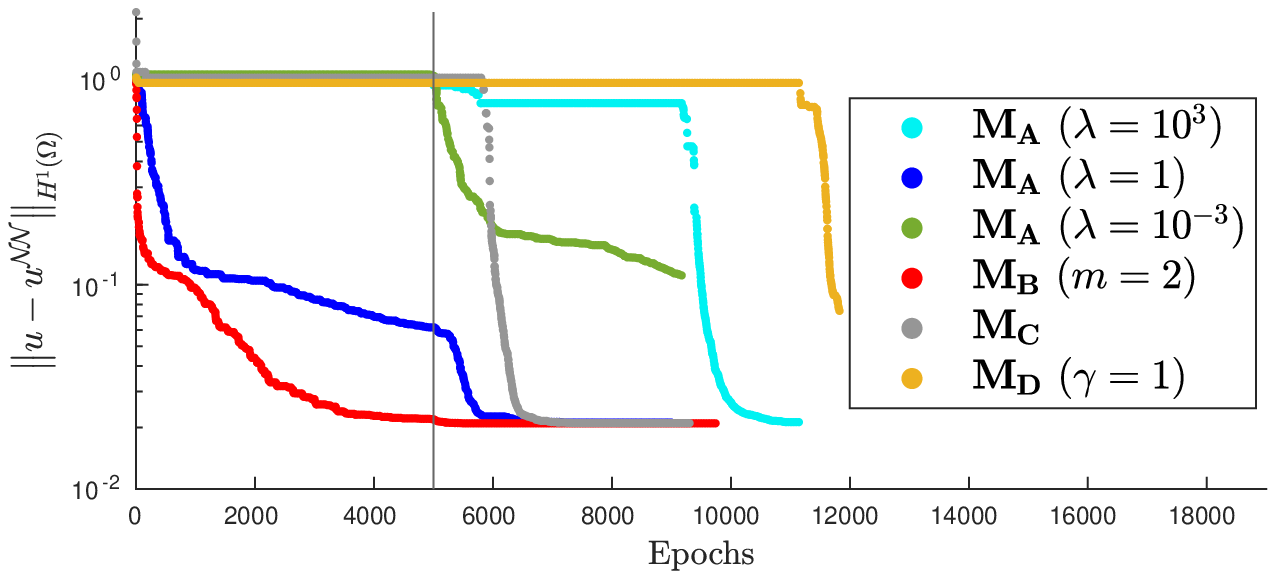} 
  \caption{}
  \label{fig:eik_vpinn_errors}
\end{subfigure}

\begin{subfigure}[t]{0.48\linewidth}
  \includegraphics[width=0.99\columnwidth,keepaspectratio,clip]{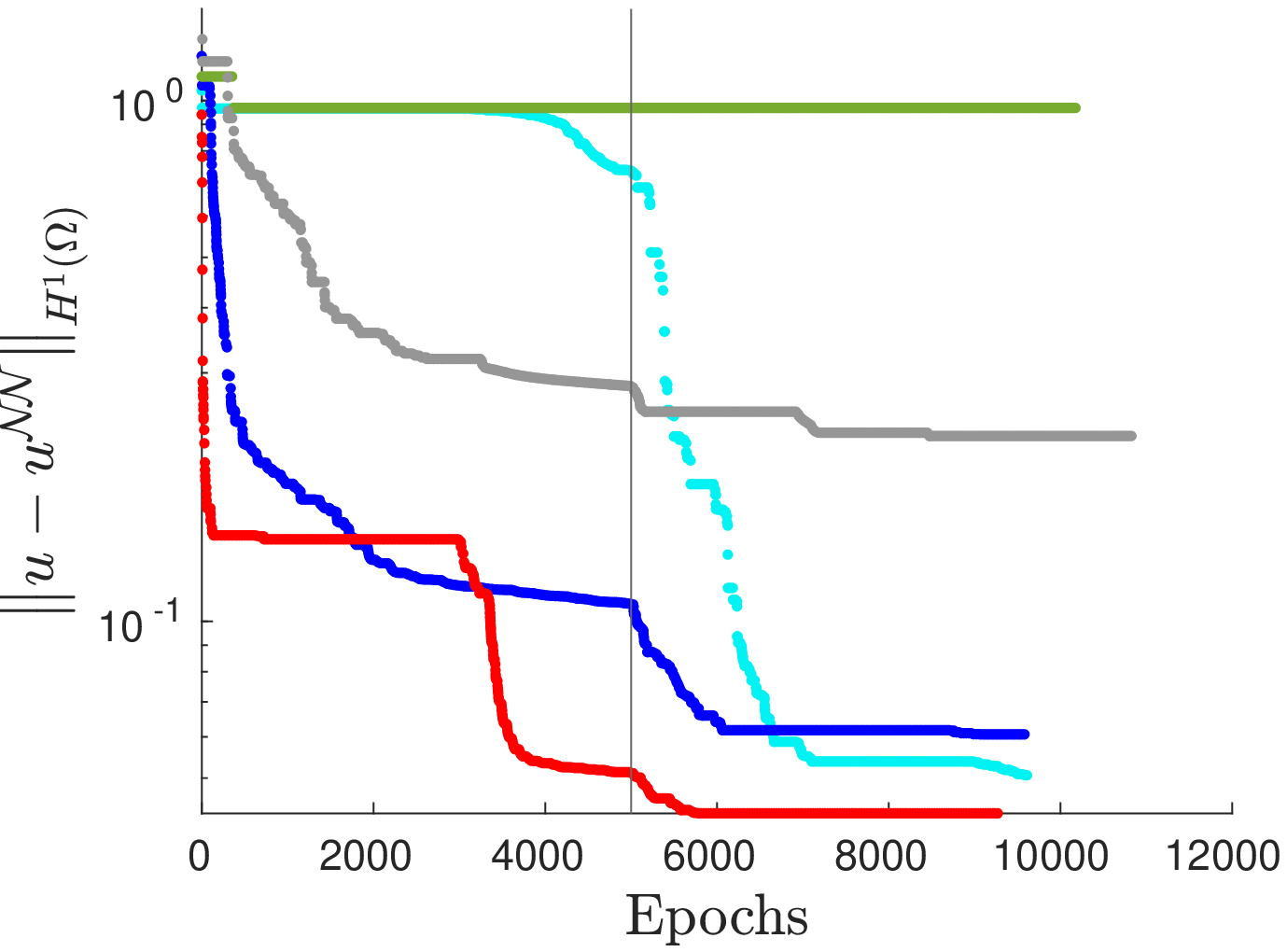} 
  \caption{}
  \label{fig:eik_pinn_errors_intppts}
\end{subfigure}
\hspace{0.015\linewidth}
\begin{subfigure}[t]{0.48\linewidth}
  \includegraphics[width=0.99\columnwidth,keepaspectratio,clip]{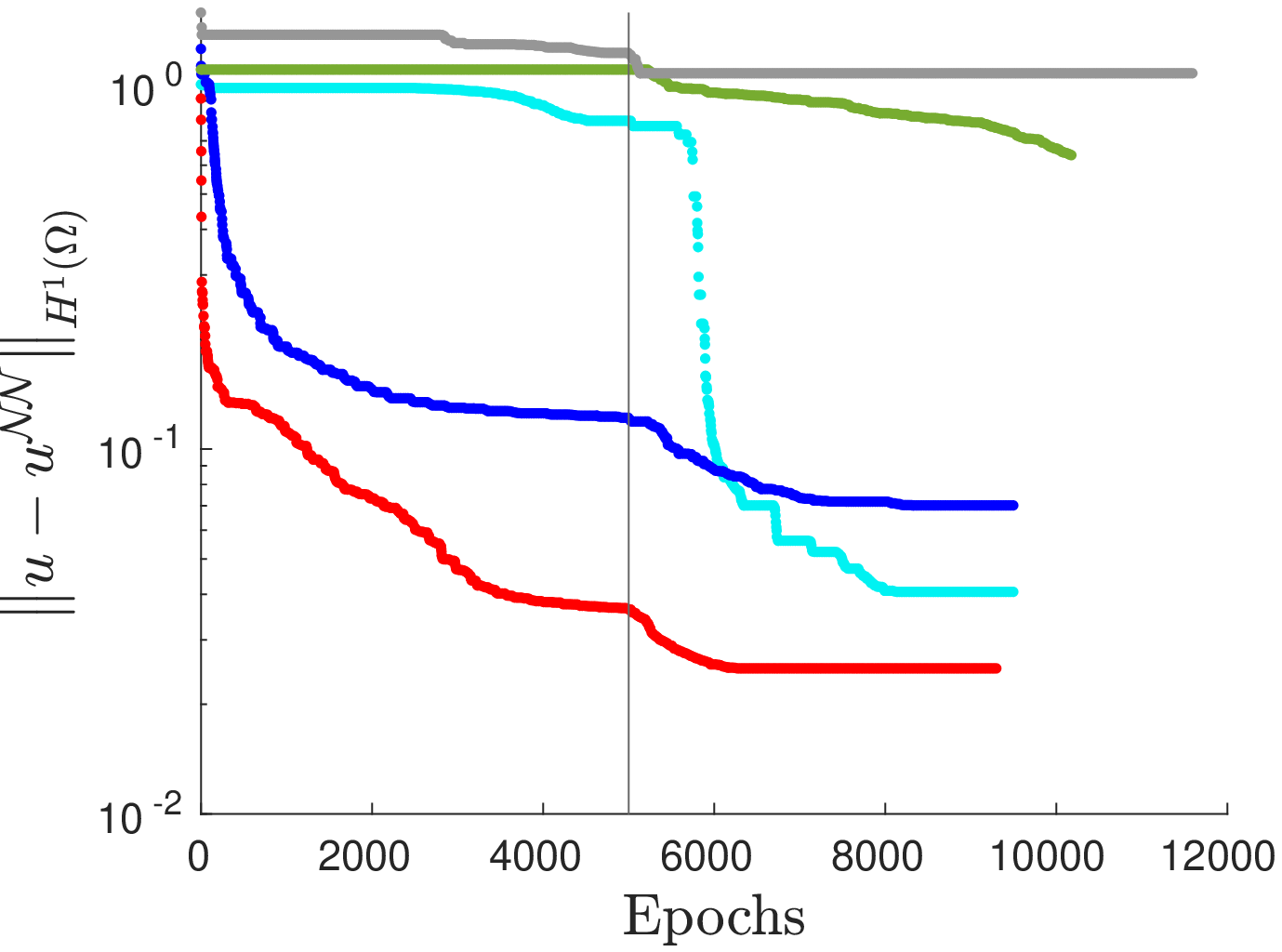} 
  \caption{}
  \label{fig:eik_pinn_errors_rndpts}
\end{subfigure}
  \caption{$H^1$ error decay during the neural network training when solving problem \eqref{eq:eikonal-model-pb}. (a) VPINN error: $H^1$ error of the most accurate solution is 0.021; (b) PINN error: model is trained with collocation points distributed on a Delaunay mesh
  and the $H^1$ error of the most accurate solution is 0.085; and (c) PINN error: model is trained with 
  collocation points from a uniform distribution and the $H^1$ error of the most accurate solution is 0.029. The legend in (a) also applies to (b) and (c).}
  \label{fig:eik_errors}
\end{figure}

As in Fig.~\ref{fig:le_errors}, in Fig.~\ref{fig:eik_errors} we show the decay of the $H^1$ error during the training for the different methods.
Again, it can be noted that the most accurate method is always $\mathbf{M_B}$; $\mathbf{M_A}$ is a valid alternative provided $\lambda$ is properly chosen. However, when the value of $\lambda$ is not suitably chosen, convergence can be completely ruined (see, for instance, the curves associated with $\lambda=10^{-3}$ in Figures {\ref{fig:eik_pinn_errors_intppts}} and  {\ref{fig:eik_pinn_errors_rndpts}}) or a second-order optimizer is required to retain convergence (see all the curves computed with $\lambda=10^3$). Moreover, similar convergence issues are present when $\mathbf{M_C}$ or $\mathbf{M_D}$ are employed.

\subsection{One-dimensional convection problem}\label{sect:transport}
As a final example,
we consider a one-dimensional convection problem on the space-time domain $\Omega:=\Omega_x\times\Omega_t=[0,1]\times[0,1]$. As discussed in \cite{krishnapriyan2021characterizing}, when solving such a hyperbolic PDE with PINN, possible failure modes may arise due to the very complex loss landscape. The 
model problem reads as: 
\begin{equation}\label{eq:transport-model-pb}
\begin{cases}
\dfrac{\partial u}{\partial t} + \beta \dfrac{\partial u}{\partial x} = 0, & \forall x\in\Omega_x=[0,1], \,
t\in\Omega_t=[0,1] , \\
u(0,t)=g(t) & \forall t\in\Omega_t , \\
u(x,0)=h(x) & \forall x\in\Omega_x .
\end{cases}
\end{equation}
Let us consider the boundary condition $g(t)=-\sin(\beta t)$ and the initial condition $h(x)=\sin(x)$. The corresponding exact solution is $u(x,t)=\sin(x-\beta t)$. We solve problem \eqref{eq:transport-model-pb} with the convection coefficient $\beta=30$.

Given a set of collocation points $(x_i,t_i)\in\Omega$, $i=1,\dots,N_I$ and a suitable set of space-time test functions $V_h:=\text{span}\{\varphi_i^v=\varphi_i^v(x,t):i\in I_h\}$, the PINN and VPINN residuals that are used to train the models are given by
\[
r_i^\text{PINN}(w) =\dfrac{\partial w}{\partial t}(x_i,t_i) + \beta \dfrac{\partial w}{\partial x}(x_i,t_i) 
\ \ \forall i=1,2,\dots,N_I
\]
and 
\[
r_{h,i}(w) = \int_\Omega \left[\dfrac{\partial w}{\partial t} + \beta \dfrac{\partial w}{\partial x}\right]  \varphi_i^v, 
\quad i\in I_h,
\]
respectively.
When the boundary conditions are exactly imposed (i.e., when $\mathbf{M_B}$ or $\mathbf{M_C}$ are used), the function $\phi=\phi(x,t)$ is constructed as $\phi(x,t):=\phi_x(x)\phi_t(t)$, where $\phi_t(t)=t$ and $\phi_x(x)$ is a function that vanishes on the Dirichlet boundary of $\Omega_x$. Note that, due to the simplicity of the spatial domain $\Omega_x$, there is no reason to distinguish between $\mathbf{M_B}$ and $\mathbf{M_C}$. Therefore, we just consider the function $\phi_x(x)=x$ in both approaches.

\begin{figure}[t!]
\centering 
\begin{subfigure}[t]{0.89\linewidth}
  \includegraphics[width=0.99\columnwidth,keepaspectratio,clip]{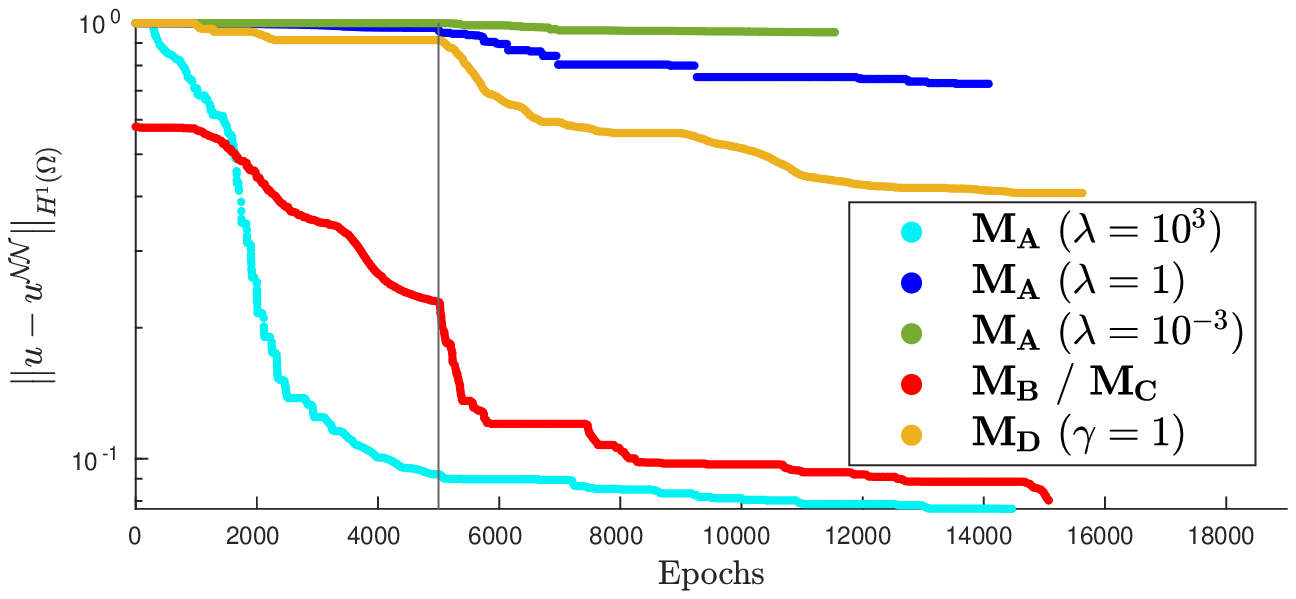} 
  \caption{}
  \label{fig:transport_vpinn_errors}
\end{subfigure}

\begin{subfigure}[t]{0.48\linewidth}
  \includegraphics[width=0.99\columnwidth,keepaspectratio,clip]{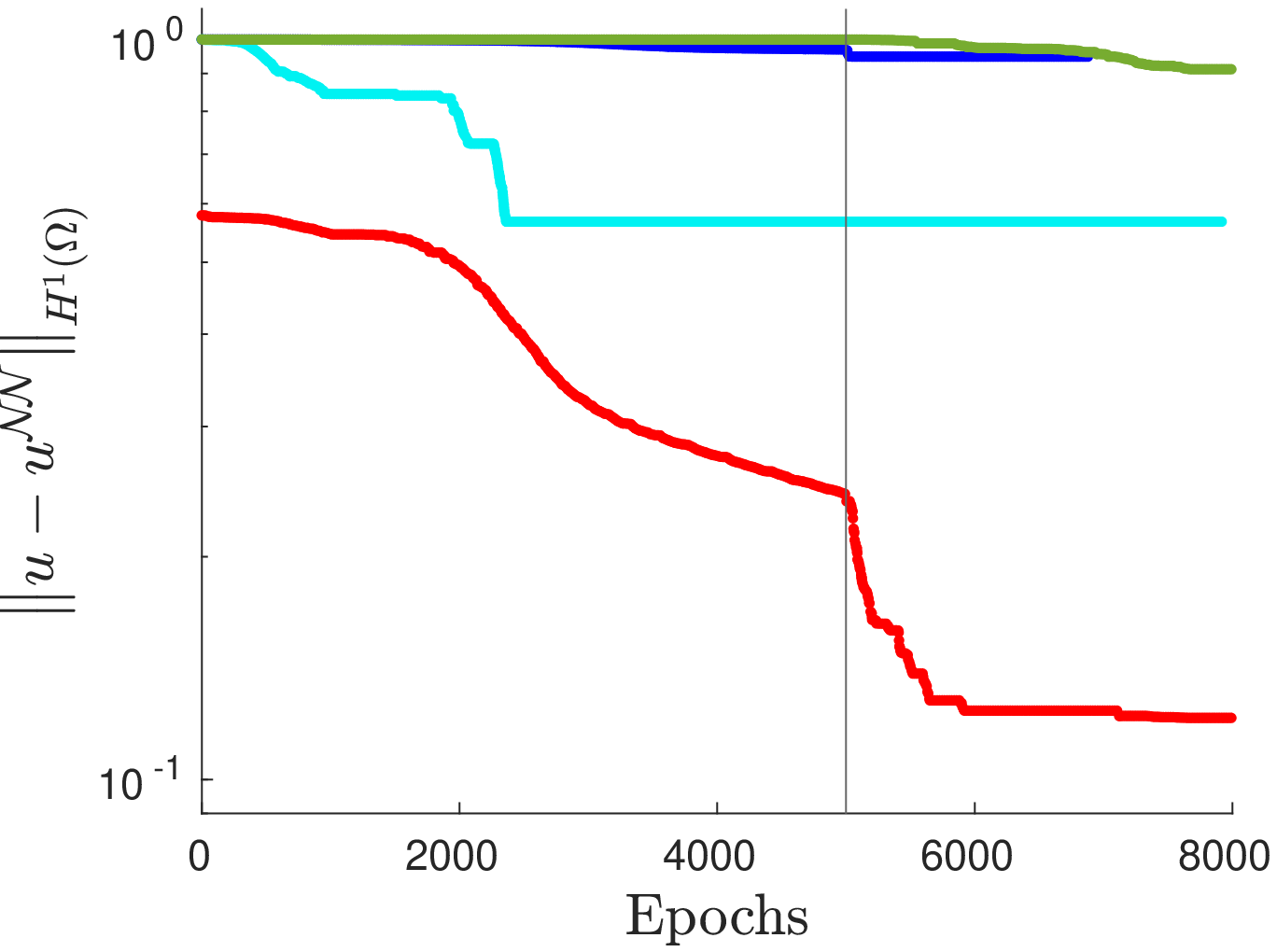} 
  \caption{}
  \label{fig:transport_pinn_errors_intppts}
\end{subfigure}
\hspace{0.015\linewidth}
\begin{subfigure}[t]{0.48\linewidth}
  \includegraphics[width=0.99\columnwidth,keepaspectratio,clip]{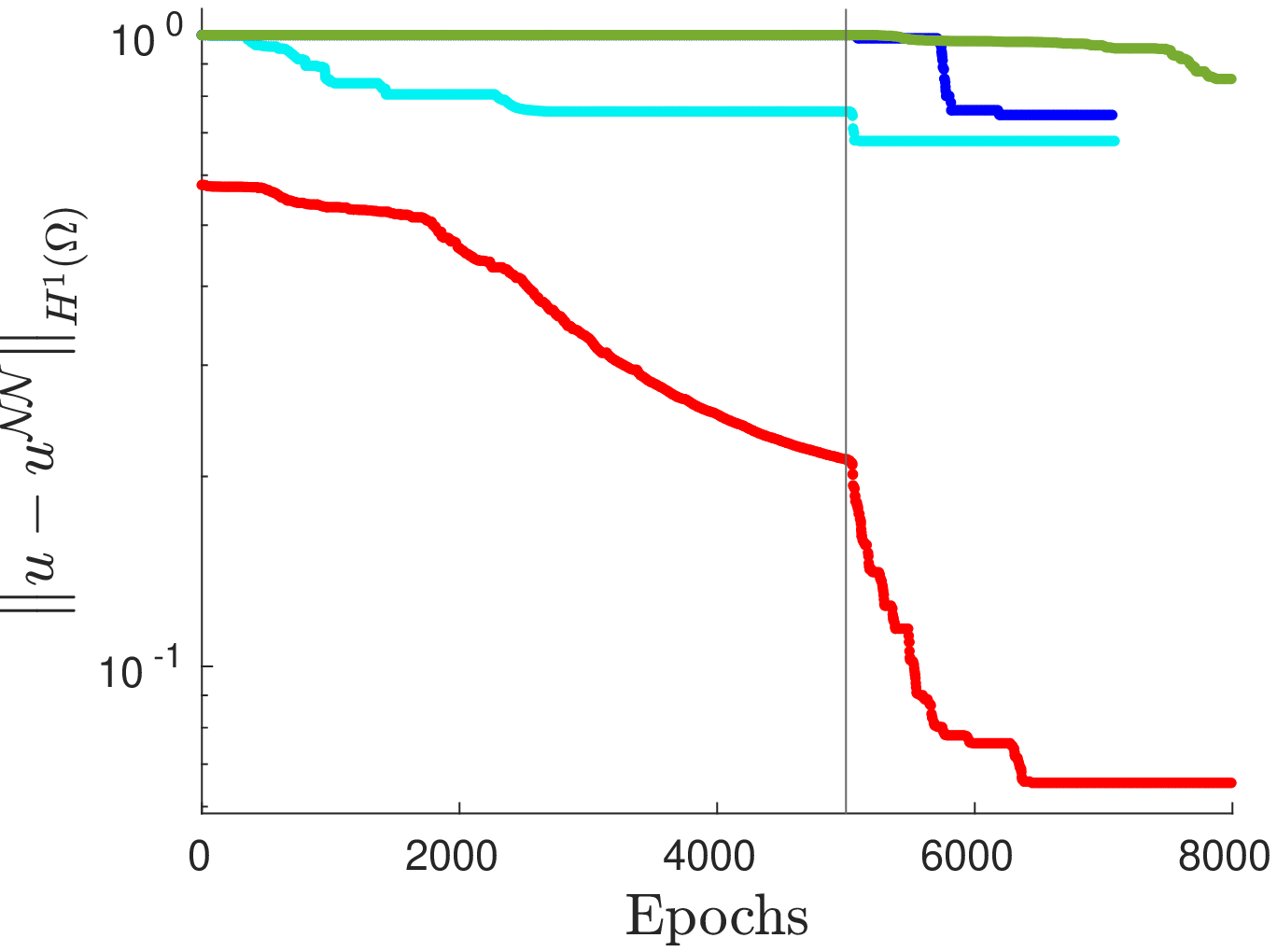} 
  \caption{}
  \label{fig:transport_pinn_errors_rndpts}
\end{subfigure}
  \caption{$H^1$ error decay during the neural network training when solving problem \eqref{eq:transport-model-pb}. (a) VPINN error: $H^1$ error of the most accurate solution is 0.077; (b) PINN error: model is trained with collocation points distributed on a Delaunay mesh
  and the $H^1$ error of the most accurate solution is 0.125; and (c) PINN error: model is trained with 
  collocation points from a uniform distribution and the $H^1$ error of the most accurate solution is 0.051. The legend in (a) also applies to (b) and (c).}
  \label{fig:transport_errors}
\end{figure}

The numerical results 
obtained using the different approaches are 
presented in Fig.~\ref{fig:transport_errors}. In Fig.~{\ref{fig:transport_vpinn_errors}}, problem \eqref{eq:transport-model-pb} is solved with the VPINN method. In this case, $\mathbf{M_A}$ is slightly more accurate and efficient than $\mathbf{M_B}$ (or $\mathbf{M_C}$ since they coincide) if $\lambda$ is chosen properly. However, when the value of $\lambda$ is not optimal, the solution is significantly less accurate. Once more, method $\mathbf{M_D}$ is not competitive with the other approaches. On the other hand, when PINN is considered, exactly imposing the boundary conditions ensures better accuracy and efficiency than using $\mathbf{M_A}$, regardless of the value of $\lambda$ (see Figures  {\ref{fig:transport_pinn_errors_intppts}} and {\ref{fig:transport_pinn_errors_rndpts}}).

\section{Conclusions}\label{sec:conclusion}
In this paper, we analyzed the formulation and the performance of four different approaches to enforce Dirichlet boundary conditions in PINNs and VPINNs on arbitrary polygonal domains. In the first approach, which is the most commonly used when training PINNs, the boundary conditions are imposed by means of additional terms in the loss function that penalize the discrepancy between the neural network output and the prescribed boundary conditions. The subsequent two approaches exactly enforce the boundary conditions and differ in the way they modify the model output in order to force it to satisfy the desired conditions. The last approach, which can be used only when the loss function is derived from the weak formulation of the PDE, is based on Nitsche's method and enforces the boundary conditions variationally.

We have shown that $\mathbf{M_B}$ and $\mathbf{M_D}$, in the considered second-order elliptic PDEs, always ensure the theoretically predicted convergence rate with respect to mesh refinement, regardless of the value of the involved parameter. Instead, method $\mathbf{M_A}$ and $\mathbf{M_C}$ ensure it only if the exact solution is not characterized by an intense oscillatory behaviour.

%
In general, we observed that the most efficient and accurate approach is the one introduced in \cite{sukumar2022exact} (method $\mathbf{M_B}$), which is based on the use of a class of approximate distance functions. A variant of this approach (method $\mathbf{M_C}$) 
leads to suboptimal results and may even ruin the convergence of the method (as in Fig. \ref{fig:eik_pinn_errors_rndpts}). Imposing the boundary conditions via additional cost (method $\mathbf{M_A}$) can be considered a valid alternative, but the choice of the additional penalization parameter is crucial because wrong values can  prevent convergence to the correct solution or dramatically slow down the training. In the proposed numerical experiments we fixed the penalization parameter. As discussed in the Introduction, we highlight that it is possible to tune it during training, but we chose to fix it in order to compare non-intrusive methods with simple implementations. Finally, we observed that Nitsche's method (method $\mathbf{M_D}$) is in some cases similar to $\mathbf{M_A}$ with an acceptable value of $\lambda$, while in other cases requires a second-order optimizer to converge to the correct solution. 

Among possible extensions of this work, we mention applications to high-dimensional PDEs over 
complex geometries, where we expect methods $\mathbf{M_B}$ and $\mathbf{M_C}$ to be even more efficient than their alternatives. In fact, such methods can enforce the correct conditions on each portion of the boundary, whereas methods $\mathbf{M_A}$ and $\mathbf{M_D}$ 
are likely to be less robust and efficient. 

\section*{Acknowledgements}
CC, SB and MP performed this research in the framework of the Italian MIUR Award ``Dipartimenti di Eccellenza 2018-2022" granted to the Department of Mathematical Sciences, Politecnico di Torino (CUP: E11G18000350001). The research leading to this paper has also been partially supported by the SmartData@PoliTO center for Big Data and Machine Learning technologies.
SB  was supported by the Italian MIUR PRIN Project 201744KLJL-004, CC was supported by the Italian MIUR PRIN Project 201752HKH8-003. 
CC, SB and MP are members of the Italian INdAM-GNCS research group.

\bibliographystyle{siam}
\bibliography{bibliography}

\appendix
\renewcommand{\thesection}{\Alph{section}.\arabic{section}}

\section{Appendix: On the Laplacian of the approximate distance function}\label{sec:appendix}
In \cite{sukumar2022exact}, the issue of the blowing-up of the Laplacian of $\phi$ in \eqref{eq:adf_gammaD} is discussed. Herein, we illustrate the same for the simple setting in which $\Gamma_D$ is composed of two edges that intersect. 

Let $\Omega$ be the non-negative quadrant $\{(x,y):x\ge0,y\ge0\}$ and let us consider the (semi-infinite) segments $s_1=\{(x,0):x\ge0\}$ and $s_2=\{(0,y):y\ge0\}$. For the sake of simplicity, we consider the exact distance functions $\phi_1(x,y)=x$ and $\phi_2(x,y)=y$. Let us compute the ADF of order $m=1$ to the boundary $\Gamma_D=s_1\cup s_2$. Substituting $m=1$ in \eqref{eq:adf_gammaD}, $\phi$ can be written as:
\begin{equation*}\label{eq:adf_xy}
\phi = \dfrac{\phi_1\phi_2}{\phi_1+\phi_2}=\dfrac{xy}{x+y}.
\end{equation*}
The gradient of $\phi$ is:
\begin{equation*}\label{eq:gradient_phi}
\nabla \phi = \left[\frac{\partial\phi}{\partial x}, \frac{\partial \phi}{\partial y}\right]^T 
= \left[ \dfrac{(x+y)y-xy}{(x+y)^2}, \dfrac{(x+y)x-xy}{(x+y)^2} \right]^T
=\left[ \dfrac{y^2}{(x+y)^2}, \dfrac{x^2}{(x+y)^2} \right]^T.
\end{equation*}
Note that the gradient is bounded on $\Gamma_D$, and in particular:
\[
\left.{\dfrac{\partial\phi}{\partial x}}\right|_{x=0} = 1,\hspace{1cm}\left.{\dfrac{\partial\phi}{\partial y}}\right|_{y=0} = 1, \hspace{1cm}
\nabla\phi|_{y=\alpha x}=\left[ \dfrac{\alpha^2}{(1+\alpha)^2}, \dfrac{1}{(1+\alpha)^2} \right]^T,
\]
where $\alpha$ is a strictly positive constant, i.e. $\phi$ is an ADF of order 1 and $\nabla \phi$ is bounded along any straight line intersecting the origin and entering inside the domain. The Laplacian of $\phi$ is:
\begin{equation*}\label{eq:laplacian_phi}
\Delta\phi = \frac{\partial^2\phi}{\partial x^2} +\frac{\partial^2\phi}{\partial y^2} 
=\frac{\partial}{\partial x}\left[\frac{y^2}{(x+y)^2}\right] + \frac{\partial}{\partial y}\left[\frac{x^2}{(x+y)^2}\right]
= -2\frac{x^2+y^2}{(x+y)^3}.
\end{equation*}
Consider the limit at $(0,0)$ along the line $y=\alpha x$, for $\alpha\ge0$:
\[
\lim_{x,y\rightarrow0} \Delta\phi = \lim_{x\rightarrow0}  -2\frac{x^2+(\alpha x)^2}{(x+(\alpha x))^3}
= \lim_{x\rightarrow0} -2\frac{(1+\alpha^2)x^2}{(1+\alpha)^3x^3} = \lim_{x\rightarrow0}-2\frac{(1+\alpha^2)}{(1+\alpha)^3}\dfrac1x
=-\infty.
\]
Therefore, $\Delta\phi$ is unbounded at the origin. 

For $m=2$, the function $\phi$ is:
\begin{equation*}\label{eq:adf_xy}
\phi=\dfrac{1}{\sqrt{\frac{1}{\phi_1^2} + \frac{1}{\phi_2^2} }} = \dfrac{\phi_1\phi_2}{\sqrt{\phi_1^2+\phi_2^2}}=\dfrac{xy}{\sqrt{x^2+y^2}}.
\end{equation*}
Its gradient is:
\begin{equation*}\label{eq:gradient_phi}
\begin{split}
\nabla \phi &= \left[\frac{y}{\sqrt{x^2+y^2}} - \frac{x^2y}{(x^2+y^2)^{3/2}}, \frac{x}{\sqrt{x^2+y^2}} - \frac{xy^2}{(x^2+y^2)^{3/2}}\right]^T ,
\end{split}
\end{equation*}
which in polar coordinates ($x=\rho\cos(\theta)$, $y=\rho\sin(\theta)$) is expressed as:
\begin{equation*}\label{eq:gradient_phi}
\begin{split}
\nabla \phi &= \left[\frac{\rho\sin(\theta)}{\rho} - \frac{\rho^3\cos^2(\theta)\sin(\theta)}{\rho^3},
\frac{\rho\cos(\theta)}{\rho} - \frac{\rho^3\cos(\theta)\sin^2(\theta)}{\rho^3}\right]^T ,\\
&=\left[\sin(\theta)-\cos^2(\theta)\sin(\theta), \cos(\theta)-\cos(\theta)\sin^2(\theta)\right]^T.
\end{split}
\end{equation*}
In polar coordinates, we can write
\begin{equation*}\label{eq:gradient_phi}
\begin{split}
\frac{\partial^2\phi}{\partial x^2} &= 3\frac{x^3y}{(x^2+y^2)^{5/2}} - 3\frac{xy}{(x^2+y^2)^{3/2}}
=3\frac{\rho^4\cos^3(\theta)\sin(\theta)}{\rho^5} - 3\frac{\rho^2\cos(\theta)\sin(\theta)}{\rho^3}\\
&=3\frac{\cos(\theta)\sin(\theta)}{\rho}\left[\cos^2(\theta)-1\right]
=-3\frac{\cos(\theta)\sin^3(\theta)}{\rho}.
\end{split}
\end{equation*}
Similarly, 
\[
\frac{\partial^2\phi}{\partial y^2} = -3\frac{\cos^3(\theta)\sin(\theta)}{\rho}
\]
holds, which implies:
\[
\Delta\phi = \frac{\partial^2\phi}{\partial x^2} + \frac{\partial^2\phi}{\partial y^2} = -3\frac{\cos(\theta)\sin(\theta)}{\rho}.
\]
Thus, as in the case $m=1$, $\Delta\phi\rightarrow-\infty$ when $\rho\rightarrow0$. 

We point out that for $\theta=\pi/2$ and $\theta=0$, respectively, we note that:
\[
\left.{\dfrac{\partial\phi}{\partial x}}\right|_{x=0} = 1,\hspace{1cm}\left.{\dfrac{\partial\phi}{\partial y}}\right|_{y=0} = 1,\hspace{1cm}\left.{\dfrac{\partial^2\phi}{\partial x^2}}\right|_{x=0,y>0} = 0,\hspace{1cm}\left.{\dfrac{\partial^2\phi}{\partial y^2}}\right|_{y=0,x>0} = 0.
\]
Therefore, $\phi$ is an ADF that is normalized up to order 2.

\end{document}